\documentclass[ps,dvips,noinfoline]{imsart}
\RequirePackage{hyperref}

\doi{10.1214/06-PS091}
\pubyear{2007}
\volume{4}
\firstpage{193}
\lastpage{267}

\startlocaldefs
\usepackage{amssymb,epsfig} 
\newtheorem{lemma}{Lemma}[section]
\newtheorem{definition}{Definition}[section]
\newtheorem{theorem}{Theorem}[section]
\newtheorem{coro}{Corollary}[section]

\newtheorem{remark}{Remark}[section]

\newcommand{\bsq}{\vrule height .9ex width .8ex depth -.1ex}
\newcommand{\RR}{{\mathbb R}}

\newcommand{\sC}{{\cal C}}
\newcommand{\sB}{{\cal B}}

\newcommand{\sF}{{\cal F}}
\newcommand{\sH}{{\cal H}}

\newcommand{\hsp}{\hspace*{\parindent}}
\newcommand{\eeq}{\end{equation}}

\newcommand{\ra}{\rightarrow}
\newcommand{\deq}{\stackrel{\rm d}{=}}
\newcommand{\beql}[1]{\begin{equation}\label{#1}}
\newcommand{\eqn}[1]{(\ref{#1})}
\newcommand{\beq}{\begin{displaymath}}
\newcommand{\eeqno}{\end{displaymath}}

\newcommand{\qandq}{\quad\mbox{and}\quad}
\newcommand{\qforq}{\quad\mbox{for}\quad}

\newcommand{\qasq}{\quad\mbox{as}\quad}
\newcommand{\qinq}{\quad\mbox{in}\quad}

\newcommand{\qforallq}{\quad\mbox{for all}\quad}

\endlocaldefs

\begin{document}

\begin{frontmatter}

\title{Martingale proofs of many-server heavy-traffic limits for Markovian queues\protect\thanksref{T1}}
\thankstext{T1}{This is an original survey paper}
\runtitle{Martingale Proofs}


\begin{aug}
\author{\fnms{Guodong} \snm{Pang}\ead[label=e1]{gp2224@columbia.edu}\ead[label=u1,url]{www.columbia.edu/$\sim$gp2224}},
\author{\fnms{Rishi} \snm{Talreja}\ead[label=e2]{rt2146@columbia.edu}\ead[label=u2,url]{www.columbia.edu/$\sim$rt2146}}
and
\author{\fnms{Ward} \snm{Whitt}\corref{}\ead[label=e3]{ww2040@columbia.edu}  
\ead[label=u3,url]{www.columbia.edu/$\sim$ww2040}}
\address{Department of Industrial Engineering and Operations Research, Columbia University,\\  New York, NY 10027-6699;
e-mail: \printead[\{gp2224, rt2146, ww2040\}@columbia.edu]{e3}}


\runauthor{G. Pang, R. Talreja and W. Whitt}
\end{aug}

\begin{abstract}
This is an expository review paper illustrating the ``martingale method''
for proving many-server heavy-traffic stochastic-process limits
for queueing models, supporting diffusion-process approximations.
Careful treatment is given to an elementary model -- the
classical infinite-server model
$M/M/\infty$, but models with finitely many servers
and customer abandonment are also treated.
 The Markovian stochastic
process representing the number of customers in the system is
constructed in terms of rate-$1$ Poisson processes in two ways:
(i) through random time changes and (ii) through random thinnings.
Associated martingale representations are obtained for these
constructions by applying, respectively:  (i) optional stopping
theorems where the random time changes are the stopping times and
(ii) the integration theorem associated with random thinning of a
counting process.  Convergence to the diffusion process limit for
the appropriate sequence of scaled queueing processes is obtained
by applying the continuous mapping theorem. A key FCLT and a key FWLLN in
this framework are established both with and without
 applying martingales.
\end{abstract}



 \received{\smonth{12} \syear{2006}}

\tableofcontents

\end{frontmatter}

\section{Introduction}\label{secIntro}
\hsp
The purpose of this paper is to illustrate how to do
martingale proofs of many-server heavy-traffic limit theorems for
queueing models, as in Krichagina and Puhalskii \cite{KP97} and
Puhalskii and Reiman \cite{PR00}. Even though the method is remarkably
effective, it is somewhat complicated.  Thus it is helpful to see
how the argument applies to simple models before considering more
elaborate models. For the more elementary Markovian models
considered here, these many-server heavy-traffic limits
 were originally established by other methods, but new methods of proof have been developed.
For the more elementary models, we show how key steps in the proof - a FCLT (\eqn{c1}) and a FWLLN (Lemma \ref{lemFluid2}) - can
be done without martingales as well as with martingales.  For the argument without martingales,
we follow Mandelbaum and Pats \cite{MP95, MP98}, but without their level of generality and without discussing
strong approximations.

Even though we focus on elementary Markov models here,
we are motivated by the desire to treat more complicated models, such as the non-Markovian $GI/GI/n$ models
in Krichagina and Puhalskii \cite{KP97}, Puhalskii and Reiman \cite{PR00}, Reed \cite{R05} and Kaspi and Ramanan \cite{KR07},
 and network generalizations of these,
such as non-Markovian generalizations of the network models considered by Dai and Tezcan \cite{DT05, DT06, DT07},
Gurvich and Whitt \cite{GW07a, GW07b, GW07c} and
references cited there.  Thus we want to do more than achieve
a quick proof for the simple models, which need not rely so much on martingales; we want to illustrate
martingale methods that may prove useful for more complicated models.

\subsection{The Classical Markovian Infinite-Server Model}\label{secClassic}
\hsp
We start by focusing on what is a candidate to be the
easiest queueing model -- the classic $M/M/\infty$ model -- which has a Poisson arrival process (the first $M$),
i.i.d. exponential service times (the second $M$), independent of the arrival process,
and infinitely many servers.
Let the arrival rate be $\lambda$ and let the individual mean service time be $1/\mu$.

Afterwards, in \S \ref{secErlangA}, we treat the associated $M/M/n/\infty +M$ (Erlang $A$ or Palm) model, which has $n$ servers,
unlimited waiting room, the first-come first-served (FCFS) service discipline and Markovian customer abandonment (the final $+M$);
customers abandon (leave) if they have to spend too long waiting in queue.
The successive customer times to abandon are assumed to be i.i.d. exponential random variables,
independent of the history up to that time.  That assumption is often reasonable with invisible queues,
as in telephone call centers.  Limits for the associated $M/M/n/\infty$ (Erlang $C$) model (without customer abandonment)
are obtained as a corollary, by simply letting the abandonment rate be zero.

For the initial $M/M/\infty$ model,
let $Q(t)$ be the number of customers in the system at time $t$, which coincides
with the number of busy servers at time $t$.
It is well known that $Q \equiv \{Q(t): t \ge 0\}$ is
a birth-and-death stochastic process and that $Q(t) \Rightarrow Q(\infty)$ as $t \ra \infty$,
where
$\Rightarrow$ denotes convergence in distribution,
provided that $P(Q(0) < \infty) = 1$, and $Q(\infty)$ has a Poisson distribution with mean
$E[Q(\infty)] = \lambda/\mu$.

We are interested in heavy-traffic limits in which the arrival rate is allowed to increase.
Accordingly, we consider a sequence of models indexed by $n$ and let the arrival rate in model $n$ be
\beql{a1}
\lambda_n \equiv n \mu, \quad n \ge 1 ~.
\eeq
Let $Q_n (t)$ be the number of customers in the system at time $t$ in model $n$.
  By the observation above,
$Q_n (\infty)$ has a Poisson distribution with mean
$E[Q_n (\infty)] = \lambda_n/\mu = n$.  Since the Poisson distribution approaches the normal distribution
as the mean increases, we know that
\beql{PoissonCLT}
\frac{Q_n (\infty) - n}{\sqrt{n}} \Rightarrow N(0,1) \qasq n \ra \infty ~,
\eeq
where $N(0,1)$ denotes
a standard normal random variables (with mean $0$
and variance $1$).

However, we want to establish a limit for the entire stochastic process $\{Q_n (t): t \ge 0\}$ as $n \ra \infty$.
For that purpose, we consider the scaled processes $X_n \equiv \{X_n (t): t \ge 0\}$ defined by
\beql{a2}
X_n (t) \equiv \frac{Q_n (t) - n}{\sqrt{n}}, \quad t \ge 0 ~.
\eeq

To establish the stochastic-process limit in \eqn{a2}, we have to be careful about the
initial conditions.  We will assume that $X_n (0) \Rightarrow X(0)$ as $n \ra \infty$.
In addition, we assume that the random initial number of busy servers, $Q_n (0)$,
is independent of the arrival process and the service times.  Since the service-time distribution is
exponential, the remaining service times of those customers initially in service
have independent exponential distributions because of
the lack-of-memory property of the exponential distribution.

The heavy-traffic limit theorem asserts that the sequence of stochastic processes
$\{X_n: n \ge 1\}$ converges in distribution in the function space $D \equiv D([0,\infty), \RR)$ to
the Ornstein-Uhlenbeck (OU) diffusion process
as $n \ra \infty$, provided that appropriate initial conditions are in force; see Billingsley \cite{B68}
and Whitt \cite{W02}
for background on such stochastic-process limits.

Here is the theorem for the basic $M/M/\infty$ model:

\begin{theorem} {\em $($heavy-traffic limit in $D$ for the $M/M/\infty$ model$)$}\label{th1}
Consider the sequence of $M/M/\infty$ models defined above.
If
\beql{a2a}
X_n (0) \Rightarrow X(0) \qinq \RR \qasq n \ra \infty ~,
\eeq
then
\beq
X_n \Rightarrow X \qinq D \qasq n \ra \infty ~,
\eeqno
where $X$ is the OU diffusion process with infinitesimal mean $m(x) = -\mu x$
and infinitesimal variance $\sigma^2 (x) = 2 \mu$.  Alternatively, $X$
satisfies the stochastic integral equation
\beql{a4}
X (t) = X (0) + \sqrt{2 \mu} B (t) - \mu \int_{0}^{t} X (s) \, d s ~, \quad t \ge 0 ~,
\eeq
where $B \equiv \{B(t): t \ge 0\}$ is a standard Brownian motion.
Equivalently, $X$ satisfies the stochastic differential equation $(SDE)$
\beql{a5}
d X (t) = - \mu X(t) d t + \sqrt{2 \mu} d B (t), \quad t \ge 0 ~.
\eeq
\end{theorem}

Much of this paper is devoted to {\em four proofs} of Theorem \ref{th1}.
It is possible to base the proof on the martingale functional central limit theorem (FCLT), as given in \S 7.1 of
Ethier and Kurtz \cite{EK86}
and here in \S \ref{secMartFCLT}, as we will show in \S \ref{secProofSecond}, but it is not necessary to do so.
Instead, it can be based on the classic FCLT for the Poisson process,
which is an easy consequence of Donsker's FCLT for sums of i.i.d. random variables,
and the continuous mapping theorem.  Nevertheless, the martingale structure can still play an important role.
With that approach,
the martingale structure can be used to establish stochastic boundedness
of the scaled queueing processes, which we show implies required fluid limits or functional weak laws of large numbers (FWLLN's)
for random-time-change stochastic processes, needed for an application of the continuous mapping theorem
with the composition map.  Alternatively (the third method of proof), the fluid limit can be
established by combining the same continuous mapping with the strong law of large numbers (SLLN)
for the Poisson process.

It is also necessary to understand the characterization of the limiting diffusion process
via \eqn{a4} and \eqn{a5}.  That is aided by the general theory of stochastic integration,
which can be considered part of the martingale theory \cite{EK86, RW87}.

\subsection{The QED Many-Server Heavy-Traffic Limiting Regime}\label{secQED}
\hsp
We will also establish many-server heavy-traffic limits for Markovian models with finitely many servers,
where the number $n$ of servers goes to infinity along with the arrival rate in the limit.  We will consider the sequence of models in
the quality-and-efficiency (QED) many-server heavy-traffic limiting regime, which is defined by the condition
\beql{QED1}
\frac{n \mu - \lambda_n}{\sqrt{n}} \ra \beta \mu \qasq n \ra \infty ~.
\eeq
This limit
in which the arrival rate and number of servers increase together according to \eqn{QED1} is just the right way so that the probability
of delay converges to a nondegenerate limit (strictly between $0$ and $1$); see Halfin and Whitt \cite{HW81}.

We will also allow finite waiting rooms of size $m_n$, where the waiting rooms grow at rate $\sqrt{n}$ as $n \ra \infty$, i.e.,
so that
\beql{QED2}
\frac{m_n}{\sqrt{n}} \ra \kappa \ge 0 \qasq n \ra \infty ~.
\eeq
With the spatial scaling by $\sqrt{n}$, as in \eqn{a2}, this scaling in \eqn{QED2} is just right
to produce a reflecting upper barrier at $\kappa$ in the limit process.

In addition, we allow Markovian abandonment, with each waiting customer abandoning at rate $\theta$.
We let the individual service rate $\mu$ and individual abandonment rate $\theta$ be fixed, independent of $n$.
These modifications produce a sequence of $M/M/n/m_n + M$ models (with $+M$ indicating abandonment).
The special case of the Erlang $A$ (abandonment), $B$ (blocking or loss) and $C$ (delay) models are obtained, respectively,
by (i) letting $m_n = \infty$,
(ii) letting $m_n = 0$, in which case the $+M$ plays no role, and (iii) letting $m_n = \infty$ and $\theta = 0$.
So all the basic many-server Markovian queueing models are covered.

Here is the corresponding theorem for the $M/M/n/m_n + M$ model.

\begin{theorem} {\em $($heavy-traffic limit in $D$ for the $M/M/n/m_n+M$ model$)$}\label{th3}
Consider the sequence of $M/M/n/m_n+M$ models defined above, with the scaling in {\em \eqn{QED1}} and {\em \eqn{QED2}}.
Let $X_n$ be as defined in {\em \eqn{a2}}.
If $X_n (0) \Rightarrow X(0)$ in $\RR$ as $n \ra \infty$,
then
$X_n \Rightarrow X$ in $D$ as $n \ra \infty$,
where the limit $X$ is the diffusion process with infinitesimal mean $m(x) = -\beta \mu -\mu x$ for $x < 0$
and $m(x) = -\beta \mu - \theta x$ for $x > 0$,
infinitesimal variance $\sigma^2 (x) = 2 \mu$ and reflecting upper barrier at $\kappa$.  Alternatively, the limit process $X$
is the unique $(-\infty, \kappa]$-valued process
satisfying the stochastic integral equation
\beql{ee6}
X (t) = X (0) - \beta \mu t + \sqrt{2 \mu} B (t) -  \int_{0}^{t} \left[\mu(X (s)\wedge 0) + \theta (X(s)\vee 0)\right] \, d s - U (t)
\eeq
for $t \ge 0$, where $B \equiv \{B(t): t \ge 0\}$ is a standard Brownian motion and $U$ is the unique nondecreasing nonnegative process in $D$
such that {\em \eqn{ee6}} holds and
\beql{ee6a}
\int_{0}^{\infty} 1_{\{X (t) < \kappa\}} \, d U(t) = 0 ~.
\eeq
\end{theorem}

\subsection{Literature Review}\label{secLit}
\hsp
A landmark in the application of martingales to queues is the book by Br\'{e}maud \cite{B81};
contributions over the previous decade are described there.
As reflected by the references cited in Krichagina and Puhalskii \cite{KP97}
and Puhalskii and Reiman \cite{PR00}, and as substantiated by Puhalskii in personal communication,
Liptser was instrumental in developing the martingale method to prove limit theorems for queueing models, leading to diffusion approximations;
e.g., see Kogan, Liptser and Smorodinskii \cite{KLS86} and \S 10.4 of Liptser and Shiryaev \cite{LS89}.

  The specific $M/M/\infty$ result in Theorem \ref{th1} was first established by
Iglehart \cite{I65}; see Borovkov \cite{B67, B84}, Whitt \cite{W82}, Glynn and Whitt \cite{GW91}, Massey and Whitt \cite{MW93},
Mandelbaum and Pats \cite{MP95, MP98},
Mandelbaum, Massey and Reiman \cite{MMR98} and
Krichagina and Puhalskii \cite{KP97}
 for further discussion and extensions.  A closely related textbook treatment appears in Chapter 6 of Robert \cite{R03}.
 Iglehart applied a different argument; in particular, he applied
Stone's \cite{S63} theorem, which shows for birth-and-death processes
that it suffices to have the infinitesimal means and infinitesimal
variances converge, plus other regularity conditions; see Garnett
et al. \cite{GMR02} and Whitt \cite{W04} for recent uses of that same
technique. Iglehart directly assumed finitely many servers and let
the number of servers increase rapidly with the arrival rate; the
number of servers increases so rapidly that it is tantamount to
having infinitely many servers. The case of infinitely many
servers can be treated by a minor modification of the same
argument.

The corresponding QED finite-server result in Theorem \ref{th3},
in which the arrival rate and number of servers increase together according to \eqn{QED1}, which is just the right way so that the probability
of delay converges to a nondegenerate limit (strictly between $0$ and $1$), was considered by Halfin and Whitt \cite{HW81}
for the case of infinite waiting room and no abandonment.
There have been many subsequent results; e.g., see
Mandelbaum and Pats \cite{MP95, MP98}, Srikant and Whitt \cite{SW96}, Mandelbaum et al. \cite{MMR98}, Puhalskii and Reiman \cite{PR00},
Garnett et al. \cite{GMR02}, Borst et al. \cite{BMR04}, Jelenkovi\'{c} et al. \cite{JMM04},
Whitt \cite{W05}, Mandelbaum and Zeltyn \cite{MZ05} and Reed \cite{R05} for further discussion and extensions.

 Puhalskii and Reiman \cite{PR00}
apply the martingale argument to establish heavy-traffic limits in the QED regime.
They
establish many-server heavy-traffic limits for the $GI/PH/n/\infty$ model, having renewal arrival process (the $GI$),
phase-type service-time distributions (the $PH$), $n$ servers, unlimited waiting space, and the first-come first-served (FCFS) service discipline.
They focus on the number of customers in each phase of service, which leads to convergence to a multi-dimensional diffusion process.
One of our four proofs is essentially their argument.

Whitt \cite{W05} applied the same martingale argument to treat the $G/M/n/m_n+M$ model, having a general stationary point process for an arrival process
(the $G$), i.i.d. exponential service times, $n$ servers, $m_n$ waiting spaces, the FCFS service discipline and customer
abandonment with i.i.d exponential times to abandon; see Theorem 3.1 in \cite{W05}.  Whitt \cite{W05} is
primarily devoted to a heavy-traffic limit for the $G/H^{*}_2/n/m_n$ model, having a special $H^{*}_2$ service-time distribution
(a mixture of an exponential and a point mass at $0$), by a different argument, but there is a short treatment of the
$G/M/n/m_n+M$ model by the martingale argument, following Puhalskii and Reiman \cite{PR00} quite closely.  The martingale
proof is briefly outlined in \S 5 there.  The extension to general arrival
processes is perhaps the most interesting contribution there, but that generalization can also be achieved in other ways.

For our proof without martingales, we follow Mandelbaum and Pats \cite{MP95, MP98}, but without using
strong approximations as they do.
Subsequent related work has established similar asymptotic results for Markovian service networks
with multiple stations, each with multiple servers, and possibly multiple customer classes; e.g.,
see Armony \cite{A05}, Armony et al. \cite{AGM06, AGM06a}, Atar \cite{Atar05a}, Atar et al. \cite{AMR04a, AMR04b},
Dai and Tezcan \cite{DT05, DT06, DT07}, Gurvich and Whitt \cite{GW07a, GW07b, GW07c}, Harrison and Zeevi \cite{HZ04}, and Tezcan \cite{T06}.
These complex stochastic networks have important applications to telephone call centers; see Gans et al. \cite{Gans03}.

\subsection{Organization}\label{secorg}
\hsp
The rest of this paper is organized as follows: We start in
\S \ref{secSample} by constructing the Markovian stochastic
process $Q$ representing the number of customers
in the system in terms of rate-$1$ Poisson processes.  We do this
in two ways:  (i) through random time changes and (ii) through
random thinnings.  We also give the construction in terms of arrival and service times used by Krichagina and Puhalskii \cite{KP97}
to treat the $G/GI/\infty$ model.

Section \S \ref{secMart} is devoted to martingales.
After reviewing basic martingale notions, we
 construct
martingale representations associated with the different constructions.  We justify the first two representations by
applying, respectively:  (i) optional stopping theorems where the
random time changes are the stopping times and (ii) the
integration theorem associated with random thinning of a counting
process. The two resulting integral representations are very
similar. These integral representations are summarized in Theorems
\ref{thMartRep} and \ref{thMartRep2}. In \S \ref{secMart} we also
present two other martingale representations, one based on constructing
counting processes associated with a continuous-time Markov chain,
exploiting the infinitesimal generator, and the other - for $G/GI/\infty$ - based
on the sequential empirical process (see \eqn{Kn4}).

Section \ref{secMain} is devoted to the main steps of the proof of Theorem \ref{th1}
using the first two martingale representations.
In \S \ref{secCont} we show that the integral representation has a unique solution and
constitutes a continuous function mapping $\RR \times D$ into $D$.  In order to establish
measurability, we establish continuity in the Skorohod \cite{S56} $J_1$ topology as well as the topology of uniform convergence over bounded intervals.
As a consequence of the continuity, it suffices to prove FCLT's for the scaled martingales in these integral representations.
For the first martingale representation,
 the scaled martingales themselves are random time changes of the scaled Poisson process.
In \S \ref{secCTMCproof} we show that a FCLT for the martingales based on this first representation, and thus the scaled queueing processes,
can be obtained by
applying the classical FCLT for the Poisson process
and the continuous mapping theorem with the composition map.

To carry out the continuous-mapping argument above, we also need to
 establish the required fluid limit or functional weak law of large numbers (FWLLN)
for the random time changes, needed in the application of the CMT with the composition map.
In \S \ref{secSLLN}, following Mandelbaum and Pats \cite{MP95, MP98}, we show that this fluid-limit step can be established by
applying the strong law of large numbers (SLLN) for the Poisson process
with the same continuous mapping determined by the integral representation.
If we do not use martingales, then we observe that it is easy to extend the stochastic-process limit
to general arrival processes satisfying a FCLT.

But we can also use martingales to establish the FCLT and the FWLLN.  For the FWLLN,
the martingale structure can be exploited via the Lenglart-Rebolledo inequality
to prove stochastic boundedness, first for the scaled martingales and then for the scaled queueing processes,
 which in turn can be used to
to establish the required FWLLN
for the scaled random time changes - Lemma \ref{lemFluid1}.

Since this martingale proof of the fluid limit relies on stochastic boundedness, which is related to tightness, we present
background on these two important concepts in \S \ref{secTightSB}.
For the proof of Theorem \ref{th1}, we conclude there
that it suffices to show that the predictable quadratic variations of the square-integrable martingales are stochastically bounded in $\RR$
in order to have the sequence of scaled queueing processes $\{X_n\}$ be stochastically bounded in $D$.

We complete the proof of Theorem \ref{th1} in \S \ref{secComplete}.  In
Lemma \ref{lemSBfluid} and \S \ref{secFluidSB} we show
that stochastic boundedness of $\{X_n\}$ in $D$ implies the desired fluid limit or FWLLN
needed for the scaled random-time-change processes needed in an application of the continuous mapping with the composition function.
In \S \ref{secSBcomp} we complete the proof by
showing that the predictable quadratic variation processes of the martingales are indeed stochastically bounded in $\RR$.
In \S \ref{secInitial} we show that it is possible to remove a moment condition imposed on $Q_n (0)$, the initial random number of
customers in the system,
in the martingale representation; in particular, we show that it is not necessary to assume that $E[Q_n (0)] < \infty$.
Finally, in \S \ref{secLimFourth} we state the $G/GI/\infty$ limit in Krichagina and Puhalskii \cite{KP97} and show that the
special case $M/M/\infty$ is consistent with Theorem \ref{th1}.

In \S \ref{secOther} we discuss stochastic-process limits for
other queueing models.  In \S \ref{secErlangA} we present a
martingale proof of the corresponding many-server heavy-traffic
limit for the $M/M/n/\infty +M$ model. Corresponding results hold
for the model without customer abandonment by setting the
abandonment rate to zero. The proof is nearly the same as for the
$M/M/\infty$ model.  The only significant difference is the use of
the optional stopping theorem for multiple stopping times with
respect to multiparameter martingales, as in Kurtz \cite{K80}, \S\S
2.8 and 6.2 of Ethier and Kurtz \cite{EK86} and \S 12 of Mandelbaum and
Pats \cite{MP98}. We discuss the extension to cover finite waiting
rooms in \S \ref{secFinite} and non-Markovian arrival processes in
\S \ref{secGen}.

We state a version of the martingale FCLT from p 339 of Ethier and Kurtz \cite{EK86} in \S \ref{secMartFCLT}.
In a companion paper, Whitt \cite{W07}, we review the proof, elaborating on the proof in Ethier and Kurtz \cite{EK86}
and presenting alternative arguments, primarily based on Jacod and Shiryayev \cite{JS87}.  We present a more extensive review of tightness
there.
In \S \ref{secApp} we show how the martingale FCLT implies both the FCLT for a Poisson process
and the required FCLT for the scaled martingales arising in the second and third martingale representations.

\section{Sample-Path Constructions}\label{secSample}
\hsp
We start by making direct sample-path constructions of the stochastic process $Q$
representing the number of customers in the system
in terms of independent Poisson processes.  We show how this can be done in two different ways.
Afterwards, we present a different construction based on arrival and service times,
which only exploits the fact that the service times are mutually independent and independent of the arrival process
and the initial number of customers in the system (with appropriate common distribution assumptions).

\subsection{Random Time Change of Unit-Rate Poisson Processes}\label{secUnit}
\hsp
We first represent arrivals and departures as random time changes of independent
unit-rate Poisson processes.
For that purpose, let $A \equiv \{A(t): t \ge 0\}$ and $S \equiv \{S(t): t \ge 0\}$
be two independent Poisson processes, each with rate (intensity) $1$.
We use the process $A$ to generate arrivals and the process $S$ to generate service completions,
and thus departures.  Let the initial number of busy servers be $Q(0)$.
We assume that $Q(0)$ is a proper random variable independent of the two Poisson processes $A$ and $S$.

The arrival process is simple.  We obtain the originally-specified arrival process with rate $\lambda$
by simply scaling time in the rate-$1$ Poisson process $A$; i.e., we use $A_{\lambda} (t) \equiv A(\lambda t)$ for $t \ge 0$.
It is elementary to see that the
stochastic process $A_{\lambda}$ is a Poisson process with rate $\lambda$.

The treatment of service completions is more complicated.  Let $D (t)$ be the number of departures (service completions)
in the interval $[0,t]$.  We construct this process in terms of $S$ by setting
\beql{b1}
D(t) \equiv S \left(\mu \int_{0}^{t} Q(s) \, d s \right), \quad t \ge 0 ~,
\eeq
but formula \eqn{b1} is more complicated than it looks.  The complication is that
$Q(t)$ appearing as part of the argument inside $S$ necessarily depends on the history
$\{Q(s): 0 \le s < t\}$, which in turn depends on the history of $S$, in particular,
upon $\{S\left(\mu \int_{0}^{s} Q(u) \, du\right): 0 \le s < t\}$.  Hence formula \eqn{b1} is recursive;
we must show that it is well defined.

Of course, the idea is not so complicated:
Formula \eqn{b1} is a consequence of the fact that the intensity of departures at time $s$
is
$\mu Q(s)$,
where the number $Q(s)$ of busy servers at time $s$ is multiplied by the individual service rate $\mu$.
The function $\mu \int_{0}^{t} Q(s) \, d s$ appearing as an argument inside $S$ serves as a random time
change; see \S II.6 of Br\'{e}maud \cite{B81}, Chapter 6 of Ethier and Kurtz \cite{EK86}
and \S 7.4 of Daley and Vere-Jones \cite{DVJ03}.

By the simple conservation of flow, which expresses the content at time $t$ as initial content plus flow in minus flow out, we
have the basic expression
\begin{eqnarray}\label{b2}
Q(t) & \equiv & Q (0) + A(\lambda t) - D(t), \quad t \ge 0 ~, \nonumber \\
& = & Q (0) + A(\lambda t) - S \left(\mu \int_{0}^{t} Q(s) \, d s \right), \quad t \ge 0 ~.
\end{eqnarray}

\begin{lemma}{\em $($construction$)$}\label{lemConstruct}
The stochastic process $\{Q(t): t \ge 0\}$ is well defined as a random element of
the function space $D$ by formula {\em \eqn{b2}}.  Moreover, it is a birth-and-death
stochastic process with constant birth rate $\lambda_k = \lambda$ and linear state-dependent
death rate $\mu_k = k \mu$.
\end{lemma}

\paragraph{Proof.}
 To construct a bonafide random element of $D$, start by conditioning upon the random variable
$Q(0)$ and the two Poisson processes $A$ and $S$.  Then, with these sample paths specified,
we recursively construct the sample path of the stochastic process $Q \equiv \{Q(t): t \ge 0\}$.
By induction over the sucessive jump times of the process $Q$,
we show that the sample paths of $Q$ are right-continuous piecewise-constant real-valued functions
of $t$.  Since the Poisson processes $A$ and $S$ have only finitely many transitions in any finite interval w.p.1,
the same is necessarily true of the constructed process $Q$.
This sample-path construction follows Theorem 4.1 in Chapter 8 on p. 327 of Ethier and Kurtz \cite{EK86}; the argument also appears
in Theorem 9.2 of Mandelbaum et al. \cite{MP98}.
Finally, we can directly
verify that the stochastic process $Q$ satisfies the differential definition of a birth-and-death stochastic process.
Let $\sF_t$ represent the history of the system up to time $t$.  That is the sigma field generated by $\{Q(s): 0 \le s \le t\}$.
It is then straightforward that the infinitesimal transition rates are as claimed:
\begin{eqnarray}
P(Q(t+h) - Q(t) = +1 | Q(t) = k, \sF_t) & = & \lambda h + o(h)  ~, \nonumber \\
P(Q(t+h) - Q(t) = -1 | Q(t) = k, \sF_t) & = & k\mu h + o(h)  ~, \nonumber \\
P(Q(t+h) - Q(t) = 0 | Q(t) = k, \sF_t) & = & 1 - \lambda h - k\mu h + o(h) ~, \nonumber
\end{eqnarray}
as $h \downarrow 0$ for each $k \ge 0$, where $o(h)$ denotes a function $f$ such that $f(h)/h \ra 0$ as $h \downarrow 0$.
Ethier and Kurtz \cite{EK86} approach this last distributional step by verifying that the uniquely constructed process
is the solution of the local-martingale problem for the generator of the Markov process,
which in our case is the birth-and-death process.~~~\bsq

Some further explanation is perhaps helpful.  Our construction above is consistent with the construction
of the queue-length process as a Markov process, specifically, a birth-and-death process.  There are other possible constructions.
A different construction would be the standard approach in discrete-event simulation, with event clocks.
Upon the arrival of each customer, we might schedule the arrival time of the next customer, by generating an exponential interarrival time.
We might also generate the required service time of the current arrival.  With the simulation approach, we have information
about the future state that we do not have with the Markov-process construction.
The critical distinction between the different constructions involves the information available at each time.
The information available is captured by the filtration, which we discuss with the martingale representations
in the next section.

The random-time-change approach we have used here is natural when applying strong approximations, as
was done by Mandelbaum, Massey and Reiman \cite{MMR98} and Mandelbaum and Pats \cite{MP95, MP98}.
They applied the strong approximation for a Poisson process,
as did Kurtz \cite{K78}.  A different use of strong approximations to establish
heavy-traffic limits for non-Markovian infinite-server models is contained in Glynn and Whitt \cite{GW91}.

\subsection{Random Thinning of Poisson Processes}\label{secThin}
\hsp
We now present an alternative construction, which follows \S II.5 of Br\'{e}maud \cite{B81} and
Puhalskii and Reiman \cite{PR00}.  For this construction, let $A_{\lambda}$ and $S_{\mu,k}$ for $k \ge 1$
be independent Poisson processes with rate $\lambda$ and $\mu$, respectively.
As before, we assume that these Poisson processes are independent of the initial number of busy servers, $Q(0)$.
This will not alter the overall system behavior, because the service-time distribution is exponential.  By the lack-of-memory property,
the remaining service times are distributed as i.i.d. exponential random variables, independent of the elapsed service times at time $0$.

We let all arrivals be generated from
arrival process $A_{\lambda}$; we let service completions from individual servers be generated from the processes $S_{\mu,k}$.
To be specific, let the servers be numbered.  With that numbering, at any time we use $S_{\mu,k}$ to generate service completions
from the busy server with the $k^{\rm th}$ smallest index among all busy servers.  (We do not fix attention on a particular server,
because we want the initial indices to refer to the busy servers at all time.)

Instead of \eqn{b1}, we define the departure process by
\beql{b3}
D(t) \equiv \sum_{k=1}^{\infty} \int_{0}^{t} 1_{\{Q(s-) \ge k\}} \, d S_{\mu,k} (s), \quad t \ge 0 ~,
\eeq
where $1_A$ is the indicator function of the event $A$; i.e., $1_A (\omega) = 1$ of $\omega \in A$
and $1_A (\omega) = 0$ otherwise.
It is important that we use the left limit $Q(s-)$ in the integrand of \eqn{b3}, so that the intensity of any service completion
does not depend on that service completion itself; i.e., the functions $Q(s-)$ and $1_{\{Q(s-) = k\}}$ are left-continuous in $s$ for
each sample point.
Then, instead of \eqn{b2}, we have
\begin{eqnarray}\label{b4}
Q(t) & \equiv & Q (0) + A(\lambda t) - \sum_{k=1}^{\infty} \int_{0}^{t} 1_{\{Q(s-) \ge k\}} \, d S_{\mu,k} (s), \quad t \ge 0 ~.
\end{eqnarray}

With this alternative construction, there is an analog of Lemma \ref{lemConstruct} proved in essentially the same way.
In this setting, it is even more evident that we can construct a sample path
of the stochastic process $\{Q(t): t \ge 0\}$ by first conditioning on a realization of $Q(0)$, $\{A_{\lambda} (t): t \ge 0\}$
and $\{S_{\mu, k} (t): t \ge 0\}$, $k \ge 1$.

Even though this construction is different that the one in \S \ref{secUnit}, it too is consistent
with the Markov-process view.  Consistent with most applications, we know what has happened up to time $t$,
but not future arrival times and service times.

\subsection{Construction from Arrival and Service Times}\label{secTimes}
\hsp
Now, following Krichagina and Puhalskii \cite{KP97},
 we construct the sample paths of the process $Q$
 from the arrival
and service times.  This construction applies to the more
general $G/GI/\infty$ system; we will show how the approach applies
to the special case $M/M/\infty$; see \S\S \ref{secFourth} and \ref{secLimFourth}.
See \cite{KP97} for references to related work.

Let $A(t)$ be the cumulative number of arrivals in the interval $[0,t]$ and let
$\tau_i$ be the time of the $i^{\rm th}$ arrival. Let all the
service times be mutually independent random variables,
independent of the arrival process and the number, $Q(0)$, of
customers in the system at time 0 (before new arrivals).
Let $Q(0)$ be independent of the arrival process.
Let the $Q(0)$ initial
customers have service times from $\{\bar{\eta}_i, i\geq 1\}$ with
cumulative distribution function (cdf) $F_0$. Let the new arrivals
have service times from $\{\eta_i, i\geq 1\}$ with cdf $F$.

Then $D(t)$, the number of customers that leave the system by
time $t$, can be expressed as
\begin{equation}
D(t) = \sum\limits_{i = 1}^{Q(0)}\mathbf{1}(\bar{\eta}_i \leq t)
+\sum\limits_{i = 1}^{A(t)} \mathbf{1}(\tau_i + \eta_i \leq t), \quad t \ge 0 ~.
\end{equation}
and
\begin{eqnarray}\label{rep4}
Q(t) &=& \sum\limits_{i = 1}^{Q(0)}\mathbf{1}(\bar{\eta}_i > t) +
\sum\limits_{i = 1}^{A(t)} \mathbf{1}(\tau_i + \eta_i > t), \quad t \ge 0 ~.
\end{eqnarray}

Since this construction does not require $A$ to be Poisson (or even renewal)
or the cdf's $F$ and $F_0$ to be exponential, this approach applies to the non-Markovian $G/GI/\infty$ model.
However, the stochastic process $Q$ itself is no longer Markov in the general setting.
With Poisson arrivals, we can extend \cite{KP97} to obtain a Markov
process by considering the two-parameter process $\{Q(t,y): 0 \le y \le t, t \ge 0\}$,
where $Q(t,y)$ is the number of customers in the system at time $t$ that have elapsed service times greater than or equal to $y$.
(For simplicity in treating the initial conditions, let the initial customers have elapsed service times $0$ at time $0$.
Since $P(A(0) = 0) = 1$, $Q(t,t)$ is (w.p.1) the number of initial customers still in the system at time $t$.)
Then
\begin{eqnarray}\label{rep4}
Q(t,y) &=& \sum\limits_{i = 1}^{Q(0)}\mathbf{1}(\bar{\eta}_i > t) +
\sum\limits_{i = 1}^{A(t-y)} \mathbf{1}(\tau_i + \eta_i > t), \,  0 \le y \le t, \, t \ge 0.
\end{eqnarray}
With renewal ($GI$) arrivals (and the extra assumption that $P(A(0) = 0) = 1$),
we can obtain a Markov process by also appending the elapsed interarrival time.
Of course, there is an alternative to \eqn{rep4} if we add remaining times instead of elapsed times,
but that information is less likely to be available as time evolves.
Heavy-traffic limits for these two-parameter processes follow from the argument of
\cite{KP97}, but we
leave detailed discussion of these extensions to future work.
For other recent constructions and limits in this spirit, see Kaspi and Ramanan \cite{KR07} and references cited there.

\section{Martingale Representations}\label{secMart}
\hsp
For each of the sample-path constructions in the previous
section, we have associated martingale representations. At a high
level, it is easy to summarize how to treat the first two sample-path
constructions, drawing only on the first two chapters of
Br\'{e}maud \cite{B81}. We represent the random time changes as
stopping times with respect to appropriate filtrations and apply
versions of the optional stopping theorem, which states that
random time changes of martingales are again martingales, under
appropriate regularity conditions; e.g., see Theorem T2 on p. 7 of
 \cite{B81}. The problems we consider tend to produce
multiple random time changes, making it desirable to apply the
optional stopping theorem for multiparameter random time changes,
as in \S\S 2.8 and 6.2 of Ethier and Kurtz \cite{EK86}, but we postpone
discussion of that more involved approach until Section
\ref{secErlangA}.

For the random thinning, we instead apply the integration theorem for martingales
associated with counting processes, as in Theorem T8 on p. 27 of Br\'{e}maud \cite{B81},
 which concludes that
integrals of predictable processes with respect to martingales associated with counting processes produce new martingales,
under regularity conditions.

We also present a third martingale representation, based on martingales
associated with generators of Markov processes.
That approach applies nicely here, because the queueing processes we consider are
birth-and-death processes.  Thus we can also apply martingales associated with Markov
chains, as on pp 5-6, 294 of \cite{B81}.
Finally, we present a fourth martingale representation associated with \S \ref{secTimes}.

\subsection{Martingale Basics}\label{secBasics}
\hsp
In this subsection we present some preliminary material on continuous-time martingales.
There is a large body of literature providing background, including the books by:
Br\'{e}maud \cite{B81}, Ethier and Kurtz \cite{EK86}, Liptser and Shiryayev \cite{LS89}, Jacod and Shiryayev \cite{JS87},
Rogers and Williams \cite{RW87, RW00}, Karatzas and Shreve \cite{KS88} and
Kallenberg \cite{K02}.
The early book by Br\'{e}maud \cite{B81} remains especially useful because of its focus on
stochastic point processes and queueing models.  More recent lecture notes by Kurtz \cite{K01}
and van der Vaart \cite{V06} are very helpful as well.

Stochastic integrals play a prominent role, but we only need relatively simple cases.
Since we will be considering martingales associated with counting processes, we
rely on the elementary theory for finite-variation processes, as in \S IV.3 of Rogers and Williams \cite{RW87}.
In that setting, the stochastic integrals reduce to ordinary Stieltjes integrals
and we exploit integration by parts.

On the other hand, as can be seen from Theorems \ref{th1} and \ref{th3}, the limiting stochastic
processes involve stochastic integrals with respect to Brownian motion, which is
substantially more complicated.  However, that still leaves us well within the classical theory.
We can apply the Ito calculus as in Chapter IV of \cite{RW87},
without having to draw upon the advanced theory in Chapter VI.

For the stochastic-process limits in the martingale setting,
it is natural to turn to Ethier and Kurtz \cite{EK86}, Liptser and Shiryayev \cite{LS89}
and Jacod and Shiryayev \cite{JS87}.  We primarily rely on Theorem 7.1 on p. 339 of
Ethier and Kurtz \cite{EK86}.
The Shiryayev books are thorough;
Liptser and Shiryayev \cite{LS89} focuses on basic martingale theory, while
Jacod and Shiryayev \cite{JS87} focuses on stochastic-process limits.

We will start by imposing regularity conditions.  We assume that all stochastic processes
$X \equiv \{X(t): t \ge 0\}$ under consideration
are measurable maps from an underlying probability space $(\Omega, \sF, P)$ to the function space
$D \equiv D ([0,\infty),\RR)$ endowed with the standard Skorohod $J_1$ topology and the associated
Borel $\sigma$-field (generated by the open subsets), which coincides with the usual
$\sigma$-field generated by the coordinate projection maps; see \S\S 3.3 and 11.5 of \cite{W02}.

Since we will be working with martingales, a prominent role is played by the \textbf{filtration}
(histories, family of $\sigma$ fields)
$\textbf{F} \equiv \{\sF_t: t \ge 0\}$ defined on the underlying probability space $(\Omega, \sF, P)$.
(We have the containments:  $\sF_{t_1} \subseteq \sF_{t_2} \subseteq \sF$ for all $t_1 < t_2$.)
As is customary, see p. 1 of \cite{LS89}, we assume that all filtrations satisfy
the \textbf{usual conditions}:

\vspace{0.1in}
(i) they are \textbf{right-continuous}:
\beq
\sF_t = \bigcap_{u: u > t} \sF_u \qforallq t, \quad 0 \le t < \infty ~, \quad and
\eeqno

(ii) \textbf{complete} ($\sF_0$, and thus $\sF_t$ contains all $P$-null sets of $\sF$).

\vspace{0.1in}
We will assume that any stochastic process $X$ under consideration is \textbf{adapted} to the filtration;
i.e., $X$ is \textbf{F}-adapted, which means that $X(t)$ is $\sF_t$-measurable for each $t$.
These regularity conditions guarantee desirable measurability properties, such as progressive measurability:
The stochastic process $X$ is {\em progressively measurable} with respect to the filtration $\textbf{F}$
if, for each $t \ge 0$ and Borel measurable subset $A$ of $\RR$, the set
$\{(s, \omega): 0 \le s \le t, \omega \in \Omega, X(s, \omega) \in A\}$ belongs to the product
$\sigma$ field $\sB([0,t]) \times \sF_t$, where $\sB([0,t])$ is the usual Borel $\sigma$ field on
the interval $[0,t]$.  See p. 5 of \cite{KS88} and Section 1.1 of \cite{LS89}.
In turn, progressive measurability implies measurability,
 regarding $X(t)$ as a map from the product space
$\Omega \times [0,\infty)$ to $\RR$.

A stochastic process $M \equiv \{M(t): t \ge 0\}$ is a
a \textbf{martingale} (submartingale) with respect to a filtration $\textbf{F} \equiv \{\sF_t: t \ge 0\}$
if $M(t)$ is adapted to $\sF_t$
and $E[M(t)] < \infty$ for each $t \ge 0$, and
\beq
E[M(t+s)|\sF_t] = (\ge) M(t)
\eeqno
with probability $1$ (w.p.1) with respect to the underlying probability measure $P$
for each $t \ge 0$ and $s > 0$.

It is often important to have a stronger property than the finite-moment condition $E[M(t)] < \infty$.
A stochastic process $X \equiv \{X(t): t \ge 0\}$ is \textbf{uniformly integrable (UI)} if
\beq
\lim_{n \ra \infty} \sup_{t \ge 0}{\{E[|X(t)| 1_{\{|X(t)| > n\}}]\}} = 0 ~;
\eeqno
see p. 286 of \cite{B81} and p. 114 of \cite{RW00}.
We remark that UI implies, but is not implied by
\beql{UI2}
\sup_{t \ge 0}{\{E[|X(t)|]\}} < \infty ~.
\eeq
A word of warning is appropriate, because the books are not consistent in their assumptions
about {\em integrability}.  A stochastic process $X \equiv \{X(t): t \ge 0\}$ may be
called integrable if $E[|X(t)|] < \infty$ for all $t$ or if the stronger \eqn{UI2} holds.
This variation occurs with {\em square-integrable}, defined below.
Similarly, the basic objects may be taken to be martingales, as defined as above, or might instead be
UI martingales, as on p. 20 of \cite{LS89}.

The stronger UI property is used in preservation theorems - theorems implying that stopped martingales
and stochastic integrals with respect to martingales - remain martingales.
In order to get this property when it is not at first present, the technique
of \textbf{localizing} is applied.
We localize by introducing
associated stopped processes, where the stopping is done with stopping times.
A nonnegative random variable $\tau$ is an $\textbf{F}$-\textbf{stopping time} if
stopping sometime before $t$ depends only on the history up to time $t$, i.e., if
\beq
\{\tau \le t\} \in \sF_t \qforallq t \ge 0~.
\eeqno

For any class $\sC$ of stochastic processes, we define the associated local class $\sC_{loc}$
as the class of stochastic processes $\{X(t): t \ge 0\}$ for which there exists a sequence of stopping times
$\{\tau_n: n \ge 1\}$ such that $\tau_n \ra \infty$ w.p.1 as $n \ra \infty$ and
the associated stopped processes $\{X(\tau_n \wedge t): t \ge 0\}$ belong to class $\sC$
for each $n$, where $a \wedge b \equiv \min{\{a,b\}}$.
We obtain the class of \textbf{local martingales} when $\sC$ is the class of martingales.
Localizing expands the scope, because if we start with martingales, then such stopped martingales remain martingales,
so that all martingales are local martingales.  Localizing is important, because the stopped processes not only are
martingales but
can be taken to be UI martingales; see p. 21 of \cite{LS89} and \S IV.12 of \cite{RW87}.
The UI property is needed for the preservation theorems.

As is customary, we will also be exploiting \textbf{predictable} stochastic processes,
which we will take to mean having left-continuous sample paths.
See p. 8 of \cite{B81} and \S 1.2 of \cite{JS87} for the
more general definition and additional discussion.
The general idea is that a stochastic process $X \equiv \{X(t): t \ge 0\}$ is
predictable if its value at $t$ is determined by values
at times prior to $t$.
But we can obtain simplification by working with stochastic processes
with sample paths in $D$.
In the setting of $D$, we can have a left-continuous process by either (i) considering the left-limit version
$\{X(t-): t \ge 0\}$ of a stochastic process $X$ in $D$ (with $X(0-) \equiv X(0)$) or
(ii) considering a stochastic
process in $D$ with continuous sample paths.  Once we restrict attention to stochastic processes with sample paths in $D$,
we do not need the more general notion, because the left-continuous version is always well defined.
If we allowed more general sample paths, that would not be the case.

We will be interested in martingales associated with counting processes, adapted to the appropriate filtration.
These counting processes
will be nonnegative submartingale processes.  Thus we will be applying the following
special case of the
Doob-Meyer decomposition.

\begin{theorem}{\em $($Doob-Meyer decomposition for nonnegative submartingales$)$}\label{thDoobMeyer}
If $Y$ is a submartingale with nonnegative sample paths,
$E[Y(t)] < \infty$ for each $t$, and $Y$ is adapted to a filtration $\textbf{F} \equiv \{\sF_t\}$,
then there exists an $\textbf{F}$-predictable process $A$, called the
\textbf{compensator} of $Y$ or the dual-predictable projection,
such that $A$ has nonnegative nondecreasing sample paths, $E[A (t)] < \infty$ for each $t$, and
$M \equiv Y - A$ is an $\textbf{F}$-martingale.  The compensator is unique in the sense that the sample paths of any two versions
must be equal w.p.1.
\end{theorem}

\paragraph{Proof.} See \S 1.4 of Karatzas and Shreve \cite{KS88}.
The DL condition in \cite{KS88} is satisfied because of the assumed nonnegativity; see Definition 4.8 and Problem 4.9 on p. 24.
For a full account, see \S VI.6 of \cite{RW87}.~~~\bsq

\subsection{Quadratic Variation and Covariation Processes}\label{secQuad}
\hsp
A central role in the martingale approach to stochastic-process limits is played by the quadratic-variation
and quadratic-covariation processes,
as can be seen from the martingale FCLT stated here in \S \ref{secMartFCLT}.
That in turn depends on the notion of square-integrability.
We say that a stochastic process $X \equiv \{X(t): t \ge 0\}$
is \textbf{square integrable} if $E[X(t)^2] < \infty$ for each $t \ge 0$.
We thus say that a martingale $M \equiv \{M(t): t \ge 0\}$
(with respect to some filtration) is square integrable if $E[M(t)^2] < \infty$ for each $t \ge 0$.
Again, to expand the scope, we can localize, and focus on the class of locally square integrable
martingales.  Because we can localize to get the square-integrability, the condition is not very restrictive.

If $M$ is a square-integrable martingale, then $M^2 \equiv \{M(t)^2: t \ge 0\}$ is necessarily
a submartingale with nonnegative sample paths, and thus satisfies the conditions of Theorem \ref{thDoobMeyer}.
The \textbf{predictable quadratic variation} (PQV) of a square-integrable martingale $M$,
denoted by $\langle M \rangle \equiv \{\langle M \rangle (t): t \ge 0\}$ the (\textbf{angle-bracket process}), is the
compensator of the submartingale $M^2$; i.e., the stochastic
process $\langle M \rangle$ is the unique
nondecreasing nonnegative predictable process such that $E[\langle M \rangle (t)] < \infty$
for each $t$ and $M^2 - \langle M \rangle$ is a martingale with respect to the reference filtration.
(Again, uniqueness holds to the extent that any two versions have the same sample paths w.p.1.)
Not only does square integrability extend by localizing, but Theorem \ref{thDoobMeyer} has a local version;
see p. 375 of \cite{RW87}.  As a consequence the PQV is well defined for any locally square-integrable martingale.

Given two locally square-integrable martingales $M_1$ and $M_2$, the \textbf{predictable
quadratic covariation} can be defined as
\beq
\langle M_1, M_2 \rangle \equiv \frac{1}{4} \left( \langle M_1 + M_2 \rangle - \langle M_1 - M_2 \rangle \right) ~;
\eeqno
see p. 48 of \cite{LS89}.  It can be characterized
as the unique (up to equality of sample paths w.p.1) nondecreasing nonnegative predictable process such
that $E[\langle M_1, M_2 \rangle (t)] < \infty$ for each $t$ and
$M_1 M_2 - \langle M_1, M_2 \rangle$
is a martingale.

We will also be interested in another quadratic variation of a square-integrable martingale $M$,
the so-called \textbf{optional quadratic variation} $[M]$ (OQV), the \textbf{square-bracket process)}.
The square-bracket process is actually more general than the angle-bracket process, because the square-bracket
process is well defined for any local martingale, as opposed to only all locally square-integrable martingales;
see \S\S IV. 18, 26 and VI. 36, 37 of \cite{RW87}.  The following is Theorem 37.8 on p. 389
of \cite{RW87}.  For a stochastic process $M$, let $\Delta M (t) \equiv M(t) - M(t-)$, the jump at $t$,
for $t \ge 0$.

\begin{theorem}{\em $($optional quadratic covariation for local martingales$)$}\label{thOQV}
Let $M_1$ and $M_2$ be local martingales with $M_1 (0) = M_2 (0) = 0$.  Then there exists a unique
process, denoted by $[M_1, M_2]$, with sample paths of finite variation over bounded intervals and $[M_1, M_2] (0) = 0$
such that
\begin{eqnarray}
&& (i)  M_1 M_2 - [M_1, M_2] \quad \mbox{is a local martingale} \nonumber \\
&& (ii) \Delta [M_1, M_2] = (\Delta M_1)(\Delta M_2) ~.  \nonumber
\end{eqnarray}
\end{theorem}
For one local martingale $M$, the optional quadratic variation is then defined as $[M] \equiv [M, M]$.

Note that, for a locally square-integrable martingale $M$, both $M^2 - \langle M \rangle$ and $M^2 - [M]$ are local martingales, but
$\langle M \rangle$ is predictable while $[M]$ is not.  Indeed, subtracting, we see that $[M] - \langle M \rangle$
is a local martingale, so that $\langle M \rangle$ is the compensator of both $M^2$ and $[M]$.

There is also an alternative definition of the OQV; see Theorem 5.57 in \S5.8 of \cite{V06}.
\begin{theorem}{\em $($alternative definition$)$}\label{thOQValt}
If $M_1$ and $M_2$ are local martingales with $M_1 (0) = M_2 (0) = 0$, then
\beq
[M_1, M_2] (t) \equiv \lim_{n \ra \infty} \sum_{i = 1}^{\infty} (M_1 (t_{n,i}) - M_1 (t_{n,i-1}))(M_2 (t_{n,i}) - M_2 (t_{n,i-1})) ~,
\eeqno
where $t_{n,i} \equiv t \wedge (i 2^{-n})$ and
the mode of convergence for the limit as $n \ra \infty$ is understood to be in probability.
The limit is independent of the way that the time points $t_{n,i}$
are selected within the interval $[0,t]$, provided that $t_{n,i} > t_{n,i-1}$
and that the maximum difference $t_{n,i} - t_{n,i-1}$ for points inside the interval $[0,t]$ goes to $0$
as $n \ra \infty$.
\end{theorem}

Unfortunately, these two quadratic-variation processes $\langle M \rangle$ and $[M ]$
associated with a locally square-integrable martingale $M$, and the more general covariation processes $\langle M_1, M_2 \rangle$ and $[M_1, M_2]$,
are somewhat elusive, since the definitions are indirect; it remains to exhibit these processes.
We will exploit our sample-path constructions in terms of Poisson processes
above to identify appropriate quadratic variation and covariation processes in following
subsections.

Fortunately, however, the story about structure is relatively simple in the two cases of interest to us:  (i)
when the martingale is a compensated counting process, and (ii) when the martingale has continuous sample paths.
The story in the second case is easy to tell:  When $M$ is continuous, $\langle M \rangle = [M]$,
and this (predictable and optional) quadratic variation process itself is continuous; see \S VI.34 of \cite{RW87}.
This case applies to Brownian motion and our limit processes.  For standard Brownian motion,
$\langle M \rangle (t)  = [M] (t) = t$, $t \ge 0$.

\subsection{Counting Processes}\label{secCount}
\hsp
The martingales we consider for our pre-limit processes will be compensated counting processes.
By a counting process (or point process), we mean a stochastic process $N \equiv \{N(t): t \ge 0\}$
with nondecreasing nonnegative-integer-valued sample paths in $D$ and $N(0) = 0$.
We say that $N$ is a unit-jump counting process if all jumps are of size $1$.
We say that $N$ is non-explosive if $N(t) < \infty$ w.p.1 for each $t < \infty$.
Equivalently,  if $\{T_n: n \ge 1\}$ is the associated sequence of points,
where
\beq
N (t) \equiv \max{\{n \ge 0: T_n \le t\}}, \quad t \ge 0 ~,
\eeqno
with $T_0 \equiv 0$, then $N$ is a unit-jump
counting process if $T_{n+1} > T_n$ for all $n \ge 0$;
while $N$ is non-explosive if $T_n \ra \infty$ w.p.1 as $n \ra \infty$; see p. 18 of Br\'{e}maud (1981).

As discussed by Br\'{e}maud \cite{B81}, the compensator of a non-explosive unit-jump counting process
is typically (under regularity conditions!) a stochastic process with sample paths that are absolutely continuous with respect to Lebesgue measure,
so that the compensator $A$ can represented as an integral
\beq
A(t) = \int_{0}^{t} X(s) \, ds, \quad t \ge 0 ~,
\eeqno
where $X \equiv \{X(t): t \ge 0\}$ is adapted to the filtration \textbf{F}.
When the compensator has such an integral representation, the integrand $X$
is called the \textbf{stochastic intensity} of the counting process $N$.

We will apply
the following extension
of Theorem \ref{thDoobMeyer}.
\begin{lemma}{\em $($PQV for unit-jump counting processes$)$}\label{lemInc}
  If $N$ is a non-explosive unit-jump
counting process adapted to $\textbf{F}$ with $E[N(t)] < \infty$ for all $t$,
and if the compensator $A$ of $N$
provided by Theorem {\em \ref{thDoobMeyer}} is continuous, then the martingale $M \equiv N - A$
is a square-integrable martingale with respect to $\textbf{F}$ with quadratic variation processes:
\beq
\langle M \rangle = A \qandq [M] = N ~.
\eeqno
\end{lemma}
We first observe that the conditions of Lemma \ref{lemInc} imply that $N$ is a nonnegative $\textbf{F}$-submartingale,
so that we can apply Theorem \ref{thDoobMeyer}.  We use the extra conditions to get more.  We prove Lemma \ref{lemInc} in the Appendix.

\subsection{First Martingale Representation}\label{secFirst}
\hsp
We start with the first sample-path construction in \S \ref{secUnit}.
As a regularity condition, here we assume that
$E[Q(0)] < \infty$.
We will show how to remove that condition later in \S \ref{secInitial}.
It could also be removed immediately if we chose to localize.

\textbf{Here is a quick summary of how martingales enter the picture:}
The Poisson processes $A (t)$ and $S (t)$ underlying the first representation of
the queueing model in \S \ref{secUnit} as well as the new processes $A(\lambda t)$
and $S \left(\mu \int_{0}^{t} Q(s) \, d s \right)$ there have nondecreasing nonnegative sample
paths.  Consequently, they are submartingales with respect to appropriate filtrations (histories,
families of $\sigma$-fields).
Thus, by subtracting the compensators, we obtain martingales.  Then the martingale $M$ so constructed turns out to be square
integrable, admitting a martingale representation $M^2 - \langle M \rangle$,
where $\langle M \rangle$ is the predictable quadratic variation, which in our context will coincide with the compensator.

In constructing this representation, we want to be careful about
the filtration (histories, family of $\sigma$ fields).  Here we will
want to use the filtration $\textbf{F} \equiv \{\sF_t: t \ge 0\}$ defined by
\beql{b5}
\sF_t \equiv  \sigma \left( Q(0), A(\lambda s), S \left( \mu \int_{0}^{s} Q(u) \, d u \right): 0 \le s \le t \right), \quad t \ge 0 ~,
\eeq
augmented by including all null sets.

The following processes will be proved to be $\textbf{F}$-martingales:
\begin{eqnarray}\label{bb3}
M_1 (t) & \equiv & A (\lambda t) - \lambda t, \nonumber \\
M_2 (t) & \equiv & S \left(\mu \int_{0}^{t} Q(s) \, d s \right) - \mu \int_{0}^{t} Q(s) \, d s, \quad t \ge 0 ~,
\end{eqnarray}
where here $A$ refers to the arrival process.  Hence,
instead of \eqn{b2}, we have the alternate martingale representation
\beql{bb4}
Q(t) = Q(0) + M_1 (t) - M_2 (t) + \lambda t - \mu \int_{0}^{t} Q(s) \, d s ~, \quad t \ge 0 ~.
\eeq

In applications of Theorem \ref{thDoobMeyer} and Lemma
\ref{lemInc}, it remains to show that the conditions are satisfied
and to identify the compensator. The following lemma fills in that
step for a random time change of a rate-$1$ Poisson process by
applying the optional stopping theorem.  At this step, it is
natural, as in \S 12 of Mandelbaum and Pats \cite{MP98}, to apply the
optional stopping theorem for martingales indexed by directed sets
(Theorem 2.8.7 on p. 87 of Ethier and Kurtz \cite{EK86}) associated
with multiparameter random time changes (\S 6.2 on p. 311 of
\cite{EK86}), but here we can use a more elementary approach.
For supporting theory at this point, see Theorem 17.24 and
Proposition 7.9 of Kallenberg \cite{K02} and \S 7.4 of Daley and Vere-Jones \cite{DVJ03}.  We use the mutiparameter
random time change in \S \ref{secErlangA}.

Let $\circ$ be the composition map applied to functions, i.e., $(x \circ y) (t) \equiv x(y(t))$.

\begin{lemma}{\em $($random time change of a rate-$1$ Poisson process$)$}\label{lemPQV}
Suppose that $S$ is a rate-$1$ Poisson process adapted to a filtration $\textbf{F} \equiv \{\sF_t: t \ge 0\}$
and $I \equiv \{I(t): t \ge 0\}$ is a stochastic process with continuous nondecreasing nonnegative
sample paths, where $I(t)$ is an $\textbf{F}$-stopping time for each $t \ge 0$.
In addition, suppose that the following moment conditions hold:
\beql{int1}
E[I(t)] < \infty \qandq E[S(I(t))] < \infty \qforallq t \ge 0.
\eeq
Then $S \circ I \equiv \{S(I(t)): t \ge 0\}$ is a non-explosive unit-jump counting process such that
$M \equiv S \circ I - I \equiv \{S(I(t)) - I(t): t \ge 0\}$ is a square-integrable martingale with respect to the
 filtration $\textbf{F}_I \equiv \{\sF_{I(t)}: t \ge 0\}$,
 having quadratic variation processes
 \beql{int2}
\langle M \rangle (t) = I (t) \qandq [M] (t) = S(I(t)), \quad t \ge 0 ~.
\eeq
\end{lemma}

\paragraph{Proof.}  Since the sample paths of $I$ are continuous, it is evident that $S \circ I$ is a
unit-jump counting process.  By condition \eqn{int1}, it is
non-explosive. In order to apply the optional stopping theorem, we
now localize by letting
\beq
 I^m (t) \equiv I(t) \wedge m,
\quad t \ge 0 ~.
\eeqno
 Since $I(t)$ is an \textbf{F}-stopping time,
$I^m (t)$ is a bounded \textbf{F}-stopping time for each $m \ge
1$. The optional stopping theorem then implies that $M^m  \equiv S
\circ I^m - I^m  \equiv \{S(I^m (t)) - I^m (t): t \ge 0\}$ is an
$\textbf{F}_I$-martingale; e.g., see p. 7 of \cite{B81} or
p. 61 of \cite{EK86}. As a consequence $I^m$ is the
compensator of $S \circ I^m$. Since we have the moment conditions
in \eqn{int1}, we can let $m \uparrow \infty$ and apply the
monotone convergence theorem with conditioning, as on p. 280 of
Br\'{e}maud, to deduce that
 $M \equiv S \circ I - I \equiv \{S(I(t)) - I(t): t \ge 0\}$ is a martingale.
Specifically, given that
\beq
E[S(I^m (t+s)) - I^m (t+s)|\sF_{I(s)}] = S(I^m (s)) - I^m (s) \quad \mbox{w.p.1}
\eeqno
for all $m$,
\beq
E[S(I^m (t+s))|\sF_{I(s)}] \ra E[S(I (t+s))|\sF_{I(s)}] \quad \mbox{w.p.1} \qasq m \ra \infty ~,
\eeqno
\beq
E[I^m (t+s)|\sF_{I(s)}] \ra E[I (t+s)|\sF_{I(s)}] \quad \mbox{w.p.1} \qasq m \ra \infty ~,
\eeqno
$S(I^m (s))) \uparrow S(I (s))$ and $I^m (s) \uparrow I (s)$ as $m \ra \infty$,
we have
\beq
E[S(I (t+s)) - I(t+s)|\sF_{I(s)}] = S(I (s)) - I (s) \quad \mbox{w.p.1}
\eeqno
  Lemma \ref{lemInc} implies
the square-integrability and identifies the quadratic variation processes.~~~\bsq

To apply Lemma \ref{lemPQV} to our $M/M/\infty$ queueing problem, we need to verify the finite-moment conditions in \eqn{int1}.
For that purpose, we use a crude inequality:
\begin{lemma}{\em $($crude inequality$)$}\label{lemCrude}
Given the representation {\em \eqn{b2}},
\beql{b101}
Q(t) \le Q(0) + A(\lambda t), \quad t \ge 0~,
\eeq
so that
\beql{b102}
\int_{0}^{t} Q(s) \, ds \le t(Q(0) + A(\lambda t)), \quad t \ge 0 ~.
\eeq
\end{lemma}

Now we want to show that our processes in \eqn{bb3} actually are martingales with respect to the filtration in \eqn{b5}.
To do so, we will apply Lemma \ref{lemPQV}.
However, to apply Lemma \ref{lemPQV}, we will first alter the filtration.
In order to focus on the service completions in the easiest way, we initially condition on the entire
arrival process, and consider the filtration $\textbf{F}^1 \equiv \{\sF^1_t: t \ge 0\}$ defined by
\beql{filt1}
\sF^1_t \equiv  \sigma \left( Q(0), \{A(u): u \ge 0\}, S \left( s \right): 0 \le s \le t \right), \quad t \ge 0 ~,
\eeq
augmented by including all null sets.
Then, as in the statement of Lemma \ref{lemPQV}, we consider the associated
filtration $\textbf{F}^{1}_I \equiv \{\sF^{1}_{I(t)}: t \ge 0\}$.
Finally, we are able to obtain the desired martingale result with respect to the desired
filtration $\textbf{F}$ in \eqn{b5}.

\begin{lemma}{\em $($verifying the conditions.$)$}\label{lemSatisfy}
Suppose that $E[Q(0)] < \infty$ in the setting of \S {\em \ref{secUnit}}
with $\{Q(t): t \ge 0\}$ defined in {\em \eqn{b2}},
\beql{bb69}
I (t) \equiv \mu \int_{0}^{t} Q(s) \, ds, \quad t \ge 0 ~,
\eeq
and the filtration being $\textbf{F}^{1}$ in {\em \eqn{filt1}}.
Then the conditions of Lemma {\em \ref{lemPQV}} are satisfied, so that $S \circ I - I$ is a square-integrable
$\textbf{F}^{1}_{I}$-martingale with $\textbf{F}^{1}_{I}$-compensator $I$ in {\em \eqn{bb69}}.
As a consequence, $S \circ I - I$ is also a square-integrable
$\textbf{F}$-martingale with $\textbf{F}$-compensator $I$ in {\em \eqn{bb69}} for filtration $\textbf{F}$ in {\em \eqn{b5}}.
\end{lemma}

\paragraph{Proof.}
First, we
can apply the crude inequality in \eqn{b102} to establish the required moment conditions:
Since $E[Q(0)] < \infty$,
\beq
E\left[\mu \int_{0}^{t} Q(s) \, ds \right]  \le \mu t(E[Q(0)] + E[A(\lambda t)])
= \mu t E[Q(0)] + \mu \lambda t^2 < \infty, \quad t \ge 0 ~,
\eeqno
and
\beq
S\left(\mu \int_{0}^{t} Q(s) \, ds \right)  \le S\left(\mu t(Q(0) + A(\lambda t)\right), \quad t \ge 0~,
\eeqno
so that
\begin{eqnarray}
E \left[ S \left(\mu \int_{0}^{t} Q(s) \, ds \right) \right] & \le &
E \left[ S \left(\mu t(Q(0) + A(\lambda t))\right) \right] ,  \nonumber \\
  & = & E\left[ E\left[ S \left(\mu t(Q(0) + A(\lambda t))\right)| Q(0) + A(\lambda t) \right] \right], \quad t \ge 0~, \nonumber \\
  & = & \mu t \left(E[Q(0)] + \lambda t  \right) < \infty,  \quad t \ge 0~. \nonumber
\end{eqnarray}

Then, by virtue of \eqn{b2} and the recursive construction in Lemma \ref{lemConstruct}, for each $t \ge 0$, $I(t)$ in \eqn{bb69}
is a stopping time relative to $\sF^1_x$ for all $x \ge 0$, i.e.,
\beq
\{I(t) \le x\} \in \sF^1_x \qforallq x \ge 0 \qandq t \ge 0 ~.
\eeqno
(This step is a bit tricky:  To know $I(t)$, we need to know $Q(s), 0 \le s < t$,
but, by (2.2), that depends on $I(s), 0 \le s < t$.  Hence, to know whether or not $\{I(t) \le x\}$ holds,
it suffices to know $S(u): 0 \le u \le x$.)
Since $\{S(t) - t: t \ge 0\}$ is a martingale with respect to $\textbf{F}^1$
and the moment conditions in \eqn{int1} are satisfied,
we can apply Lemma \ref{lemPQV} to deduce that
$\{S(I(t)) - I(t): t \ge 0\}$ is a square-integrable martingale with respect to the filtration $\textbf{F}^{1}_I \equiv \{\sF^{1}_{I(t)}: t \ge 0\}$
augmented by including all null sets and that $I$ in \eqn{bb69} is both the compensator and the predictable quadratic variation.
Finally, since
the process representing the arrivals after time $t$, i.e., the stochastic process $\{A(t+s) - A(t): s \ge 0\}$, is independent
of $Q(s)$, $0 \le s \le t$, by virtue of the recursive construction in Lemma \ref{lemConstruct} (and the assumption that $A$ is a Poisson process),
we can replace the filtration $\textbf{F}^{1}_I$ by the smaller filtration $\textbf{F}$ in \eqn{b5}.  That completes
the proof.~~~\bsq

We now introduce corresponding processes associated with the \textbf{sequence of models} indexed by $n$.
We have
\beql{b105}
Q_n (t) = Q_n (0) + M^{*}_{n,1} (t) - M^{*}_{n,2} (t) + \lambda_n t - \mu \int_{0}^{t} Q_n (s) \, d s ~, \quad t \ge 0 ~.
\eeq
where
\begin{eqnarray}\label{b106}
M^{*}_{n,1} (t) & \equiv & A (\lambda_n t) - \lambda_n t, \nonumber \\
M^{*}_{n,2} (t) & \equiv & S \left(\mu \int_{0}^{t} Q_n (s) \, d s \right) - \mu \int_{0}^{t} Q_n (s) \, d s,
\end{eqnarray}
The filtrations change with $n$ in the obvious way

We now introduce the scaling, just as in \eqn{a2}.  Let the scaled martingales be
\beql{b107}
M_{n,1} (t) \equiv \frac{M^{*}_{n,1} (t)}{\sqrt{n}} \qandq
M_{n,2} (t) \equiv \frac{M^{*}_{n,2} (t)}{\sqrt{n}}, \quad t \ge 0 ~.
\eeq
Then, from \eqn{b105}--\eqn{b107}, we get
\begin{eqnarray}\label{b108}
X_n (t) & \equiv & \frac{Q_n (t) - n}{\sqrt{n}} \nonumber \\
& = & \frac{Q_n (0) - n}{\sqrt{n}} + \frac{M^{*}_{n,1} (t)}{\sqrt{n}}
- \frac{M^{*}_{n,2} (t)}{\sqrt{n}}  + \frac{\lambda_n t - n \mu t}{\sqrt{n}} \nonumber \\
& & \quad \quad - \mu \int_{0}^{t} \left(\frac{Q_n (s)- n}{\sqrt{n}}\right) \, d s ~, \nonumber \\
& = & \frac{Q_n (0) - n}{\sqrt{n}} + \frac{M^{*}_{n,1} (t)}{\sqrt{n}}
- \frac{M^{*}_{n,2} (t)}{\sqrt{n}}   - \mu \int_{0}^{t} \left(\frac{Q_n (s)- n}{\sqrt{n}}\right) \, d s ~, \nonumber \\
& = & X_n (0) + M_{n,1} (t) - M_{n,2} (t)   - \mu \int_{0}^{t} X_n (s) \, d s ~, \quad t \ge 0 ~.
\end{eqnarray}

Now we summarize this martingale representation for the scaled processes, depending upon the index $n$.
Here is the implication of the analysis above:

\begin{theorem}{\em $($first martingale representation for the scaled
processes$)$}\label{thMartRep}\break
If $E[Q_n (0)] < \infty$
for each $n \ge 1$, then the scaled processes $X_n$ in {\em \eqn{a2}} have the martingale representation
\beql{b109}
X_n (t) \equiv X_n (0) + M_{n,1} (t) - M_{n,2} (t)   - \mu \int_{0}^{t} X_n (s) \, d s ~, \quad t \ge 0 ~,
\eeq
where $M_{n,i}$ are given in {\em \eqn{b106}} and {\em \eqn{b107}}.
These processes $M_{n,i}$ are square-integrable martingales with
respect to the filtrations $\textbf{F}_n \equiv \{\sF_{n,t}: t \ge 0\}$ defined by
\beq
\sF_{n,t} \equiv \sigma\left(Q_n (0), A(\lambda_n s), S\left(\mu \int_{0}^{s} Q_n (u) \, d u \right): 0 \le s \le t\right), \quad t \ge 0 ~,
\eeqno
augmented by including all null sets.  Their associated
predictable
quadratic variations are
$\langle M_{n,1} \rangle (t)  =  \lambda_n t/n$, $t \ge 0$, and
\beql{b110}
\langle M_{n,2} \rangle (t)  =  \frac{\mu}{n} \int_{0}^{t} Q_n (s) \, d s, \quad t \ge 0 ~,
\eeq
where $E [\langle M_{n,2} \rangle (t)] < \infty$ for all $t \ge 0$ and $n \ge 1$.  The associated optional quadratic variations are
$[M_{n,1}] (t)  =  A (\lambda_{n} t)/n$, $t \ge 0$ and
\beq
[M_{n,2}] (t)   =  \frac{S \left( \mu \int_{0}^{t} Q_n (s) \, d s \right) }{n}, \quad t \ge 0 ~.
\eeqno
\end{theorem}

 Note that $X_n$ appears on both sides of the integral representation \eqn{b109}, but $X_n (t)$ appears on the left,
while $X_n (s)$ for $0 \le s \le t$ appears on the right.  In \S \ref{secCont} we show how to
work with this integral representation.

\subsection{Second Martingale Representation}\label{secSecond}
\hsp
We can also start with the second sample-path construction and obtain another
integral representation of the form \eqn{b109}.

Now we start with the martingales:
\begin{eqnarray}\label{b201}
M_A (t) & \equiv & A_{\lambda} (t) - \lambda t, \nonumber \\
M_{S_{\mu,k}} (t) & \equiv & S_{\mu,k} (t) - \mu t, \nonumber \\
M_S (t) & \equiv & \sum_{k=1}^{\infty} \int_{0}^{t} 1_{\{Q(s-) \ge k\}} \, d S_{\mu,k} (s)
- \sum_{k=1}^{\infty} \int_{0}^{t}  \mu 1_{\{Q(s-) \ge k\}} \, ds, \nonumber \\
& = & \sum_{k=1}^{\infty} \int_{0}^{t} 1_{\{Q(s-) \ge k\}} \, d S_{\mu,k} (s)
- \int_{0}^{t}  \sum_{k=1}^{\infty} \mu  1_{\{Q(s-) \ge k\}} \, ds, \nonumber \\
& = & \sum_{k=1}^{\infty} \int_{0}^{t} 1_{\{Q(s-) \ge k\}} \, d S_{\mu,k} (s)
- \mu \int_{0}^{t}  Q(s-) \, ds,
\end{eqnarray}
so that, instead of \eqn{b2}, we have the alternate representation
\beql{b202}
Q(t) = Q(0) + M_A (t) - M_S (t) + \lambda t - \mu \int_{0}^{t} Q(s-) \, ds, ~, \quad t \ge 0 ~,
\eeq
where $M_A$ and $M_S$ are square-integrable martingales with respect to the filtration
$\textbf{F} \equiv \{\sF_t: t \ge 0\}$ defined by
\beql{b203}
\sF_t \equiv  \sigma \left(Q(0), A_{\lambda} (s), S_{\mu,k} (s), \quad k \ge 1: 0 \le s \le t \right), \quad t \ge 0 ~,
\eeq
again augmented by the null sets.
Notice that this martingale representation is very similar to the martingale representation in
\eqn{bb4}.  The martingales in \eqn{bb4} and \eqn{b202} are different and the filtrations in \eqn{b5} and \eqn{b203} are different,
but the predictable quadratic variations are the same and the form of the integral representation is the same.
Thus, there is an analog of Theorem \ref{thMartRep} in this setting.

We now provide theoretical support for the claims above.  First, we put ourselves in
the setting of Lemma \ref{lemInc}.

\begin{lemma}{\em $($a second integrable counting process with unit jumps$)$}\label{lemYInc}
If $E[Q(0)] < \infty$, then, in addition to being adapted to the filtration $\{\sF_t\}$ in {\em \eqn{b203}}, the stochastic process $Y$ defined by
\beql{x1}
Y(t) \equiv \sum_{k=1}^{\infty} \int_{0}^{t} 1_{\{Q(s-) \ge k\}} \, d S_k (s), \quad t \ge 0,
\eeq
is a unit-jump counting process such that
$E[Y(t)] < \infty$ for all $t \ge 0$.
\end{lemma}

\paragraph{Proof.}  It is immediate that $Y$ is a counting process with unit jumps, but there is some question
about integrability
To establish integrability, we apply the crude inequality \eqn{b101}
to get
\beq
Y(t)  \le  \sum_{k=1}^{\infty} \int_{0}^{t} 1_{\{Q(0) + A (\lambda t) \ge k\}} \, d S_k (s)
      \le  \sum_{k=1}^{\infty}  1_{\{Q(0) + A (\lambda t) \ge k\}} S_k (t) ~,
     \eeqno
     so that
\begin{eqnarray}
     E[Y(t)] & \le & \sum_{k=1}^{\infty}  P(Q(0) + A (\lambda t) \ge k)E[S_k (t)] \nonumber \\
             & \le & \sum_{k=1}^{\infty}  P(Q(0) + A (\lambda t) \ge k)\mu t \nonumber \\
             & \le &  \mu t E[Q(0) + A(\lambda t)] =  \mu t(E[Q(0)] + \lambda t) < \infty ~.\nonumber ~~~\bsq
\end{eqnarray}

Given Lemmas \ref{lemInc} and \ref{lemYInc}, it only remains to identify the compensator of the counting process $Y$,
which we call $\tilde{Y}$ since $A$ is used to refer to the arrival process.
For that purpose,
we can apply the integration theorem, as on p. 10 of Br\'{e}maud \cite{B81}.
But our process $Y$ in \eqn{x1} is actually a sum involving infinitely many Poisson processes,
 so we need to be careful.

\begin{lemma}{\em $($identifying the compensator of $Y)$}\label{lemCompensator}
The compensator of $Y$ in {\em \eqn{x1}} is given by
\beql{x4}
\tilde{Y} (t) \equiv  \sum_{k=1}^{\infty} \int_{0}^{t} 1_{\{Q(s-) \ge k\}} \mu \, d s =  \mu \int_{0}^{t}  Q(s-) \, ds, \quad t \ge 0 ~;
\eeq
i.e., $Y - \tilde{Y}$ is an $\textbf{F}$-martingale for the filtration in {\em \eqn{b203}}.
\end{lemma}

\paragraph{Proof.} As indicated above, we can apply the integration theorem on p. 10 of Br\'{e}maud, but
we have to be careful because $Y$ involves
infinitely many Poisson processes.
Hence we first consider the first $n$ terms in the sum.  With that restriction, since the integrand is an indicator function for each $k$,
we consider the integral of a bounded predictable process
with respect to the martingale $\{\sum_{k=1}^{n}(S_{\mu,k} (t) - \mu t): t \ge 0\}$, which is a martingale of ``integrable bounded
variation,'' as required (and defined) by Br\'{e}maud.  As a consequence, $\{M^{n}_{S} (t): t \ge 0\}$ is an $\textbf{F}$-martingale,
where
\begin{eqnarray}\label{x4a}
M^{n}_{S} (t) & \equiv &
  \sum_{k=1}^{n} \int_{0}^{t} 1_{\{Q(s-) \ge k\}} \, d S_{\mu,k} (s)
- \sum_{k=1}^{n} \int_{0}^{t}  \mu 1_{\{Q(s-) \ge k\}} \, ds  ~.
\end{eqnarray}
However, given that $E[Y(t)] < \infty$, we can apply the monotone convergence theorem to each of the two
terms in \eqn{x4a} in order to
take the limit as $n \ra \infty$ to deduce that $E[\tilde{Y} (t)] < \infty$ and
$M_S$ itself, as defined in \eqn{b201}, is an $\textbf{F}$-martingale, which implies
that the compensator of $Y$ in \eqn{x1} is indeed given by \eqn{x4}.~~~\bsq

By this route we obtain another
integral representation for the scaled processes of exactly the same form as in Theorem \ref{thMartRep}.
As before in \eqn{b105}-\eqn{b108}, we introduce the sequence of models indexed by $n$.
The martingales and filtrations are slightly different, but in the end the predictable quadratic variation processes are
essentially the same.

\begin{theorem}{\em $($second martingale representation for the scaled processes$)$}\label{thMartRep2}
If $E[Q_n (0)] < \infty$ for each $n \ge 1$, then the scaled processes $X_n$ in {\em \eqn{a2}} have the martingale representation
\beql{xb109}
X_n (t) \equiv X_n (0) + M_{n,1} (t) - M_{n,2} (t)   - \mu \int_{0}^{t} X_n (s) \, d s ~, \quad t \ge 0 ~,
\eeq
where $M_{n,i}$ are given in {\em \eqn{b107}}, but instead of
{\em \eqn{b106}}, we have
\begin{eqnarray}
M^{*}_{n,1} (t) & \equiv & A_{\lambda_n} (t) - \lambda_n t, \nonumber \\
M^{*}_{n,2} (t) & \equiv & \sum_{k=1}^{\infty} \int_{0}^{t} 1_{\{Q_n (s-) \ge k\}} \, d S_{\mu,k} (s)
- \mu \int_{0}^{t}  Q_n (s-) \, ds ~. \nonumber
\end{eqnarray}
These processes $M_{n,i}$ are square-integrable martingales with respect to the filtrations
$\textbf{F}_n \equiv \{\sF_{n,t}: t \ge 0\}$ defined by
\beq
\sF_{n,t} \equiv  \sigma \left(Q_n (0), A_{\lambda_n} (s), S_{\mu,k} (s), \quad k \ge 1: 0 \le s \le t \right), \quad t \ge 0 ~,
\eeqno
augmented by including all null sets.  Their associated
predictable
quadratic variations are
$\langle M_{n,1} \rangle (t) =  \lambda_n t/n$, $t \ge 0$, and
\beql{xb110}
\langle M_{n,2} \rangle (t)  =  \frac{\mu}{n} \int_{0}^{t} Q_n (s) \, d s, \quad t \ge 0 ~,
\eeq
where $E[\langle M_{n,2} \rangle (t)] < \infty$ for all $t \ge 0$ and $n \ge 1$
and $\langle M_{n,1} \rangle (t)  =  (\lambda_n t/n) = \mu t$.
The associated optional quadratic variations are
$[ M_{n,1} ] (t)  =  A_{\lambda_n} (t)/n$, $t \ge 0$, and
\beql{xb666}
[ M_{n,2} ] (t)  =  \frac{\sum_{k=1}^{\infty} \int_{0}^{t} 1_{\{Q_n (s-) \ge k\}} \, d S_{\mu,k} (s)}{n}, \quad t \ge 0 ~.
\eeq
\end{theorem}

\subsection{Third Martingale Representation}\label{secThird}
\hsp
We can also obtain a martingale representation for the stochastic process $Q$ by exploiting
the fact that the stochastic process $Q$ is a birth-and-death
process.
We have the basic representation
\beq
Q(t)  = Q (0) + A(t) - D(t), \quad t \ge 0 ~,
\eeqno
where $A$ is the arrival process and $D$ is the departure process, as in \eqn{b2}.
Since $Q$ is a birth-and-death process, we can apply the L\'{e}vy and Dynkin formulas,
as on p. 294 of Br\'{e}maud \cite{B81} to obtain martingales associated with various counting processes associated with
$Q$, including the counting processes $A$ and $D$.
Of course, $A$ is easy, but the Dynkin formula immediately yields the desired martingale for $D$:
\beq
M_D (t) \equiv  D(t) - \mu \int_{0}^{t} Q(s) \, ds, \quad t \ge 0 ~,
\eeqno
where the compensator of $M_D$ is just as in the first two martingale representation, i.e., as in \eqn{bb3}, \eqn{bb69}. \eqn{b201} and \eqn{x4};
see pp 6 and 294 of \cite{B81}.
We thus again obtain the martingale representation of the form \eqn{bb4} and \eqn{b202}.
Here, however, the filtration can be taken to be
\beq
\sF_t \equiv \sigma \left(Q(s): 0 \le s \le t\right), \quad t \ge 0 ~.
\eeqno
The proof of Theorem \ref{th1} is then the same as for the second representation,
which will be by an application of the martingale FCLT in \S \ref{secMartFCLT}.

\subsection{Fourth Martingale Representation}\label{secFourth}
\hsp
In this section, we present the martingale representation for
the construction in terms of arrival and service times in \S \ref{secTimes},
but without any proofs.  Consider a sequence of
$G/GI/\infty$ queues indexed by $n$ and let $Q_n (0)$, $Q_n$, $A_n$,
and $D_n$ be the corresponding quantities in the
$n^{\rm th}$ queueing system, just as defined in \S \ref{secTimes}.
For any cdf
$F$, let the associated complementary cdf (ccdf) be $F^c \equiv 1 - F$.

Given representation \eqn{rep4},
the Krichagina and Puhalskii \cite{KP97} insight is to write the process $Q_n$ as
\begin{eqnarray}\label{QnK4}
Q_n (t) & = & \sum\limits_{i =
1}^{Q_n (0)}(\mathbf{1}(\bar{\eta}_i > t)-F_0^c(t)) + Q_n (0) F_0^c(t) \nonumber \\
&& \quad + n \int_0^t \int_0^{\infty} \mathbf{1}( s+x > t) d K_n \left(\frac{A_n(s)}{n},x\right) ,
\end{eqnarray}
 where
 \beql{Kn4} K_n (t,x) \equiv
\frac{1}{n}\sum\limits_{i = 1}^{\lfloor n t
\rfloor}\mathbf{1}(\eta_i \leq x), \quad t \geq 0, \quad x\geq 0 ~,
\eeq is a sequential empirical process (a random field, having two parameters), so that
 \beql{KnA4}
 K_n \left(\frac{A_n (t)}{n},x\right) = \frac{1}{n}\sum\limits_{i =
1}^{A_n (t)}\mathbf{1}(\eta_i \leq x), \quad t \geq 0, \quad x\geq
0 .\eeq

The division by $n$ in \eqn{Kn4} provides a law-of-large-numbers
(LLN) or fluid-limit scaling. To proceed, we define associated
queueing processes with LLN scaling. In particular, define the
normalized processes $\bar{Q}_n \equiv \{\bar{Q}_n (t), t\geq
0\}$, $\bar{A}_n \equiv \{\bar{A}_n (t), t\geq 0\}$ and $\bar{D}_n
\equiv \{\bar{D}_n (t), t\geq 0\}$ as

\beql{fq4} \bar{Q}_n (t) \equiv \frac{1}{n}Q_n(t), \quad
\bar{A}_n(t) \equiv \frac{1}{n}A_n(t), \quad \bar{D}_n(t) \equiv
\frac{1}{n}D_n(t), \quad t \ge 0 ~. \eeq

For our general arrival process, we assume that $\bar{A}_n(t) \ra
a (t) \equiv \mu t$ w.p.1 as $n \ra \infty$. For the $M/M/\infty$
special case, that follows from \eqn{a1}.

Next write equation \eqn{KnA4} as
\begin{eqnarray*}
K_n\left(\frac{A_n (t)}{n},x\right) &=&
\frac{1}{\sqrt{n}}\Big[\frac{1}{\sqrt{n}}\sum\limits_{i =
1}^{A_n(t)}(\mathbf{1}(\eta_i \leq x)-F(x))\Big] +
\frac{1}{n}A_n(t)F(x) \\
&=& \frac{1}{\sqrt{n}}\Big[\frac{1}{\sqrt{n}}\sum\limits_{i =
1}^{A_n (t)}(\mathbf{1}(\eta_i \leq x)-F(x))\Big]  {} \nonumber \\
& & {} +  \frac{1}{\sqrt{n}} \Big[\sqrt{n}(\bar{A}_n
(t)-a(t))\Big]F(x) + a(t) F(x).
\end{eqnarray*}

Now we introduce stochastic processes with central-limit-theorem (CLT) scaling.
In particular, let
\beq V_n(t,x) \equiv \frac{1}{\sqrt{n}}\sum\limits_{i =
1}^{A_n(t)}(\mathbf{1}(\eta_i \leq x)-F(x)),
\eeqno
and
\beql{fq4a}
\hat{A}_n (t)\equiv \sqrt{n}(\bar{A}_n (t) -a(t)).
\eeq
Then
$$
K_n(\bar{A}_n (t),x) =
 \frac{1}{\sqrt{n}}V_n(t,x) +
\frac{1}{\sqrt{n}}\hat{A}_n(t) F(x) + a(t) F(x),
$$
so that the process $Q_n$ in \eqn{QnK4} can be written as
\begin{eqnarray} \label{QnM4}
Q_n (t) &=& \sum\limits_{i = 1}^{Q_n (0)}(\mathbf{1}(\bar{\eta}_i >
t)-F_0^c(t)) + Q_n (0) F_0^c(t)  \nonumber \\[-2pt]
&&{}+ \sqrt{n}\int_0^t \int_0^{\infty}
\mathbf{1}( s+x > t) d V_n (s,x) \nonumber \\[-2pt]
&&{} + \sqrt{n}\int_0^t \int_0^{\infty} \mathbf{1}( s+x > t)d
\hat{A}_n (s) d F(x) \nonumber \\[-2pt]
&&{} + n \int_0^t \int_0^{\infty} \mathbf{1}( s+x
> t) d a(s) d F(x) \nonumber \\[-2pt]
 &=& \sum\limits_{i = 1}^{Q_n (0)}(\mathbf{1}(\bar{\eta}_i >
t)-F_0^c(t)) + Q_n (0) F_0^c(t)  \\[-2pt]
&&{} + \sqrt{n}\int_0^t \int_0^{\infty}
\mathbf{1}( s+x > t) d V_n(s,x)  {} \nonumber \\[-2pt]
& & {} + \sqrt{n} \int_0^t F^{c}(t-s)d \hat{A}_n(s)  + n
\int_0^t F^{c}(t-s) da(s)\nonumber\\
&=& \sum\limits_{i = 1}^{Q_n (0)}(\mathbf{1}(\bar{\eta}_i >
t)-F_0^c(t)) + Q_n (0) F_0^c(t) + n \int_0^t F^{c}(t-s) da(s) \nonumber \\
&& \quad  +  \sqrt{n}(M_{n,1} (t) - M_{n,2} (t)) ~, \nonumber
\end{eqnarray}
where
\begin{eqnarray}\label{Mn14}
M_{n,1} (t) & \equiv & \int_0^t F^c(t-s)d \hat{A}_n (s)
\end{eqnarray}
and
\begin{eqnarray} \label{Mn24}
M_{n,2} (t) &\equiv & - \int_0^t \int_0^{\infty} \mathbf{1}( s+x
> t) d V_n(s,x) \nonumber \\
&=& \int_0^t \int_0^{\infty} \mathbf{1}( s+x \leq t) d V_n(s,x).
\end{eqnarray}

In contrast to previous representations, note that, except for
$Q_n (0)$ which can be regarded as known, $Q_n (s)$ for $s < t$
does not appear on the righthand side of representation
\eqn{QnM4}.  Instead of the integral representations in Theorems
\ref{thMartRep} and \ref{thMartRep2}, here we have a direct
expression of $Q_n (t)$ in terms of other model elements, but we
will see that some of these model elements in turn do have
integral representations.

By equations \eqn{fq4} and \eqn{QnM4}, we have
\begin{eqnarray}\label{qnM4}
\bar{Q}_n (t) & = & \frac{1}{n}\sum\limits_{i = 1}^{Q_n
(0)}(\mathbf{1}(\bar{\eta}_i >
t)-F_0^c(t)) + \bar{Q}_n (0) F_0^c(t) +  \int_0^t F^c(t-s) da(s) \nonumber \\
&& \quad + \frac{1}{\sqrt{n}}(M_{n,1}(t) - M_{n,2} (t)), \quad t
\ge 0.
\end{eqnarray}

From equation \eqn{qnM4}, we can prove the following FWLLN. We
remark that we could allow more general limit functions $a$ for
the LLN-scaled arrival process.

\begin{theorem}{\em $($FWLLN for the fourth martingale representation$)$}\label{FL4}
If there is convergence $(\bar{Q}_n (0), \bar{A}_n) \Rightarrow (q (0), a)$ in
$\RR \times D$ as $n \rightarrow \infty$, where $a(t) \equiv \mu
t$, $t \ge 0$, then $\bar{Q}_n \Rightarrow q$, where
\begin{equation} \label{q4}
q(t) = q(0) F_0^c(t)+ \int_0^t F^c(t-s) d a(s), \quad t \ge 0~.
\end{equation}
For the $M/M/\infty$ special case,
\beq
 q(t) = q(0) e^{-\mu t}+
\mu\int_0^t e^{-\mu(t-s)} d s  = 1 -(1- q (0)) e^{-\mu t}, \quad t
\ge 0.
\eeqno
If, in addition, $q(0) = 1$, then $q(t) = 1$ for $t \ge 0$.
\end{theorem}

Let the scaled process $X_n$ be defined
by \beql {Xnqnq4} X_n (t) = \sqrt{n} (\bar{Q}_n (t) - q(t)), \quad t
\ge 0. \eeq
If $q(t) = 1$ for $t \ge 0$, then \eqn{Xnqnq4} coincides with \eqn{a2}.
By equations \eqn{Xnqnq4}, \eqn{qnM4} and \eqn{q4}, we obtain the
following theorem for the scaled processes.
\begin{theorem} {\em $($fourth martingale representation for the scaled
processes$)$}\label{MRS4} The scaled process $X_n$ in
{\em \eqn{Xnqnq4}} has the representation
 \begin{eqnarray}\label{XMARS4}
X_n (t) & = & \frac{1}{\sqrt{n}}\sum\limits_{i =
1}^{Q_{n}(0)}(\mathbf{1}(\bar{\eta}_i > t)-F_0^c(t)) +
\sqrt{n}(\bar{Q}_n (0)-q(0)) F_0^c(t) \nonumber \\
&& \quad + M_{n,1}(t) - M_{n,2}(t),
\quad t \ge 0, \end{eqnarray}
 where $M_{n,1}$ and $M_{n,2}$ are defined as in
{\em \eqn{Mn14}} and {\em \eqn{Mn24}}, respectively.
\end{theorem}

The situation is more complicated here, because
the processes $M_{n,1}$ and $M_{n,2}$ in \eqn{Mn14} and \eqn{Mn24} are not naturally martingales
for the $G/GI/\infty$ model or even the $M/M/\infty$ special case, with respect to the obvious
filtration, but they can be
analyzed by martingale methods.
In particular, associated martingales can be exploited to
establish stochastic-process limits.
In particular,
the proof of the FCLT for the
processes $X_n$ in \eqn{Xnqnq4} - see \S \ref{secLimFourth} - exploits
semimartingale decompositions of the following two-parameter
process $U_n \equiv \{U_n (t,x), t\geq 0, 0\leq x \leq 1\}$ (and
related martingale properties):

\beql{Un4} U_n(t,x) \equiv \frac{1}{\sqrt{n}}\sum\limits_{i =
1}^{\lfloor n t\rfloor}(\mathbf{1}(\zeta_i \leq x)-x)~, \eeq where
the $\zeta_i$ are independent and uniformly distributed on
$[0,1]$.

Extending Bickel and Wichura \cite{BW71}, Krichagina and Puhalskii
\cite{KP97} proved that the sequence of processes $\{U_n, n\geq 1\}$
converges in distribution to the Kiefer process $U$ in
$D([0,\infty),D([0,1]))$. For properties of Kiefer processes, we
refer to Cs\"{o}rg\'{o} M. and P. R\'{e}v\'{e}z \cite{CR81} and Khoshnevisan \cite{K02}.  The
importance of the Kiefer process for infinite-server queues was
evidently first observed by Louchard \cite{L88}.

The process $U_n$ has the following semimartingale decomposition
(See Chapter IX of Jacod and Shiryaev \cite{JS87}):
$$
U_n (t,x) = -\int_0^x \frac{U_n (t,y)}{1-y} d y + M_{n,0} (t,x),
$$
where
$$
M_{n,0} (t,x) \equiv \frac{1}{\sqrt{n}}\sum\limits_{i = 1}^{\lfloor n
t\rfloor}\Big(\mathbf{1}(\zeta_i \leq x)-\int_0^{x\wedge \zeta_i}
\frac{1}{1-y} d y\Big),
$$
is a square-integrable martingale relative to the filtration $F^n
= \bigvee_{i \leq \lfloor n t\rfloor} \mathcal{F}^i(x)$ and
$\mathcal{F}^i(x) = \sigma(\mathbf{1}(\zeta_i \leq y), 0 \leq y
\leq x)\vee \mathcal{N}$ for all $x\in[0,1]$.

Hence $V_n (t,x) = U_n (a_n(t),F(x))$ can be written as
$$
V_n(t,x) = -\int_0^x \frac{V_{n, -} (t,y)}{1-F_{-}(y)} d F(y) +
L_n (t,x) ~,
$$
where
$$
L_n (t,x) \equiv \frac{1}{\sqrt{n}}\sum\limits_{i =
1}^{A_n(t)}(\mathbf{1}(\eta_i \leq x)-\int_0^{x\wedge\eta_i}
\frac{1}{1-F_{-}(y)} d F(y),
$$
$F_{-}(y)$ is the left-continuous version of $F$, $F_{-} (0)
\equiv 0$ and $V_{n,-}$ is the left-continuous version of $V_n$ in the
second argument. Therefore, $M_{n,2}$ can be written as
$$
M_{n,2} (t) = G_n (t) + H_n (t)~,
$$
where
\begin{eqnarray*}
G_n(t) &\equiv & \int_0^t \int_0^{\infty} \mathbf{1}( s+x \leq t) d\Big(
-\int_0^x \frac{V_{n, -} (s,y)}{1-F_{-}(y)} d F(y) \Big) \nonumber \\
&=& - \int_0^t \frac{V_{n,-} (t-x,x)}{1-F_{-}(x)} d F(x)
\end{eqnarray*}
and
$$
H_n(t) \equiv \int_0^t \int_0^{\infty}\mathbf{1}( s+x \leq t)d L_n(s,x).
$$

In closing this subsection, we remark that an associated representation holds for the
two-parameter process $Q (t, y)$ in \eqn{rep4}.
Let the associated scaled two-parameter process be defined
by
\beql{Two1}
X_n (t, y) \equiv \sqrt{n} (\bar{Q}_n (t, y) - q(t, y)), \quad t
\ge 0,
\eeq
where $\bar{Q}_n (t, y) \equiv Q_n (t, y)/n$, $\bar{Q}_n \Rightarrow q$ as $n \ra \infty$ and
\beql{Two2}
q(t,y) = q(0) F^{c}_{0} (t) + \int_{0}^{t-y} F^c (t-s) \, d a(s) ~.
\eeq

\begin{coro}{\em $($associated representation for the scaled two-parameter
processes$)$}\label{MRS4cor} Paralleling {\em \eqn{XMARS4}}, the scaled process in {\em \eqn{Two1}} has the representation
\begin{eqnarray}\label{Two3}
X_n (t, y) & = & \frac{1}{\sqrt{n}}\sum\limits_{i =
1}^{Q_{n}(0)}(\mathbf{1}(\bar{\eta}_i > t)-F_0^c(t)) +
\sqrt{n}(\bar{Q}_n (0)-q(0)) F_0^c(t) \nonumber \\
&& \quad + M_{n,1}(t, y) - M_{n,2}(t, y),
\quad t \ge 0,
\end{eqnarray}
where, paralleling {\em \eqn{Mn14}} and {\em \eqn{Mn24}},
\begin{eqnarray}
M_{n,1} (t,y) & = & \int_0^{t-y} F^c(t-s)d \hat{A}_n (s) \nonumber
\end{eqnarray}
and
\begin{eqnarray}
M_{n,2} (t,y) &=& \int_0^{t-y} \int_0^{\infty} \mathbf{1}( s+x \leq t) d V_n(s,x). \nonumber
\end{eqnarray}
\end{coro}

\section[Main Steps in the Proof of Theorem 1.1]{Main Steps in the Proof of Theorem \ref{th1}}\label{secMain}
\hsp
In this section we indicate the main steps in the proof of Theorem \ref{th1} starting from one of
the first three martingale representations in the previous section.
First, in \S \ref{secCont} we show that the integral representation appearing in both Theorems
\ref{thMartRep} and \ref{thMartRep2} has a unique solution, so that it constitutes
a continuous function from $D$ to $D$.
Next, in \S \ref{secCTMCproof} we show how the limit can be obtained from
the functional central limit theorem for the Poisson process and the continuous mapping theorem,
but a fluid limit (Lemma \ref{lemFluid1} or Lemma \ref{lemFluid2}) remains to be verified.
In \S \ref{secSLLN} we show how the proof can be completed without martingales by directly establishing
that associated fluid limit.  In \S\S \ref{secTightSB}-\ref{secComplete} we show how martingales can achieve the same result.
In \S \ref{secLimFourth} we indicate how to complete the proof with the fourth martingale representation.

\subsection{Continuity of the Integral Representation}\label{secCont}
\hsp
We apply the \textbf{continuous-mapping theorem (CMT)} with the integral representation in \eqn{b109} and \eqn{xb109}
in order to establish the desired convergence;
for background on the CMT, see \S 3.4 of \cite{W02}.
In subsequent sections we will show that the scaled martingales converge weakly to independent Brownian motions, i.e.,
\beql{c1}
(M_{n,1}, M_{n,2}) \Rightarrow (\sqrt{\mu} B_1, \sqrt{\mu} B_2) \qinq D^2 \equiv D \times D  \qasq n \ra \infty ~,
\eeq
where $B_1$ and $B_2$ are two independent standard Brownian motions,
from which an application of the CMT with subtraction yields
\beql{c1a}
M_{n,1}- M_{n,2} \Rightarrow \sqrt{\mu} B_1 - \sqrt{\mu} B_2 \deq \sqrt{2 \mu} B \qinq D  \qasq n \ra \infty ~,
\eeq
where $B$ is a single standard Brownian motion.

  We then apply the CMT
with the function $f: D \times \RR \ra D$ taking $(y,b)$ into $x$ determined by the integral representation
\beql{c2}
x(t) =  b + y (t) - \mu \int_{0}^{t} x (s) \, ds ~, \quad t \ge 0 ~.
\eeq
In the pre-limit, the function $y$ in \eqn{c2} is played by $M_{n,1} - M_{n,2} \equiv \{M_{n,1} (t) - M_{n,2} (t): t \ge 0\}$ in \eqn{c1a},
while $b$ is played by $X_n (0)$.
In the limit, the function $y$ in \eqn{c2} is played by the limit $\sqrt{2 \mu} B$ in \eqn{c1a},
while $b$ is played by $X (0)$.  (The constant $b$ does not play an essential role in \eqn{c2};
it is sometimes convenient when we want to focus on the solution $x$ as a function of the initial conditions.)

For our application, the limiting stochastic process in \eqn{c1a} has continuous sample paths.
Moreover, the function $f$ in \eqn{c2} maps continuous functions into continuous functions,
as we show below.  Hence, it suffices to show that the map $f: D \times \RR \ra D$ is measurable
and continuous at continuous limits.  Since the limit is necessarily continuous as well,
the required continuity follows from continuity when the function space $D$ appearing in both the domain
and the range is endowed with the topology of uniform convergence on bounded intervals.
However, if we only establish such continuity, then that leaves open the issue of measurability.
It is significant that the $\sigma$ field on $D$ generated by the topology of uniform convergence on bounded intervals
is not the desired customary $\sigma$ field on $D$, which is generated by the coordinate projections
or by any of the Skorohod topologies; see \S 11.5 of \cite{W02} and \S 18 of Billingsley \cite{B68}.
We prove measurability with respect to the appropriate $\sigma$ field on $D$ (generated by the $J_1$ topology)
by proving continuity when the function space $D$ appearing in both the domain
and the range is endowed with the Skorohod $J_1$
topology.  That implies the required measurability.  At the same time, of course, it
provides continuity in that setting.

We now establish
 the basic continuity result.  We establish a slightly more general form than needed here in order to be able to
 treat other cases.  In particular, we introduce a Lipschitz function $h: \RR \ra \RR$;
i.e., we assume that there exists a constant $c > 0$ such that
\beql{lip2}
| h (s_1) - h(s_2)| \le c |s_1 - s_2| \qforallq s_1, s_2 \in \RR ~.
\eeq
We apply the more general form to treat the Erlang $A$ model in \S \ref{secErlangA}.  Theorem \ref{thCMT2} in \S \ref{secFinite} involves an even
more general version in which $h: D \ra D$.

\begin{theorem}{$($continuity of the integral representation$)$}\label{thCMT}
Consider the integral representation
\beql{lip1}
x(t) =  b + y (t) +  \int_{0}^{t} h(x (s)) \, ds ~, \quad t \ge 0 ~,
\eeq
where $h: \RR \ra \RR$ satisfies $h(0) = 0$ and is a Lipschitz function as defined in {\em \eqn{lip2}}.
The integral representation in {\em \eqn{lip1}} has a unique solution $x$,
so that the integral representation constitutes a function $f: D \times \RR \ra D$ mapping $(y,b)$ into $x \equiv f(y,b)$.
In addition, the function $f$ is continuous provided that the function space $D$ $($in both the domain and range$)$ is endowed with
either: $(i)$ the topology of uniform
convergence over bounded intervals or $(ii)$ the Skorohod $J_1$ topology.
Moreover, if $y$ is continuous, then so is $x$.
\end{theorem}

\paragraph{Proof.}  If $y$ is a piecewise-constant function, then we can directly construct the solution $x$ of the
integral representation by doing an inductive construction, just as in
Lemma \ref{lemConstruct}.  Since any element $y$ of $D$ can be represented as the limit
of piecewise-constant functions, where the convergence is uniform over bounded intervals,
using endpoints that are continuity points of $y$, we can then extend the function $f$
to arbitrary elements of $D$, exploiting continuity in the topology of uniform convergence
over bounded intervals, shown below.  Uniqueness follows from the fact that the only function
$x$ in $D$ satisfying the inequality
\beq
| x(t)| \le  c \int_{0}^{t} |x(s)| \, ds, \quad t \ge 0 ~,
\eeqno
is the zero function, which is a consequence of Gronwall's inequality,
which we re-state in Lemma \ref{lemGronwall} below in the form needed here.

For the remainder of the proof, we
 apply Gronwall's inequality again.
We introduce the norm
\beq
|| x||_T \equiv \sup_{0 \le t \le T}{|x(t)|} ~.
\eeqno
First consider the case of the topology of uniform convergence over bounded intervals.
We need to show that, for any $\epsilon > 0$, there exists a $\delta > 0$
such that $|| x_1 - x_2||_T < \epsilon$ when $|b_1 - b_2| + || y_1 - y_2||_T < \delta$,
where $(y_i, x_i)$ are two pairs of functions satisfying the relation \eqn{lip1}.
From \eqn{lip1}, we have
\begin{eqnarray}\label{c4}
|x_1 (t) - x_2 (t)| & \le & |b_1 - b_2| + |y_1 (t) - y_2 (t)| +  \int_{0}^{t} |h(x_1 (s)) - h(x_2 (s))| \, ds ~, \nonumber \\
& \le & |b_1 - b_2| + |y_1 (t) - y_2 (t)| +  c \int_{0}^{t} |x_1 (s) - x_2 (s)| \, ds ~.
\end{eqnarray}
Suppose that $|b_1 - b_2| + ||y_1 - y_2 ||_{T} \le \delta$.
By Gronwall's inequality,
\beq
|x_1 (t) - x_2 (t)| \le \delta e^{c t} \qandq  ||x_1  - x_2 ||_{T} \le \delta e^{c T} ~.
\eeqno
Hence it suffices to let $\delta = \epsilon e^{-c T}$.

We now turn to the Skorohod $J_1$ topology; see \S\S 3.3 and 11.5 and Chapter 12
of \cite{W02} for background.  To treat this non-uniform topology, we will use the fact that
the function $x$ is necessarily bounded.  That is proved later in Lemma \ref{lemSBIntRep}.
We want to show that $x_n \ra x$ in $D([0,\infty), \RR, J_1)$ when
$b_n \ra b$ in $\RR$ and $y_n \ra y$ in $D([0,\infty), \RR, J_1)$.  For $y$ given, let the interval
right endpoint $T$ be a continuity point of $y$.  Then there exist increasing homeomorphisms $\lambda_n$ of
the interval $[0,T]$ such that $\| y_n - y \circ \lambda_n \|_T \ra 0$ and $\| \lambda_n - e\|_T \ra 0$ as $n \ra \infty$.
Moreover, it suffices to consider homeomorphisms $\lambda_n$ that are absolutely continuous with respect to Lebesgue measure on $[0,T]$
having derivatives $\dot{\lambda}_n$ satisfying $\| \dot{\lambda}_n - 1\|_T \ra 0$ as $n \ra \infty$.
The fact that the topology is actually unchanged is a consequence
of Billingsley's equivalent complete metric $d_0$ on pp 112--114 of Billingsley \cite{B68}.
Hence, for $y$ given, let $M \equiv \sup_{0 \le t \le T}{\{ |x(t)|\}}$.
Since $h$ in \eqn{lip2} is Lipschitz, we have
\beq
\sup_{0 \le t \le T}{\{ |h(x(t))|\}} \le h(0) + \sup_{0 \le t \le T}{\{ |h(x(t)) - h(0)|\}} \le h(0) + c M = cM ~.
\eeqno
Thus we have
\begin{eqnarray}
|x_n (t) - x (\lambda_n (t))| & \le & |b_n - b| + \| y_n - y \circ \lambda_n\|_T \nonumber \\
&& \quad + \left| \int_{0}^{t} h(x_n (u)) \, du - \int_{0}^{\lambda_n (t)} h(x (u)) \, du \right| \nonumber \\
& \le &  |b_n - b| + \| y_n - y \circ \lambda_n\|_T  \nonumber \\
&& \quad + \left| \int_{0}^{t} h(x_n (u)) \, du - \int_{0}^{t} h(x (\lambda_n (u))) \dot{\lambda}_n (u) \, du \right| \nonumber \\
& \le &  |b_n - b| + \| y_n - y \circ \lambda_n\|_T + \|\dot{\lambda}_n - 1\|_T \int_{0}^{T} | h(x (u))|  \, du \nonumber \\
&& \quad +   \int_{0}^{t} |h(x_n (u)) - h(x (\lambda_n (u)))| \, du  \nonumber \\
& \le &  |b_n - b| + \| y_n - y \circ \lambda_n\|_T + \|\dot{\lambda}_n - 1\|_T (c M T) \nonumber \\
&& \quad +  c \int_{0}^{t} |x_n (u) - x (\lambda_n (u))| \, du ~. \nonumber
\end{eqnarray}

Choose $n_0$ such that $\|\dot{\lambda}_n - 1\|_T < \delta/(2 c M T)$
and $|b_n - b| + \| y_n - y \circ \lambda_n\|_T < \delta/2$.  Then Gronwall's inequality implies that
\beq
|x_n (t) - x (\lambda_n (t))| \le \delta e^{ct} \qforallq t, \quad 0 \le t \le T ~,
\eeqno
so that
\beq
\|x_n  - x \circ \lambda_n \|_T \le \delta e^{cT} ~.
\eeqno
Hence, for $\epsilon > 0$ given, choose $\delta < \epsilon e^{-cT}$ to have
$\|x_n  - x \circ \lambda_n \|_T \le \epsilon$ for $n \ge n_0$.  If necessary,
choose $n$ larger to make $\| \lambda_n - e\|_T < \epsilon$ and $\| \dot{\lambda}_n - 1\|_T < \epsilon$ as well.
Finally, for the inheritance of continuity, note that
\beq
 x(t+s) - x(t) = y(t+s) - y(t) + \int_{t}^{t +s} h(x (u)) \, du ~,
\eeqno
so that
\beq
| x(t+s) - x(t)| \le | y(t+s) - y(t)| + \int_{t}^{t +s}| h(x (u))| \, du ~.
\eeqno
Since $x$ is bounded over $[0,T]$, $x$ is continuous if $y$ is continuous.~~~\bsq

In our case we can simply let $h(s) = \mu s$, but we will need the more complicated function $h$ in \eqn{lip1} and \eqn{lip2}
 in \S \ref{secErlangA}.
To be self-contained, we now state a version of Gronwall's inequality; see p. 498 of \cite{EK86}.
See \S 11 of \cite{MMR98} for other versions of Gronwall's inequality.

\begin{lemma}{$($version of Gronwall's inequality$)$}\label{lemGronwall}
Suppose that $g: [0,\infty) \ra [0,\infty)$ is a Borel-measurable function
such that
\beq
0 \le g(t) \le \epsilon + M \int_{0}^{t} g(s) \, ds, \quad 0 \le t \le T ~,
\eeqno
for some positive finite $\epsilon$ and $M$.
Then
\beq
g(t) \le \epsilon e^{Mt}, \quad 0 \le t \le T ~.
\eeqno
\end{lemma}

It thus remains to establish the limit in \eqn{c1}.  Our proof based on the first martingale representation
in Theorem \ref{thMartRep}
relies on a FCLT for the Poisson process and the CMT
with the composition map.
The application of the CMT with the composition map requires a fluid limit, which requires further argument.
That is contained in subsequent sections.

\subsection{Poisson FCLT Plus the CMT}\label{secCTMCproof}
\hsp
As a consequence of the last section, it suffices to show that the scaled martingales
converge, as in \eqn{c1}.
From the martingale perspective, it is natural to achieve that goal by directly applying the martingale FCLT,
as in \S 7.1
of Ethier and Kurtz \cite{EK86}, and as reviewed here in \S \ref{secMartFCLT}, and that works.
In particular, the desired limit \eqn{c1} follows from Theorems
\ref{thMartRep} and \ref{thMart} (ii)
(or Theorems \ref{thMartRep2} and \ref{thMart} (ii))
plus Lemma \ref{lemFluid1} below.
Lemma \ref{lemFluid1} shows
 that the scaled
predictable quadratic variation processes in \eqn{b110} and \eqn{xb110} converge,
as required in condition \eqn{ek6} of Theorem \ref{thMart} here; see \S \ref{secProofSecond}.

However, starting with the first martingale representation in Theorem \ref{thMartRep},
we do not need to apply the martingale FCLT.
Instead, we can justify the martingale limit in \eqn{c1}
 by yet another application of
the CMT, using the composition map associated with the random time changes, in addition
to a functional central limit theorem (FCLT) for scaled Poisson processes.
Our approach also requires establishing a limit for the sequence of scaled
predictable quadratic variations associated with the martingales,
so the main steps of the
argument become the same as when applying the martingale FCLT.

The FCLT for
Poisson processes is a classical result.  It is a special case of the FCLT for a renewal process,
appearing as Theorem 17.3 in Billingsley \cite{B68}.  It and generalizations are also discussed
extensively in \cite{W02}; see \S\S 6.3, 7.3, 7.4, 13.7 and 13.8.
The FCLT for a Poisson process can also be obtained via a strong approximation,
as was done by Kurtz \cite{K78}, Mandelbaum and Pats \cite{MP95, MP98} and Mandelbaum, Massey and Reiman \cite{MMR98}.
Finally, the FCLT for a Poisson process itself can be obtained as an easy application of the martingale FCLT, as we show in
\S \ref{secMartFCLT}.

We start with the scaled Poisson processes

\beql{c8}
M_{A,n} (t)  \equiv  \frac{A (n t) - n t}{\sqrt{n}}  \qandq
M_{S,n} (t)  \equiv  \frac{S (n t) - n t}{\sqrt{n}},  \quad t \ge 0 ~.
\eeq

We employ the following basic FCLT:  Since $A$ and $S$ are independent rate-$1$ Poisson processes, we have

\begin{theorem}{$($FCLT for independent Poisson processes$)$}\label{thFCLTpoisson}
If $A$ and $S$ are independent rate-$1$ Poisson processes, then
\beql{c9}
(M_{A,n}, M_{S,n}) \Rightarrow (B_1, B_2) \qinq D^2 \equiv D \times D \qasq n \ra \infty ~,
\eeq
where $M_{A,n}$ and $M_{S,n}$ are the scaled processes in {\em \eqn{c8}},
while $B_1$ and $B_2$ are independent standard Brownian motions.
\end{theorem}

We can prove the desired limit in \eqn{c1} for both martingale representations,
but we will only give the details for the first martingale representation in Theorem \ref{thMartRep}.
In order to get the desired limit in \eqn{c1}, we introduce a deterministic and a random time change.
For that purpose, let $e:[0,\infty) \ra [0,\infty)$
be the identity function in $D$, defined by $e(t) \equiv t$ for $t \ge 0$.
Then let
\begin{eqnarray}\label{c10}
\Phi_{A,n} (t) & \equiv &   \frac{\lambda_n t}{n} = \mu t \equiv (\mu e) (t),  \nonumber \\
\Phi_{S,n} (t) & \equiv & \frac{\mu}{n} \int_{0}^{t} Q_n (s) \, d s, \quad t \ge 0 ~.
\end{eqnarray}

We will establish the following fluid limit, which can be regarded as a functional weak law of large numbers (FWLLN).
Here below, and frequently later, we have convergence in distribution to a deterministic limit;
that is \textbf{equivalent to convergence in probability}; see p. 27 of \cite{B68}.
\begin{lemma}{$($desired fluid limit$)$}\label{lemFluid1}
Under the conditions of Theorem {\em \ref{th1}},
\beql{c11}
\Phi_{S,n} \Rightarrow \mu e \qinq D \qasq n \ra \infty ~,
\eeq
where $\Phi_{S,n}$ is defined in {\em \eqn{c10}}.
\end{lemma}

For that purpose, it suffices to establish another more basic fluid limit.
Consider the stochastic process
\beql{c11a}
\Psi_{S,n} (t)  \equiv  \frac{Q_n (t)}{n}  ~, \quad t \ge 0 ~.
\eeq
Let $\omega$ be the function that is identically $1$ for all $t$.

\begin{lemma}{$($basic fluid limit$)$}\label{lemFluid2}
Under the conditions of Theorem {\em \ref{th1}},
\beql{c11b}
\Psi_{S,n} \Rightarrow \omega \qinq D \qasq n \ra \infty ~,
\eeq
where $\Psi_{S,n}$ is defined in {\em \eqn{c11a}}
and $\omega (t) = 1$, $t \ge 0$.
\end{lemma}

\paragraph{Proof of Lemma \ref{lemFluid1}.}
The desired fluid limit in Lemma \ref{lemFluid1} follows from the basic fluid limit
in Lemma \ref{lemFluid2} by applying the CMT with the function
$h: D \ra D$ defined by
\beql{c11c}
h(x) (t) \equiv \mu \int_{0}^{t} x(s) \, ds ~, \quad t \ge 0 ~.~~~\bsq
\eeq

We thus have the following result
\begin{lemma}{\em $($all but the fluid limit$)$}\label{lemFluid3}
If the limit in {\em \eqn{c11b}} holds, then
\beql{c11d}
(M_{n,1}, M_{n,2}) \Rightarrow ( \sqrt{\mu} B_1 , \sqrt{\mu} B_2 ) \qinq D^2
\eeq
as required to complete the proof of Theorem {\em \ref{th1}}.
\end{lemma}

\paragraph{Proof.}
From the limit in \eqn{c9}, the desired fluid limit in \eqn{c11} and
Theorem 11.4.5 of \cite{W02}, it follows that
\beql{c12}
(M_{A,n}, \mu e, M_{S,n},\Phi_{S,n}) \Rightarrow (B_1, \mu e, B_2, \mu e) \qinq D^4 \equiv D \times \cdots \times D
\eeq
as $n \ra \infty$.
From the CMT with the composition map, as in \S\S 3.4 and 13.2 of \cite{W02} - in particular, with Theorem 13.2.1 -
we obtain the desired limit in \eqn{c1}:
\beql{c13}
(M_{n,1}, M_{n,2}) \equiv (M_{A,n} \circ \mu e, M_{S,n} \circ \Phi_{S,n})
\Rightarrow ( B_1 \circ \mu e, B_2 \circ \mu e) \qinq D^2
\eeq
as $n \ra \infty$.
By basic properties of Brownian motion,
\beq
( B_1 \circ \mu e, B_2 \circ \mu e) \deq  ( \sqrt{\mu} B_1 , \sqrt{\mu} B_2 ) \qinq D^2 ~.~~~\bsq
\eeqno

It thus remains to establish the key fluid limit in Lemma \ref{lemFluid2}.  In the next section
we show how to do that directly, without martingales, by applying the continuous mapping provided
by Theorem \ref{thCMT} in the fluid scale or, equivalently, by applying Gronwall's inequality again.
We would stop there if we only wanted to analyze the $M/M/\infty$ model,
but in order to illustrate other methods used in Krichagina and Puhalskii \cite{KP97}
and Puhalskii and Reiman \cite{PR00}, we also apply martingale methods.
Thus, in the subsequent
 four sections
we show how to establish that fluid limit using martingales.  Here is an \textbf{outline of the remaining martingale argument:}

\vspace{0.1in}
(1) To prove the needed Lemma \ref{lemFluid2}, it suffices to demonstrate that $\{X_n\}$
is stochastically bounded in $D$. (Combine Lemma \ref{lemSBfluid} and \S \ref{secFluidSB}.)

\vspace{0.1in}
(2)  However, $\{X_n\}$ is stochastically bounded in $D$ if the sequences of martingales $\{M_{n,1}\}$ and $\{M_{n,2}\}$
are stochastically bounded in $D$. (Combine Theorem \ref{thMartRep} and Lemma \ref{lemSBIntRep}.)

\vspace{0.1in}
(3)  But then the sequences of martingales $\{M_{n,1}\}$ and $\{M_{n,2}\}$
are stochastically bounded in $D$ if the associated sequences of predictable quadratic variations
$\{\langle M_{n,1} \rangle (t) \}$ and $\{\langle M_{n,2} \rangle (t)\}$ are stochastically bounded in $\RR$ for each $t > 0$
(Apply Lemma \ref{lemSBLenglart}.  One of these is trivial because it is deterministic.)

\vspace{0.1in}
(4)  Finally, we establish stochastic boundedness of $\{\langle M_{n,2} \rangle (t)\}$
(the one nontrivial case) through a crude bound in \S \ref{secSBcomp}.

\vspace{0.1in}
This alternate route to the fluid limit is much longer, but all the steps might be considered well known.
We remark that the fluid limit seems to be required by any of the proofs, including
by the direct application of the martingale FCLT.

\subsection{Fluid Limit Without Martingales}\label{secSLLN}
\hsp
In this section we prove Lemma \ref{lemFluid2} without using martingales.
We do so by
establishing a stochastic-process limit in the fluid scale which is similar to
the corresponding stochastic-process limit with the more refined scaling.
This is a standard line of reasoning for heavy-traffic stochastic-process limits;
e.g., see the proofs of Theorems 9.3.4, 10.2.3 and 14.7.4 of Whitt \cite{W02}.
The specific argument here follows \S 6 of Mandelbaum and Pats \cite{MP95}.
With this approach, even though we exploit the martingale representations,
we do not need to mention martingales at all.  We are only applying the continuous mapping theorem.

By essentially the same reasoning as in \S \ref{secFirst}, we obtain a fluid-scale analog of \eqn{b108} and \eqn{b109}:
\begin{eqnarray}\label{st1}
\bar{X}_n (t) & \equiv & \frac{Q_n (t) - n}{n} \nonumber \\
& = & \frac{Q_n (0) - n}{n} + \frac{M^{*}_{n,1} (t)}{n}
- \frac{M^{*}_{n,2} (t)}{n}  + \frac{\lambda_n t - n \mu t}{n} \nonumber \\
&& \quad - \mu \int_{0}^{t} \left(\frac{Q_n (s)- n}{n}\right) \, d s ~, \nonumber \\
& = & \bar{X}_n (0) + \bar{M}_{n,1} (t) - \bar{M}_{n,2} (t)   - \mu \int_{0}^{t} \bar{X}_n (s) \, d s ~, \quad t \ge 0 ~,
\end{eqnarray}
where
\beql{st2}
\bar{M}_{n,1} (t) \equiv \frac{M^{*}_{n,1} (t)}{n} \qandq
\bar{M}_{n,2} (t) \equiv \frac{M^{*}_{n,2} (t)}{n}, \quad t \ge 0 ~,
\eeq
with $M^{*}_{n,i} (t)$ defined in \eqn{b106}.

Notice that the limit
\beql{st3}
\bar{X}_n \Rightarrow \eta \qinq D^k \qasq n \ra \infty ~,
\eeq
where
\beql{st4}
\eta (t) \equiv 0, \quad t \ge 0 ~,
\eeq
 is equivalent to
the desired conclusion of Lemma \ref{lemFluid2}.
Hence we will prove the fluid limit in \eqn{st3}.

The assumed limit in \eqn{a2a} implies that $\bar{X}_n (0) \Rightarrow 0$ in $\RR$ as $n \ra \infty$.
We can apply Theorem \ref{thCMT} or directly Gronwall's inequality in Lemma \ref{lemGronwall}
to deduce the desired limit \eqn{st3} if we can establish the following
lemma.

\begin{lemma}{$($fluid limit for the martingales$)$}\label{lemFluid11}
Under the conditions of Theorem {\em \ref{th1}},
\beql{st5}
\bar{M}_{n,i} \Rightarrow \eta  \qinq D \quad \mbox{w.p.1} \qasq n \ra \infty ~,
\eeq
for $i = 1, 2$, where $\bar{M}_{n,i}$ is defined in {\em \eqn{st2}} and $\eta$ is defined in {\em \eqn{st4}}.
\end{lemma}

\paragraph{Proof of Lemma \ref{lemFluid11}.}
We can apply the SLLN for the Poisson process, which is equivalent to the
more general functional strong law of large numbers (FSLLN); see \S 3.2 of \cite{W02IS}.
(Alternatively, we could apply the FWLLN, which is a corollary to the FCLT.)
First, the SLLN for the Poisson process states that
\beq
\frac{A(t)}{t} \ra 1 \qandq \frac{S(t)}{t} \ra 1 \quad \mbox{w.p.1} \qasq t \ra \infty ~,
\eeqno
which implies the corresponding FSLLN's
\beql{st7}
\sup_{0 \le t \le T}{\{\frac{A(nt)}{n} - t \}} \ra 0 \qandq \sup_{0 \le t \le T}{\{\frac{S(nt)}{n} - t \}} \ra 0 \quad \mbox{w.p.1}
\eeq
as $n \ra \infty$ for each $T$ with $0 < T < \infty$.
We thus can treat $\bar{M}_{n,1}$ directly.  To treat $\bar{M}_{n,2}$,
we combine \eqn{st7} with the crude inequality in \eqn{b102} and the representation in \eqn{st1}--\eqn{st2} in order
to obtain the desired limit \eqn{st5}.  To elaborate, the crude inequality in \eqn{b102}
implies that, for any $T_1 > 0$,
there exists $T_2$ such that
\beq
P\left(\frac{\mu}{n} \int_{0}^{T_1} Q_n (s) \, ds > T_2\right) \ra 0 \qasq n \ra \infty ~.
\eeqno
That provides the key, because
\beq
P\left( \| \bar{M}_{n,2} \|_{T_1} >  \epsilon \right) \le P\left(\frac{\mu}{n} \int_{0}^{T_1} Q_n (s) \, ds > T_2\right)
+ P\left(\| \bar{S}_n \|_{T_2} > \epsilon/2\right) ~,
\eeqno
where
\beq
\bar{S}_n (t) \equiv \frac{S(nt) - nt}{n}, \quad t \ge 0 ~.~~~\bsq
\eeqno


\section{Tightness and Stochastic Boundedness}\label{secTightSB}

\subsection{Tightness}\label{secTight}
\hsp
As indicated at the end of \S \ref{secCTMCproof}, we can also use a stochastic-boundedness argument in order to establish the desired fluid limit.
Since stochastic boundedness
is closely related to tightness, we start
by reviewing tightness concepts.  In the next section we apply the tightness notions
to stochastic boundedness.  The next three sections contain extra material not really needed for the current
proofs.  Additional material on tightness criteria appears in Whitt \cite{W07}.

We work in the setting of a complete separable metric space (CSMS),
also known as a Polish space; see \S\S 13 and 19 of Billingsley \cite{B68}, \S\S 3.8-3.10 of Ethier and Kurtz \cite{EK86}
and \S\S 11.1 and 11.2 of \cite{W02}.
(The space $D^k \equiv D([0, \infty), \RR)^k$ is made a CSMS in a standard way and the space of probability measures on $D^k$
becomes a CSMS as well.)  Key concepts are:  closed, compact, tight, relatively compact and sequentially compact.
We assume knowledge of metric spaces and compactness in metric spaces.

\begin{definition} {\em $($tightness$)$}\label{defTight}
A set $A$ of probability measures on a metric space $S$ is \textbf{tight} if, for all $\epsilon > 0$,
there exists a compact subset $K$ of $S$ such that
\beq
P(K) > 1 - \epsilon \qforallq P \in A ~.
\eeqno
A set of random elements of the metric space $S$ is tight if the associated set of their probability laws on $S$ is tight.
Consequently, a
sequence $\{X_n: n \ge 1\}$ of random elements of the metric space $S$ is tight if, for all $\epsilon > 0$,
there exists a compact subset $K$ of $S$ such that
\beq
P(X_n \in K) > 1 - \epsilon \qforallq n \ge 1 ~.
\eeqno
\end{definition}

Since a continuous image of a compact subset is compact, we have the following lemma.

\begin{lemma}{\em $($continuous functions of random elements$)$}\label{lemTightCont}
Suppose that $\{X_n: n \ge 1\}$ is a tight sequence of random elements of the metric space $S$.
If $f: S \ra S'$ is a continuous function mapping the metric space $S$ into another metric space $S'$,
then $\{f(X_n): n \ge 1\}$ is a tight sequence of random elements of the metric space $S'$.
\end{lemma}

\paragraph{Proof.}  As before, let $\circ$ be used for composition:  $(f \circ g) (x) \equiv f(g(x))$.
 For any function $f: S \ra S'$ and any subset $A$ of $S$, $A \subseteq f^{-1} \circ f (A)$.
Let $\epsilon > 0$ be given.
Since $\{X_n: n \ge 1\}$ is a tight sequence of random elements of the metric space $S$, there exists
a compact subset $K$ of $S$ such that
\beq
P(X_n \in K) > 1 - \epsilon \qforallq n \ge 1 ~.
\eeqno
Then $f(K)$ will serve as the desired compact set in $S'$, because
\beq
P(f(X_n) \in f(K)) = P(X_n \in (f^{-1} \circ f)(K)) \ge  P(X_n \in K) > 1 - \epsilon 
\eeqno
for all $n \ge 1$.~~~\bsq

We next observe that on products of separable metric spaces tightness is characterized by tightness of the components;
see \S 11.4 of \cite{W02}.

\begin{lemma}{\em $($tightness on product spaces$)$}\label{lemTightProd}
Suppose that $\{(X_{n,1}, \ldots , X_{n,k}) : n \ge 1\}$ is a sequence of random elements of the product space $S_1 \times \cdots \times S_k$,
where each coordinate space $S_i$ is a separable metric space.
The sequence $\{(X_{n,1}, \ldots , X_{n,k}) : n \ge 1\}$ is tight if and only if the sequence $\{X_{n,i} : n \ge 1\}$ is tight for each $i$,
$1 \le i \le k$.
\end{lemma}

\paragraph{Proof.}
The implication from the random vector to the components follows from Lemma \ref{lemTightCont} because
the component $X_{n,i}$ is the image of the projection map $\pi_i: S_1 \times \cdots \times S_k \ra S_i$
taking $(x_1, \ldots , x_k)$ into $x_i$, and the projection map is continuous.  Going the other way, we use the
fact that
\beq
A_1 \times \cdots \times A_k = \bigcap_{i = 1}^{k} \pi^{-1}_i (A_i) = \bigcap_{i = 1}^{k} \pi^{-1}_i \circ \pi_i (A_1 \times \cdots \times A_k)
\eeqno
for all subsets $A_i \subseteq S_i$.  Thus, for each $i$ and any $\epsilon > 0$, we can choose $K_i$
such that $P(X_{n,i} \notin K_i) < \epsilon/k$ for all $n \ge 1$.  We then let $K_1 \times \cdots \times K_k$
be the desired compact for the random vector.  We have
\begin{eqnarray}
P\left((X_{n,1}, \ldots , X_{n,k}) \notin K_1 \times \cdots \times K_k\right) & = & P\left(\bigcup_{i=1}^{k} \{X_{n,i} \notin K_i\} \right) \nonumber \\
& \le & \sum_{i=1}^{k} P\left( X_{n,i} \notin K_i \right) \le \epsilon ~.~~~\bsq \nonumber
\end{eqnarray}




Tightness goes a long way toward establishing convergence because of Prohorov's theorem.  It involves the notions of sequential compactness and relative
compactness.

\begin{definition} {\em $($relative compactness and sequential compactness$)$}\label{defRelCompact}
A subset $A$ of a metric space $S$ is \textbf{relatively compact} if every sequence $\{x_n: n \ge 1\}$ from $A$
has a subsequence that converges to a limit in $S$ $($which necessarily belongs to the closure $\bar{A}$ of $A)$.
\end{definition}

We can now state Prohorov's theorem; see \S 11.6 of \cite{W02}.  It relates compactness of sets of measures to
compact subsets of the underlying sample space $S$ on which the probability measures are defined.

\begin{theorem}{\em $($Prohorov's theorem$)$}\label{thmProhorov}
A subset of probability measures on a CSMS is tight if and only if it is
relatively compact.
\end{theorem}

We have the following elementary corollaries:

\begin{coro}{\em $($convergence implies tightness$)$}\label{corConvTight}
If $X_n \Rightarrow X$ as $n \ra \infty$ for random elements of
a CSMS, then the sequence $\{X_n: n \ge 1\}$ is tight.
\end{coro}

\begin{coro}{\em $($individual probability measures$)$}\label{corIndiv}
Every individual probability measure on
a CSMS is tight.
\end{coro}

As a consequence of Prohorov's Theorem, we have the following method for establishing convergence of random elements:

\begin{coro}{\em $($convergence in distribution via tightness$)$}\label{corProhorov}
Let $\{X_n: n \ge 1\}$ be a sequence of random elements of a CSMS $S$.  We have
\beq
X_n \Rightarrow X \qinq S \qasq n \ra \infty
\eeqno
if and only if (i) the sequence $\{X_n: n \ge 1\}$ is tight and (ii) the limit of every convergent subsequence of $\{X_n: n \ge 1\}$
is the same fixed random element $X$ (has a common probability law).
\end{coro}

In other words, once we have established tightness, it only remains to show that the limits of all converging subsequences must be the same.
With tightness, we only need to uniquely determine the limit.  When proving Donsker's theorem, it is natural
to uniquely determine the limit through the finite-dimensional distributions.  Convergence of all the finite-dimensional distributions
is not enough to imply convergence on $D$, but it does uniquely determine the distribution of the limit; see pp 20 and 121 of
Billingsley \cite{B68} and
Example 11.6.1 in \cite{W02}.

This approach is applied to prove the martingale FCLT stated in \S \ref{secMartFCLT}; see \cite{W07}.
In the martingale setting it is natural instead to use the martingale characterization of Brownian motion,
originally established by L\'{e}vy \cite{L48} and proved by Ito's formula by Kunita and Watanabe \cite{KW67};
see p. 156 of Karatzas and Shreve \cite{KS88}, and various extensions, such as to continuous processes
with independent Gaussian increments, as in Theorem 1.1 on p. 338 of Ethier and Kurtz \cite{EK86}.
A thorough study of martingale characterizations appears in Chapter 4 of Liptser and Shiryayev \cite{LS89}
and in Chapters VIII and IX of Jacod and Shiryayev \cite{JS87}.

We have not discussed conditions to have tightness; they are reviewed in \cite{W07}.

\subsection{Stochastic Boundedness}\label{secSB}
\hsp
We start by defining stochastic boundedness and relating it to tightness.
We then discuss situations in which stochastic boundedness is preserved.
Afterwards, we give conditions for a sequence of martingales to be stochastically bounded in $D$
involving the stochastic boundedness of appropriate sequences of $\RR$-valued random variables.
Finally, we show that the FWLLN follows from stochastic boundedness.

\subsubsection{Connection to Tightness}

 For random elements of $\RR$ and $\RR^k$, stochastic boundedness and tightness
 are equivalent, but tightness is stronger than stochastic boundedness for
random elements of the functions spaces $C$ and $D$ (and the associated product spaces $C^k$ and $D^k$).

\begin{definition} {\em $($stochastic boundedness for random vectors$)$}\label{defSB}
A sequence $\{X_n: n \ge 1\}$ of random vectors taking values in $\RR^k$ is \textbf{stochastically bounded (SB)} if the sequence is tight, as defined in
Definition {\em \ref{defTight}}.
\end{definition}

The notions of tightness and stochastic boundedness thus agree for random elements of $\RR^k$, but these
notions differ for stochastic processes.
For a function $x \in D^k \equiv D([0, \infty), \RR)^k$, let
\beq
\| x\|_T \equiv \sup_{0 \le t \le T}{\{ | x (t) |\}} ~,
\eeqno
where $|b|$ is a norm of $b \equiv (b_1, b_2, \ldots , b_k)$ in $\RR^k$ inducing the Euclidean topology, such as the maximum norm:
$|b| \equiv \max{\{|b_1|, |b_2|, \ldots , |b_k|\}}$.  (Recall that all norms on Euclidean space $\RR^k$ are equivalent.)

\begin{definition} {\em $($stochastic boundedness for random elements of $D^k)$}\label{defSBD}
A sequence $\{X_n: n \ge 1\}$ of random elements of $D^k$ is \textbf{stochastically bounded in $D^k$} if the sequence
of real-valued random variables $\{\| X_n\|_T: n \ge 1\}$
is stochastically bounded in $\RR$ for each $T > 0$, using Definition {\em \ref{defSB}}.
\end{definition}

For random elements of $D^k$, tightness is a strictly stronger concept than stochastic boundedness.
Tightness of $\{X_n\}$ in $D^k$
implies stochastic boundedness, but not conversely; see \S 15 of Billingsely \cite{B68}.  However, stochastic boundedness is sufficient
for us, because it alone implies the desired fluid limit.

\subsubsection{Preservation}

We have the following analog of Lemma \ref{lemTightProd}, which characterizes tightness for sequences of random vectors in terms
of tightness of the associated sequences of components.

\begin{lemma}{\em $($stochastic boundedness on $D^k$ via components$)$}\label{lemSBprod}
A sequence
$$\{(X_{n,1}, \ldots , X_{n,k}) : n \ge 1\} \qinq D^k \equiv D \times \cdots \times D$$
is stochastically bounded in $D^k$ if and only if the sequence $\{X_{n,i} : n \ge 1\}$
is stochastically bounded in $D \equiv D^1$ for each $i$,
$1 \le i \le k$.
\end{lemma}

\paragraph{Proof.}  Assume that we are using the maximum norm on product spaces.  We can apply Lemma \ref{lemTightProd} after noticing that
\beq
\| (x_1, \ldots , x_k)\|_T = \max{\{ \| x_i\|_T: 1 \le i \le k\}}
\eeqno
for each element $(x_1, \ldots , x_k)$ of $D^k$.  Since other norms are equivalent, the result applies more generally.~~~\bsq

\begin{lemma}{\em $($stochastic boundedness in $D^k$ for sums$)$}\label{lemSBsums}
Suppose that
\beq
Y_n (t) \equiv X_{n,1} (t) + \cdots + X_{n,k} (t), \quad t \ge 0,
\eeqno
for each $n \ge 1$, where
$\{(X_{n,1}, \ldots , X_{n,k}) : n \ge 1\}$ is a sequence of random elements of the product space $D^k \equiv D \times \cdots \times D$.
If $\{X_{n,i} : n \ge 1\}$
is stochastically bounded in $D$ for each $i$,
$1 \le i \le k$, then the sequence $\{Y_n: n \ge 1\}$
is stochastically bounded in $D$.
\end{lemma}

Note that the converse is not true:  We could have $k=2$ with $X_{n,2} (t) = - X_{n,1} (t)$ for all $n$ and $t$.
In that case we have $Y_n (t) = 0$ for all $X_{n,1} (t)$.

We now provide conditions for the stochastic boundedness of integral representations such as \eqn{b109}.

\begin{lemma}{\em $($stochastic boundedness for integral representations$)$}\label{lemSBIntRep}
Suppose that
\beq
X_n (t) \equiv X_n (0) + Y_{n,1} (t) + \cdots + Y_{n,k} (t) +  \int_{0}^{t} h(X_n (s)) \, ds, \quad t \ge 0 ~,
\eeqno
 where $h$ is a Lipschitz function as in {\em \eqn{lip2}}
 and $(X_n (0), (Y_{n,1}, \ldots , Y_{n,k}) )$ is a random element of $\RR \times D^k$ for each $n \ge 1$.
 If the sequences $\{X_n (0): n\ge 1\}$ and $\{Y_{n,i}: n\ge 1\}$ are stochastically bounded (in $\RR$ and $D$, respectively,)
 for $1 \le i \le k$, then the sequence $\{X_n: n \ge 1\}$ is stochastically bounded in $D$.
 \end{lemma}

 \paragraph{Proof.}  For the proof here, it is perhaps easiest to imitate the proof of Theorem \ref{thCMT}.
 In particular, we will construct the following bound
 \beq
||X_n ||_{T} \le K e^{c T},
\eeqno
assuming that $|X_n (0)| + ||Y_{n,1} ||_{T} + \cdots + ||Y_{n,k} ||_{T} + T |h(0)| \le K$.
As before, we apply Gronwall's inequality, Lemma \ref{lemGronwall}.
Inserting absolute values into all the terms of the integral representation,
we have
\begin{eqnarray}
|X_n (t)| & \le & |X_n (0)| + |Y_{n,1} (t)| + \cdots + |Y_{n,k} (t)| +  \int_{0}^{t}|h(X_n (s))| \, ds,  \nonumber \\
          & \le & |X_n (0)| + \|Y_{n,1} \|_T + \cdots + \|Y_{n,k} \|_T + T h(0) + c \int_{0}^{t}|X_n (s)| \, ds ~, \nonumber
\end{eqnarray}
for $t \ge 0$, where $c$ is the Lipschitz modulus.
Suppose that $|X_n (0)| + ||Y_{n,1} ||_{T} + \cdots + ||Y_{n,k} ||_{T} + T |h(0)| \le K$.
By Gronwall's inequality,
\beq
|X_n (t)| \le K e^{c t} \qandq  ||X_n ||_{T} \le K e^{c T} ~.~~~\bsq
\eeqno

\subsubsection{Stochastic Boundedness for Martingales}

We now provide ways to get stochastic boundedness for sequences of martingales in $D$
from associated sequences of random variables.  Our first result exploits the classical submartingale-maximum inequality;
e.g., see p. 13 of Karatzas and Shreve \cite{KS88}.  We say that a function $f: \RR \ra \RR$ is {\em even} if
$f(-x) = f(x)$ for all $x \in \RR$.

\begin{lemma}{\em $($SB from the maximum inequality$)$}\label{lemSBmax}
Suppose that, for each $n \ge 1$, $M_n \equiv \{M_n (t): t \ge 0\}$ is a martingale $($with respect to a specified filtration$)$
with sample paths in $D$.
Also suppose that, for each $T > 0$, there exists an even nonnegative convex function $f: \RR \ra \RR$
with first derivative $f' (t) > 0$ for $t > 0$ $($e.g., $f(t) \equiv t^2)$,
there exists a positive constant $K \equiv K(T, f)$, and there exists
an integer $n_0 \equiv n_0 (T,f, K)$,
 such that
\beq
E[f(M_n(T))] \le K \qforallq n \ge n_0 ~.
\eeqno
Then the sequence of stochastic processes $\{M_n: n \ge 1\}$ is stochastically bounded in $D$.
 \end{lemma}

 \paragraph{Proof.}  Since any set of finitely many random elements of $D$ is automatically tight,
 Theorem 1.3 of Billingsley \cite{B68}, it suffices to consider $n \ge n_0$.
 Since $f$ is continuous and $f'(t) >0 $ for $t > 0$, $t > c$ if and only if
 $f(t) > f(c)$ for $t > 0$.  Since $f$ is even,
 \beq
 E[f(M_n(t))] = E[f(| M_n(t)| )] \le E[f(|M_n(T)|)] = E[f(M_n(T))]\le K
 \eeqno
 for all $t$, $0 \le t \le T$.
 Since these moments are finite and $f$ is convex,
 the stochastic process $\{f(M_n(t)) : 0 \le t \le T\}$ is a submartingale for each $n \ge 1$,
 so that we can apply the submartingale-maximum inequality to get
 \beq
 P(\| M_n \|_T > c) =  P(\| f \circ M_n \|_T > f(c)) \le \frac{E[f(M_n(T))]}{f(c)} \le \frac{K}{f(c)}
 \eeqno
 for all $n \ge n_0$.
 Since $f(c) \ra \infty$ as $c \ra \infty$, we have the desired conclusion.~~~\bsq

We now establish another sufficient condition for stochastic boundedness of
square-integrable martingales by applying the Lenglart-Rebolledo inequality;
see p. 66 of Liptser and Shiryayev \cite{LS89} or p. 30 of Karatzas and Shreve \cite{KS88}.

\begin{lemma}{\em $($Lenglart-Rebolledo inequality$)$}\label{lemLenglart}
Suppose that $M \equiv \{M (t): t \ge 0\}$ is a square-integrable martingale $($with respect to a specified filtration$)$
with predictable quadratic variation $\langle M \rangle \equiv \{\langle M \rangle (t): t \ge 0\}$, i.e.,
such that $M^2 - \langle M \rangle \equiv \{ M (t)^2 - \langle M \rangle (t): t \ge 0\}$ is a martingale
by the Doob-Meyer decomposition.  Then, for all $c>0$ and $d>0$,
\beql{lenglart}
P \left( \sup_{0 \le t \le T}{\{ |M(t)|\}} > c \right) \le \frac{d}{c^2} + P\left(\langle M \rangle (T) > d \right) ~.
\eeq
 \end{lemma}

As a consequence we have the following criterion for stochastic boundedness of a sequence of square-integrable martingales.

\begin{lemma}{\em $($SB criterion for square-integrable martingales$)$}\label{lemSBLenglart}
Suppose that, for each $n \ge 1$, $M_n \equiv \{M_n (t): t \ge 0\}$ is a square-integrable martingale $($with respect to a specified filtration$)$
with predictable quadratic variation $\langle M_n \rangle \equiv \{\langle M_n \rangle (t): t \ge 0\}$, i.e.,
such that $M_n^2 - \langle M_n \rangle \equiv \{ M_n (t)^2 - \langle M_n \rangle (t): t \ge 0\}$ is a martingale
by the Doob-Meyer decomposition.  If the sequence of random variables $\{\langle M_n \rangle (T): n \ge 1\}$
is stochastically bounded in $\RR$ for each $T > 0$, then the sequence of stochastic processes $\{M_n: n \ge 1\}$ is stochastically bounded
in $D$.
 \end{lemma}

\paragraph{Proof.}
For $\epsilon > 0$ given, apply the assumed stochastic boundedness of
the sequence $\{\langle M_n \rangle (T): n \ge 1\}$ to obtain a constant $d$ such that
\beq
P\left(\langle M_n \rangle (T) > d \right) < \epsilon/2 \qforallq n \ge 1 ~.
\eeqno
Then for that determined $d$, choose $c$ such that $d/c^2 < \epsilon/2$.  By the
Lenglart-Rebolledo inequality \eqn{lenglart}, these two inequalities imply that
\beq
P \left( \sup_{0 \le t \le T}{\{ |M_n(t)|\}} > c \right) < \epsilon ~.~~~\bsq
\eeqno

\subsubsection{FWLLN from Stochastic Boundedness}

We will want to apply stochastic boundedness in $D$ to
imply the desired fluid limit in
Lemmas \ref{lemFluid1} and \ref{lemFluid2}.  The fluid limit corresponds
to a \textbf{functional weak law of large numbers (FWLLN)}.

\begin{lemma}{\em $($FWLLN from stochastic boundedness in $D^k)$}\label{lemSBfluid}
Let $\{X_n: n \ge 1\}$ be a sequence of random elements of $D^k$.  Let $\{a_n: n \ge 1\}$ be a sequence of
positive real numbers such that $a_n \ra \infty$ as $n \ra \infty$.  If the sequence $\{X_n: n \ge 1\}$
is stochastically bounded in $D^k$ , then
\beql{m1}
\frac{X_n}{a_n} \Rightarrow \eta \qinq D^k \qasq n \ra \infty ~,
\eeq
where $\eta (t) \equiv (0, 0, \ldots , 0)$, $t \ge 0$.
\end{lemma}

\paragraph{Proof.}
As specified in Definition \ref{defSBD}, stochastic boundedness of the sequence $\{X_n: n \ge 1\}$
in $D^k$ corresponds to stochastic boundedness
of the associated sequence $\{\| X_n\|_T: n \ge 1\}$ in $\RR$ for each $T > 0$.  By Definition \ref{defSB}, stochastic boundedness in $\RR$
is equivalent to tightness.
It is easy to verify directly that we then have tightness (or, equivalently, stochastic boundedness)
for the associated sequence $\{\| X_n\|_T/a_n: n \ge 1\}$ in $\RR$.
By Prohorov's theorem, Theorem \ref{thmProhorov}, tightness on $\RR$ (or any CSMS)
is equivalent to relative compactness.  Hence consider a convergent subsequence $\{\| X_{n_k}\|_T/a_{n_k}: k \ge 1\}$ of
the sequence $\{\| X_n\|_T/a_n: n \ge 1\}$ in $\RR$:  $\| X_{n_k}\|_T/a_{n_k} \rightarrow L$ as $k \ra \infty$.
It suffices to show that $P(L = 0) = 1$; then all convergent subsequences will have the same limit,
which implies convergence to that limit.  For that purpose, consider the associated subsequence
$\{\| X_{n_k}\|_T: k \ge 1\}$ in $\RR$.  It too is tight.  So by Prohorov's theorem again,
it too is relatively compact.  Thus there exists a convergent sub-subsequence: $\| X_{n_{k_l}}\|_T \Rightarrow L'$ in $\RR$.
It follows immediately that
\beql{m34}
\frac{\| X_{n_{k_l}}\|_T}{a_{n_{k_l}}} \Rightarrow 0 \qinq \RR \qasq l \ra \infty ~.
\eeq
This can be regarded as a consequence of the \textbf{generalized continuous mapping theorem (GCMT)},
Theorem 3.4.4 of \cite{W02}, which involves a sequence of continuous functions:  Consider the functions $f_n: \RR \ra \RR$ defined by
$f_n (x) \equiv x/a_{n}$, $n \ge 1$, and the limiting zero function $f: \RR \ra \RR$ defined by
$f (x) \equiv 0$ for all $x$.  It is easy to see that $f_n (x_n) \ra f(x) \equiv 0$
whenever $x_n \ra x$ in $\RR$.  Thus the GCMT implies the limit in \eqn{m34}.
Consequently, the limit $L$ we found for the subsequence $\{\| X_{n_k}\|_T/a_{n_k}: k \ge 1\}$ must actually be $0$.
Since that must be true for all convergent subsequences, we must have the claimed convergence in \eqn{m1}.~~~\bsq

\section[Completing the Proof of Theorem 1.1]{Completing the Proof of Theorem \ref{th1}}\label{secComplete}
\hsp
In the next three subsections we complete the proof of Theorem \ref{th1} by the martingale argument,
as outlined at the end of \S \ref{secCTMCproof}.
In \S \ref{secFluidSB} we show that the required fluid limit in Lemma \ref{lemFluid2}
follows from the stochastic boundedness of the the sequence of stochastic processes $\{X_n: n \ge 1\}$.
In \S \ref{secSBcomp} we finish the proof of the fluid limit by proving that the
associated predictable quadratic variation processes are stochastically bounded.
Finally, in \S \ref{secInitial} we show how to remove the condition on the initial conditions.

\subsection{Fluid Limit from Stochastic Boundedness in D}\label{secFluidSB}
\hsp
We now show how to apply stochastic boundedness to imply the desired fluid limit in
Lemma \ref{lemFluid2}.  Here we simply apply Lemma \ref{lemSBfluid}
with $a_n \equiv \sqrt{n}$ for $n \ge 1$
to our particular sequence of stochastic processes $\{X_n: n \ge 1\}$
in \eqn{a2} or \eqn{b109}.

\begin{lemma}{\em $($application to queueing processes$)$}\label{lemQfluid}
Let $X_n$ be the random elements of $D$ defined in {\em \eqn{a2}} or {\em \eqn{b109}}.
If the sequence $\{X_n: n \ge 1\}$
is stochastically bounded in $D$, then
\beql{y2}
\|n^{-1} Q_n - \omega \|_T \Rightarrow 0 \qasq n \ra \infty \qforallq T > 0 ~,
\eeq
where $\omega (t) \equiv 1$, $t \ge 0$.
Equivalently,
\beql{y3}
\sup_{0 \le t \le T}{\{| (Q_n (t)/n) - 1|\}} \Rightarrow 0 \qasq n \ra \infty \qforallq T > 0 ~.
\eeq

\end{lemma}

\paragraph{Proof.}  As a consequence of Lemma \ref{lemSBfluid}, From the stochastic boundedness of
$\{X_n\}$ in $D$, we will obtain
\beq
\frac{X_n}{\sqrt{n}} \Rightarrow \eta \qinq D \qasq n \ra \infty ~,
\eeqno
where $\eta$ is the zero function defined above.  This limit is equivalent to \eqn{y2} and \eqn{y3}.~~~\bsq

In \S \ref{secCTMCproof} we already observed that Lemma \ref{lemFluid2} implies the desired Lemma \ref{lemFluid1}.
So we have completed the required proof.

\paragraph{Recap:  How this all goes together.}

The consequence of the last three sections is as follows:
In order to complete the martingale argument to establish the required fluid limit
in Lemmas \ref{lemFluid1} and \ref{lemFluid2}, it suffices to establish the stochastic boundedness of two sequences of random variables:
$\{\langle M_{n,1} \rangle (T): n \ge 1\}$ and $\{\langle M_{n,2} \rangle (T): n \ge 1\}$ for any $T > 0$.

Here is an explanation:
For each $n \ge 1$, the associated stochastic processes $\{\langle M_{n,1} \rangle (t) : t \ge 0\}$ and $\{\langle M_{n,2} \rangle (t) : t \ge 0\}$
are the predictable quadratic variations (compensators)
of the scaled martingales  $M_{n,1}$ and $M_{n,2}$
in \eqn{b106}, \eqn{b107} and \eqn{c1}.  First, by Lemma \ref{lemSBLenglart}, the stochastic boundedness
of these two sequences of PQV random variables implies stochastic boundedness in $D$ of the two sequences of
scaled martingales $\{M_{n,1}: n \ge 1\}$ and $\{M_{n,2}: n \ge 1\}$
in \eqn{b107}.   That, in turn, under condition \eqn{a2a}, implies the stochastic boundedness in $D$ of the sequence
of scaled queue-length processes $\{X_n: n \ge 1\}$ in \eqn{b109} and \eqn{xb109} by Lemma \ref{lemSBIntRep}.
Finally, by Lemma \ref{lemSBfluid} the stochastic boundedness of $\{X_n\}$ in $D$ implies the required fluid limit in \eqn{y2}.
The fluid limit in \eqn{y2} was just what was needed in Lemmas \ref{lemFluid1} and \ref{lemFluid2}.

 \subsection{Stochastic Boundedness of the Quadratic Variations}\label{secSBcomp}
 \hsp
We have observed at the end of the last section that it only remains to
establish the stochastic boundedness of two sequences of PQV random variables:
$\{\langle M_{n,1} \rangle (t): n \ge 1\}$ and $\{\langle M_{n,2} \rangle (t): n \ge 1\}$ for any $t > 0$.
For each $n \ge 1$, the associated stochastic processes $\{\langle M_{n,1} \rangle (t) : t \ge 0\}$ and $\{\langle M_{n,2} \rangle (t) : t \ge 0\}$
are the predictable quadratic variations (compensators)
of the scaled martingales  $M_{n,1}$ and $M_{n,2}$
in \eqn{b110}.  We note that the stochastic boundedness of $\{\langle M_{n,1} \rangle (t) : t \ge 0\}$ is trivial because it is deterministic.

\begin{lemma}{\em $($stochastic boundedness of $\{\langle M_{n,2} \rangle\}$}\label{lemFinal}
Under the assumptions in Theorems {\em \ref{th1}} and {\em \ref{thMartRep}},
the sequence $\{\langle M_{n,2} \rangle (t) \}: n \ge 1\})$ is stochastically bounded
for each $t > 0$, where
\beq
\langle M_{n,2} \rangle (t) \equiv \frac{\mu}{n} \int_{0}^{t} Q_n (s) \, d s, \quad t \ge 0 ~.
\eeqno
\end{lemma}

\paragraph{Proof.}  It suffices to apply the crude inequality in Lemma \ref{lemCrude}:
\beq
\frac{\mu}{n} \int_{0}^{t} Q_n (s) \, d s \le \frac{\mu t}{n} (Q_n (0) + A(\lambda n t))~.
\eeqno
By Lemma \ref{lemSBsums}, it suffices to show that the two sequences $\{Q_n (0)/n\}$ and $\{A(\lambda n t)/n\}$
are stochastically bounded.  By condition \eqn{a2a} and the CMT, we have the WLLN's
\beq
\frac{Q_n (0) - n}{n} \Rightarrow 0 \qasq n \ra \infty \quad \mbox{so that} \quad \frac{Q_n (0)}{n} \Rightarrow 1 \qasq n \ra \infty ~,
\eeqno
but that convergence implies stochastic boundedness; see Corollary \ref{corConvTight}.  Turning to the other term,
since $\lambda_n = n \mu$, we have
\beq
\frac{A(\lambda_n t)}{n} \ra \mu t \qasq n \ra \infty \quad \mbox{w.p.1}
\eeqno
by the SLLN for Poisson processes.
But convergence w.p.1 implies convergence in distribution, which in turn implies stochastic boundedness; again apply Corollary \ref{corConvTight}.
Hence we have completed the proof.~~~\bsq

\subsection{The Initial Conditions}\label{secInitial}
\hsp
The results in \S\S \ref{secMart}--\ref{secFluidSB} have used the finite moment condition
$E[Q_n (0)] < \infty$, but we want to establish Theorem \ref{th1} without this condition.
In particular, note that this moment condition appears prominently in Lemmas \ref{lemSatisfy} and \ref{lemYInc} and
in the resulting final martingale representations for the scaled processes, e.g., as stated in Theorems \ref{thMartRep}
and \ref{thMartRep2}.

However, for the desired Theorem \ref{th1}, we do not need to directly impose this moment condition.
We do need the assumed convergence of the scaled initial conditions in \eqn{a2a},
but we can circumvent this moment condition by defining bounded initial conditions that converge to the same limit.
For that purpose, we can work with
\beq
\hat{Q}_n (0) \equiv Q_n (0) \wedge 2n, \quad n \ge 0 ~,
\eeqno
and
\beq
\hat{X}_n (0) \equiv \frac{\hat{Q}_n (0) - n}{\sqrt{n}}, \quad n \ge 1 ~.
\eeqno
Then, for each $n \ge 1$, we use $\hat{Q}_n (0)$ and $\hat{X}_n (0)$ instead of
$Q_n(0)$ and $X_n (0)$.\break  We then obtain Theorem \ref{th1} for these modified processes:
$\hat{X}_n \Rightarrow X \qinq D \qasq\!\!\!\!\!\!\break n \ra \infty$.
However,
$P(X_n \not= \hat{X}_n) \ra 0 \qasq n \ra \infty$.
Hence, we have $X_n \Rightarrow X$ as well.

\subsection{Limit from the Fourth Martingale Representation}\label{secLimFourth}
\hsp
 In this section we state the FCLT limit for the
general $G/GI/\infty$ queue stemming from the fourth martingale representation in \S \ref{secFourth}, but
we omit proofs and simply refer to Krichagina and Puhalskii \cite{KP97}.
(We do not call the FCLT a diffusion limit because the limit is not a diffusion process.)
Even in the $M/M/\infty$ special case, the limit process has a different representation
from the representation
in Theorem \ref{th1}.  We will show that
the two representations are actually equivalent.

\begin{theorem}{\em $($FCLT from the fourth martingale representation$)$}\label{T:XnFCLT}
Let $X_n$ be defined in {\em \eqn{Xnqnq4}} and let $\hat{A}_n$ be defined in {\em \eqn{fq4}} and {\em \eqn{fq4a}}.  If
$X_n (0) \Rightarrow X(0)$ in $\RR$ and $\hat{A}_n \Rightarrow Z$ in $D$,
where $Z$ is a process with continuous sample paths and
$Z(0) = 0$,
then $X_n \Rightarrow X$ in $D$, where
\beql{X4}
X(t) \equiv F^c_0(t)X(0) + \sqrt{q(0)}W^0(F_0(t)) + M_1(t) -
M_2(t),\qquad t\geq0,
\eeq
where $W^0 \equiv \{W^0 (t): 0 \le t \le 1\}$ is a
Brownian bridge, $U \equiv \{U(t,x), t\geq 0, 0 \le x \le 1\}$ is a Kiefer
process, $X(0)$, $W^0$, $Z$ and $U$ are independent,
\beq
M_1 (t) = \int_0^t F^c(t-s)d Z(s), \qquad t\geq0,
\eeqno
 and
\beq
  M_2(t) =
\int_0^t\int_0^t\mathbf{1}(s+x\leq t) d U(a(s),F(x)), \qquad
t\geq0.
\eeqno

For $M/M/\infty$ queues, $Z = \sqrt{\mu} B$, where
$B$ is a standard Brownian motion and the limit process is
\begin{eqnarray}\label{XMM4} X(t) &=& e^{-\mu t}X(0) +
\sqrt{q(0)}W^0(1-e^{-\mu t}) +
 \sqrt{\mu} \int_0^t e^{-\mu
(t-s)}d B (s)  {} \nonumber \\
& & {} - \int_0^t\int_0^t\mathbf{1}(s+x\leq t) dU(\mu s,1-e^{-\mu
x}), \qquad t\geq0.\end{eqnarray}

\end{theorem}

\begin{remark}  A corresponding limit holds for the two-parameter processes $X_n \equiv \{X_n (t,y)\}$ in Corollary \ref{MRS4cor}
by a minor variation of the same argument.
\end{remark}

\paragraph{Connection to Theorem \ref{th1}.}

We now show that the two characterizations of the limit
$X$ in \eqn{XMM4} and \eqn{a4} are equivalent for the $M/M/\infty$ special case.  For that purpose,
express the process $X$ in \eqn{XMM4} as $X = Z_1 + Z_2 + Z_3 + Z_4$, where
$$
Z_1(t) = e^{-\mu t}X(0), \quad
Z_2(t) = \sqrt{q(0)}W^0(1-e^{-\mu t}),
$$
$$
Z_3(t) = \sqrt{\mu}\int_0^t e^{-\mu (t-s)} d B (s),
$$
$$
Z_4(t) = - \int_0^t\int_0^t\mathbf{1}(s+x\leq t) dU(\mu
s,1-e^{-\mu x}), \quad t \ge 0 ~.
$$
Clearly, $Z_1(t)$ can be written as
\beql{Z14}
Z_1(t) = -\mu \int_0^t Z_1(s)ds + X(0).
\eeq
By the solution to the linear
SDE, \S 5.6 in \cite{KS88}, we have
\beql{Z34}
Z_3(t) = -\mu \int_0^t Z_3(s)ds +  \sqrt{\mu}B (t),
\eeq
where $B$ is again standard Brownian motion.

Recall - see \S 5.6.B,
\cite{KS88} - that the Brownian bridge $W^0$ is the unique strong
solution to the one-dimensional SDE
$$
dY(t) = -\frac{Y(t)}{1-t} dt + d B_2(t), \qquad Y(0) = 0, \qquad
0\leq t\leq 1,
$$
where $B_2$ is second independent standard Brownian motion. So we can write
 \beql{W04}
 W^0(x) = -\int_0^x \frac{W^0(y)}{1-y}dy +
B_2(x)
\eeq
and it follows that
\begin{eqnarray}\label{Z24}
Z_2(t) &=&  \sqrt{q(0)} \Big(-\int_0^{1-e^{-\mu t}}
\frac{W^0(y)}{1-y}dy +
B_2(1-e^{-\mu t}) \Big)\nonumber\\
&=&  \sqrt{q(0)} \Big(-\mu \int_0^t W^0(1-e^{-\mu s}) ds + B_2(1-e^{-\mu t}) \Big) \nonumber\\
&=& -\mu \int_0^t Z_2(s) ds + \sqrt{q(0)} B_2(1-e^{-\mu t}).
\end{eqnarray}

Paralleling \eqn{W04}, the Kiefer process $U$ is related to the
Brownian sheet $W \equiv \{W(t,x): t\geq 0, 0 \leq x \leq 1\}$ by
 \beql{Usmd4}
U(t,x) = -\int_0^x\frac{U(t,y)}{1-y} dy + W(t,x), \qquad t\geq 0,
\qquad 0 \leq x \leq 1~.
\eeq
Hence, for $t\geq 0, x \geq 0$,
we have
\begin{eqnarray*}
U(\mu t, 1- e^{-\mu x}) &=& -\int_0^{1- e^{-\mu x}}\frac{U(\mu
t,y)}{1-y} dy + W(\mu t,1- e^{-\mu x}) \nonumber \\
&=& -\mu \int_0^xU(\mu t, 1- e^{-\mu y}) dy + W(\mu t,1- e^{-\mu
x}).
\end{eqnarray*}
Next, by similar reasoning, it can be shown that
\beql{Z44}
Z_4(t) = -\mu \int_0^t Z_4(s) ds +
\int_0^t \int_0^t \mathbf{1}(s+x\leq t) d W(\mu s, 1-e^{-\mu x}).
\eeq
By \eqn{Z14}, \eqn{Z34}, \eqn{Z24} and \eqn{Z44}, we have
\begin{eqnarray}\label{Z44a}
X(t) &=& X_0 - \mu \int_0^t X(s) ds  +  \sqrt{\mu} B(t)
+\sqrt{q(0)} B_2(1-e^{-\mu t})  {} \nonumber \\
& & {} + \int_0^t \int_0^t \mathbf{1}(s+x\leq t) d W(\mu s,
1-e^{-\mu x}),
\end{eqnarray}
where $B$ and $B_2$ are independent standard Brownian motions.  Let $\hat{B}$
be the sum of the last three components in \eqn{Z44a}, i.e.,
$$\hat{B}(t) \equiv \sqrt{\mu} B (t) +\sqrt{q(0)} B_2(1-e^{-\mu
t})+ \int_0^t \int_0^t \mathbf{1}(s+x\leq t) d W(\mu s, 1-e^{-\mu
x}).$$
It is evident that the process  $\hat{B}$
is a continuous Gaussian process with mean 0 and, for $s < t$,
$$
E[\hat{B}(t) - \hat{B}(s)]^2 = 2 \mu (t-s) - (1-q(0))(e^{-\mu s} -
e^{-\mu t}),
$$
and
\beq
E[\hat{B}(t) \hat{B}(s)] = 2\mu s -
(1-q(0))(1-e^{-\mu s}).
\eeqno
However, in the $M/M/\infty$ case, we have $q(0) = 1$, so that $E[\hat{B}(t) \hat{B}(s)] = 2\mu (s \wedge
t)$, which implies that the process $\{\hat{B}(t): t\geq 0\}$ is distributed the same as
$\{\sqrt{2\mu}B(t): t\geq 0\}$, where $B$ is a standard
Brownian motion. Therefore, we have shown that the $M/M/\infty$ case of Theorem \ref{T:XnFCLT} is consistent with
Theorem \ref{th1}.

\section{Other Models}\label{secOther}
\hsp
In this section we discuss how to treat other models closely related to our
initial $M/M/\infty$ model.  We consider the Erlang-$A$ model in \S \ref{secErlangA}; we also consider limits
for the waiting time there.
We consider finite waiting rooms in \S \ref{secFinite}.
Finally, we indicate how to treat general non-Poisson arrival processes
in \S \ref{secGen}.

\subsection{Erlang A Model}\label{secErlangA}
\hsp
In this section we prove the corresponding many-server heavy-traffic
limit for the $M/M/n/\infty+M$ (or Erlang-$A$ or Palm) model in the QED regime.
As before, the arrival rate is $\lambda$ and the individual service rate is $\mu$.
Now there are $n$ servers and unlimited waiting room with the FCFS service discipline.
Customer times to abandon are i.i.d. exponential random variables with a mean of $1/\theta$.
Thus individual customers waiting in queue abandon at a constant rate $\theta$.  The
Erlang-$C$ model arises when there is no abandonment, which occurs when $\theta = 0$.
The Erlang $C$ model is covered as a special case of the result below.

Let $Q(t)$ denote the number of customers in the system at time $t$, either waiting or being served.
It is well known
that the stochastic process $Q \equiv \{Q(t): t \ge 0\}$ is a birth-and-death stochastic process
with constant birth rate $\lambda_k = \lambda$ and state-dependent death rate
$\mu_k = (k \wedge n) \mu + (k - n)^{+} \theta$, $k \ge 0$, where $a \wedge b \equiv \min{\{a,b\}}$,
$a \vee b \equiv \max{\{a,b\}}$
and $(a)^{+} \equiv a \vee 0$ for real numbers $a$ and $b$.

As in Theorem \ref{th1}, the many-server heavy-traffic limit involves a
sequence of Erlang-$A$ queueing models.  As before, we let this
sequence be indexed by $n$, but now this $n$ coincides with the number of servers.
Thus now we are letting the number of servers be finite, but then letting that number
grow.  At the same time, we let the arrival rate increase.  As before, we let the arrival rate in model $n$
be $\lambda_n$, but now we stipulate that $\lambda_n$ grows with $n$.  At the same time, we hold the individual service rate $\mu$ and
abandonment rate $\theta$ fixed.  Let $\rho_n \equiv \lambda_n/n \mu$ be the \textbf{traffic intensity}
in model $n$.  We stipulate that
\beql{e1}
(1- \rho_n)\sqrt{n} \ra \beta \qasq n \ra \infty ~,
\eeq
where $\beta$ a (finite) real number.  That is equivalent to assuming that
\beql{e2}
\frac{n \mu - \lambda_n}{\sqrt{n}} \ra \beta \mu \qasq n \ra \infty ~,
\eeq
as in \eqn{QED1}.
Conditions \eqn{e1} and \eqn{e2} are known to characterize the QED many-server
heavy-traffic regime;
see Halfin and Whitt \cite{HW81} and Puhalskii and Reiman \cite{PR00}.
The many-server heavy-traffic limit theorem for the Erlang-$A$ model was proved
by Garnett, Mandelbaum and Reiman \cite{GMR02}, exploiting Stone \cite{S63}.
See Zeltyn and Mandelbaum \cite{ZM05} and Mandelbaum and Zeltyn \cite{MZ05} for
extensions and elaboration.  For related results for single-server models with
customer abandonment, see Ward and Glynn \cite{WG03a, WG03b, WG05}.

Here is the QED many-server heavy-traffic limit theorem for the $M/M/n/\infty+M$ model:

\begin{theorem} {\em $($heavy-traffic limit in $D$ for the $M/M/n+M$ model$)$}\label{th2}
Consider the sequence of $M/M/n/\infty+M$ models defined above, with the scaling in {\em \eqn{e1}}.
Let $X_n$ be as defined in {\em \eqn{a2}}.
If
$X_n (0) \Rightarrow X(0)$ in $\RR$ as $n \ra \infty$,
then
$X_n \Rightarrow X$ in $D$ as $n \ra \infty$,
where $X$ is the diffusion process with infinitesimal mean $m(x) = -\beta \mu -\mu x$ for $x < 0$
and $m(x) = -\beta \mu - \theta x$ for $x > 0$,
and infinitesimal variance $\sigma^2 (x) = 2 \mu$.  Alternatively, the limit process $X$
satisfies the stochastic integral equation
\beq
X (t) = X (0) - \mu \beta t + \sqrt{2 \mu} B (t) -  \int_{0}^{t} \left[\mu(X (s)\wedge 0) + \theta (X(s)\vee 0)\right] \, d s
\eeqno
for $t \ge 0$, where $B$ is a standard Brownian motion.
Equivalently, $X$ satisfies the SDE
\beq
d X (t) = -\mu \beta - \mu (X(t)\wedge 0) d t - \theta (X(t) \vee 0) d t + \sqrt{2 \mu} d B (t), \quad t \ge 0 ~.
\eeqno
\end{theorem}

\paragraph{Proof.}
In the rest of this section we very quickly summarize the proof;
the argument mostly differs little from what we did before.
Indeed, if we use the second martingale representation as in \S\S
\ref{secThin} and \ref{secSecond}, then there is very little
difference.  However, if we use the first martingale
representation, as in \S\S \ref{secUnit} and \ref{secFirst}, then
there is a difference, because now we want to use the optional
stopping theorem for multiparameter random time changes, as in
\S\S 2.8 and 6.2 of Ethier and Kurtz \cite{EK86}.  That approach
follows Kurtz \cite{K80}, which draws on Helms \cite{H74}.  That approach
has been applied in \S 12 of Mandelbaum and Pats \cite{MP98}.  To
illustrate this alternate approach, we use the random time change
approach here.

Just as in \S \ref{secUnit}, we can construct the stochastic process $Q$ in terms of rate-$1$ Poisson processes.
In addition to the two Poisson processes $A$ and $S$ introduced before, now we have an extra rate-$1$ Poisson process
$R$ used to generate abandonments.
Instead of \eqn{b2}, here we have representation
\begin{eqnarray}\label{ErlArep}
Q(t) & \equiv & Q(0) + A(\lambda t) - D(t) - L(t), \quad t \ge 0 ~, \nonumber \\
& = & Q(0) + A(\lambda t) - S\left(\mu\int_{0}^{t} (Q(s)\wedge n) \, d s \right) \nonumber \\
&& \quad - R \left(\theta \int_{0}^{t} (Q(s)- n)^{+} \, d s \right) ~,
\end{eqnarray}
for $t \ge 0$, where $D(t)$ is the number of departures (service completions) in the time interval $[0,t]$,
while $L(t)$ is the number of customers lost because of abandonment in the time interval $[0,t]$.
Since there are only $n$ servers, the instantaneous overall service rate at time $s$ is $\mu (Q(s)\wedge n)$,
while the instantaneous overall abandonment rate (which is only from waiting customers, not yet in service) at time $s$ is
$\theta (Q(s)- n)^{+}$.

Paralleling \eqn{b4}, we have the martingale representation
\begin{eqnarray}\label{e8}
Q(t) & = & Q(0) + M_1 (t) - M_2 (t) - M_3 (t) + \lambda t \nonumber \\
&& \quad - \mu \int_{0}^{t} (Q(s)\wedge n) \, d s
- \theta \int_{0}^{t} (Q(s)- n)^{+} \, d s
\end{eqnarray}
for $t \ge 0$, where
\begin{eqnarray}\label{e9}
M_1 (t) & \equiv & A (\lambda t) - \lambda t,  \nonumber \\
M_2 (t) & \equiv & S \left(\mu \int_{0}^{t} (Q(s)\wedge n) \, d s \right) - \mu \int_{0}^{t} (Q(s)\wedge n) \, d s ~, \nonumber \\
M_3 (t) & \equiv & R \left(\theta \int_{0}^{t} (Q(s)- n)^{+} \, d s \right) - \theta \int_{0}^{t} (Q(s)- n)^{+} \, d s
\end{eqnarray}
for $t \ge 0$ and the filtration is $\textbf{F} \equiv \{\sF_t: t \ge 0\}$ defined by
\begin{eqnarray}\label{e10}
&& \sF_t \equiv
 \sigma \left( Q(0), A(\lambda s), S \left( \mu \int_{0}^{s} (Q(u)\wedge n) \, d u \right), \right. \\
&& \quad \left. R \left( \theta \int_{0}^{s} (Q(u)- n)^{+} \, d u  \right): 0 \le s \le t \right) ~, \nonumber
\end{eqnarray}
for $t \ge 0$, augmented by including all null sets.

We now want to justify the claims in \eqn{e8}--\eqn{e10}.
Just as before, we can apply Lemmas \ref{lemInc} and \ref{lemCrude} to justify this martingale representation,
but now we need to replace Lemmas \ref{lemPQV} and \ref{lemSatisfy} by corresponding lemmas
involving the optional stopping theorem with multiparameter random time changes, as in \S\S 2.8 and 6.2 of \cite{EK86}.
We now sketch the approach:  We start with the three-parameter filtration
\begin{eqnarray}\label{mult1}
\textbf{H} & \equiv & \sH (t_1, t_2, t_3)  \\
& \equiv & \sigma\left(Q(0), A(s_1), S(s_2), R(s_3): 0 \le s_1 \le t_1, 0 \le s_2 \le t_2, 0 \le s_3 \le t_3\right) \nonumber
\end{eqnarray}
augmented by adding all null sets.
Next introduce the three nondecreasing nonnegative stochastic processes
\begin{eqnarray}
I_1 (t) & \equiv & \lambda t, \nonumber \\
I_2 (t) & \equiv & \mu \int_{0}^{t} (Q(s)\wedge n) \, ds, \nonumber \\
I_3 (t) & \equiv & \theta \int_{0}^{t} (Q(s) - n)^{+} \, ds, \quad t \ge 0 ~. \nonumber
\end{eqnarray}
Then observe that the vector $(I_1 (t), I_2 (t), I_3 (t))$ is an \textbf{H}-stopping time.
(This is essentially by the same arguments as we used in Lemma \ref{lemSatisfy}.
But here we directly gain control of the arrival process, because the event $\{I_1 (t) \le x_1\}$
coincides with the requirement that $\{\lambda t \le x_1\}$, which ensures that we always have
enough of the arrival process to construct $Q(s)$, $0 \le s \le t$.)

Moreover, since $A$, $S$ and $R$ are assumed to be independent rate-$1$ Poisson processes,
the stochastic process
\beq
\bar{M} \equiv (\bar{M}_1, \bar{M}_2, \bar{M}_3) \equiv \{(\bar{M}_1 (s_1), \bar{M}_2 (s_2), \bar{M}_3 (s_3): s_1 \ge 0, s_2 \ge 0, s_3 \ge 0\}
\eeqno
where
\begin{eqnarray}
\bar{M}_1 (t) & \equiv & A(t) - t, \nonumber \\
\bar{M}_2 (t) & \equiv & S(t) - t, \nonumber \\
\bar{M}_3 (t) & \equiv & R(t) - t, \quad t \ge 0, \nonumber
\end{eqnarray}
is an \textbf{H}-multiparameter martingale.  As a consequence of the optional stopping theorem, Theorem 8.7 on p. 87
of \cite{EK86},
\beq
(\bar{M}_1 \circ I_1, \bar{M}_2 \circ I_2, \bar{M}_3 \circ I_3) \equiv \{\bar{M}_1(I_1(t)), \bar{M}_2(I_2(t)), \bar{M}_3(I_3(t)): t \ge 0\}
\eeqno
is a martingale with respect to the filtration $\textbf{F} \equiv \{\sF_t: t \ge 0\}$ in \eqn{e10}, because
\begin{eqnarray}
&& \sH (I_1 (t), I_2 (t), I_3 (t))  \nonumber \\
& = & \sigma\left(Q(0), A(s_1), S(s_2), R(s_3): 0 \le s_1 \le I_1 (t), 0 \le s_2 \le I_2 (t), 0 \le s_3 \le I_3 (t) \right) \nonumber \\
& = & \sigma\left(Q(0), A\left( \lambda s\right), S\left(\mu \int_{0}^{s} (Q(u)\wedge n) \, d u\right) \right., \nonumber \\
&& \quad \quad \left. R\left(\mu \int_{0}^{t} (Q(s) - n)^{+} \, du\right): 0 \le s \le t  \right)
 = \sF_t \qforallq t \ge 0 ~. \nonumber
\end{eqnarray}

As in \S \ref{secFirst}, we use the crude inequality in Lemma \ref{lemCrude} to
guarantee that the moment conditions are satisfied:
\beq
E[I_j (t)] < \infty \qandq E[M_j (t)] < \infty \qforq j = 1, 2, 3 ~.
\eeqno

Just as before, we then consider the sequence of models indexed by $n$.
Just as in \eqn{b105}-\eqn{b108}, we define associated scaled processes:
\begin{eqnarray}\label{mult7}
M_{n,1} (t) & \equiv & n^{-1/2}\left[ A (\lambda_n t) - \lambda_n t\right],  \\
M_{n,2} (t) & \equiv &  n^{-1/2}\left[ S \left(\mu \int_{0}^{t} (Q_n (s)\wedge n) \, d s \right)
- \mu \int_{0}^{t} (Q_n (s) \wedge n) \, d s\right] ~, \nonumber \\
M_{n,3} (t) & \equiv &  n^{-1/2}\left[ R \left(\theta \int_{0}^{t} (Q_n (s)- n)^{+} \, d s \right)
- \theta \int_{0}^{t} (Q_n (s)- n)^{+} \, d s\right], \nonumber
\end{eqnarray}
for $t \ge 0$.

We thus obtain the following analog of Theorem \ref{thMartRep}:

\begin{theorem}{\em $($first martingale representation for the scaled processes in the Erlang A model$)$}\label{thMartRep3}
If $E[Q_n (0)] < \infty$, then the scaled processes have the martingale representation
\begin{eqnarray}\label{e11}
X_n (t) & \equiv & X_n (0) + M_{n,1} (t) - M_{n,2} (t) - M_{n,3} (t)  + \frac{(\lambda_n - \mu n) t}{\sqrt{n}}  \nonumber \\
&&  \quad  -  \int_{0}^{t} \left[\mu (X_n (s)\wedge 0)+ \theta X_n (s)^{+}\right] \, d s, \quad t \ge 0  ~,
\end{eqnarray}
where $M_{n,i}$ are the scaled martingales in {\em \eqn{mult7}}.
These processes $M_{n,i}$ are square-integrable martingales with
respect to the filtrations $\textbf{F}_n \equiv \{\sF_{n,t}: t \ge 0\}$ defined by
\begin{eqnarray}\label{e12}
\sF_{n,t} & \equiv & \sigma \left( Q_n (0), A(\lambda_n s), S\left(\mu \int_{0}^{s} (Q_n (u) \wedge n) \, d u \right) \right., \nonumber \\
&& \quad \quad \left.  R \left(\theta \int_{0}^{s} (Q_n (u)- n)^{+} \, d u  \right): 0 \le s \le t\right) ~, \quad t \ge 0 ~,
\end{eqnarray}
augmented by including all null sets.  Their associated
predictable
quadratic variations are
\begin{eqnarray}\label{e13}
\langle M_{n,1} \rangle (t) & = & \frac{\lambda_n t}{n} ~,  \nonumber \\
\langle M_{n,2} \rangle (t) & = & \frac{\mu}{n} \int_{0}^{t} (Q_n (s)\wedge n) \, d s,  \nonumber \\
\langle M_{n,3} \rangle (t) & = & \frac{\theta}{n} \int_{0}^{t} (Q_n (s)- n)^{+} \, d s, \quad t \ge 0 ~,
\end{eqnarray}
where $E[\langle M_{n,i} \rangle (t)] < \infty$ for all $i$, $t \ge 0$ and $n \ge 1$.
\end{theorem}

The representation in \eqn{e11} satisfies the Lipschitz condition in Theorem \ref{thCMT}, so that the integral representation is
again a continuous mapping.  In particular, here we have
\begin{eqnarray}\label{e14}
X_n (t) & \equiv & X_n (0) + M_{n,1} (t) - M_{n,2} (t) - M_{n,3} (t)  + \frac{(\lambda_n - \mu n) t}{\sqrt{n}} \nonumber \\
&& \quad +  \int_{0}^{t} h(X_n (s)) \, d s  ~, \quad t \ge 0 ~,
\end{eqnarray}
where $h: \RR \ra \RR$ is the function
\beql{newlip}
h(s) = -\mu(s\wedge 0) - \theta (s)^{+}, \quad s \in \RR ~,
\eeq
so that $h$ is Lipschitz as required for Theorem \ref{thCMT}:
\beq
|h(s_1) - h(s_2)| \le (\mu \vee \theta) |s_1 - s_2| \qforallq s_1, s_2 \in \RR ~.
\eeqno
Hence the proof can proceed exactly as before.  Note that we have convergence
\beq
\frac{(\lambda_n - \mu n)}{\sqrt{n}} \ra - \mu \beta  \qasq n \ra \infty
\eeqno
for the deterministic term in \eqn{e14} by virtue of the QED scaling assumption in \eqn{e1}.

The analog of Theorem \ref{thFCLTpoisson} is the corresponding FCLT for three independent rate-$1$ Poisson processes,
now including $R$ as well as $A$ and $S$.
We now have three random-time-change processes:  $\Phi_{A,n}$, $\Phi_{S,n}$ and $\Phi_{R,n}$, which here take the form:
\begin{eqnarray}\label{e17}
\Phi_{A,n} (t) & \equiv & \langle M_{n,1} \rangle (t) =  \frac{\lambda_n t}{n},  \nonumber \\
\Phi_{S,n} (t) & \equiv & \langle M_{n,2} \rangle (t) =  \frac{\mu}{n} \int_{0}^{t} (Q(s)\wedge n) \, d s,  \nonumber \\
\Phi_{R,n} (t) & \equiv & \langle M_{n,3} \rangle (t) =  \frac{\theta}{n} \int_{0}^{t} (Q(s)- n)^{+} \, d s ~,
\end{eqnarray}
drawing upon \eqn{e13}.
By the same line of reasoning as before, we obtain the deterministic limits
\begin{eqnarray}\label{e18}
\Phi_{A,n} (t) & \Rightarrow &   \mu e,  \quad
\Phi_{S,n} (t)  \Rightarrow  \mu e  \qandq
\Phi_{R,n} (t)  \Rightarrow  \eta ~,
\end{eqnarray}
where, as before, $e$ is the identity map $e(t) \equiv t, \quad t \ge 0$, and
$\eta$ is the zero function $\eta (t) \equiv 0$, $t \ge 0$.
In particular, we again get the sequence $\{X_n: n \ge 1\}$ stochastically bounded in $D$ by
the reasoning in \S\S \ref{secSB} and \ref{secSBcomp}.  Then, by Lemma \ref{lemSBfluid} and \S \ref{secFluidSB}, we
get the FWLLN corresponding to Lemma \ref{lemFluid2}.  Finally, from Lemma \ref{lemFluid2},
we can prove \eqn{e18}, just as we proved Lemma \ref{lemFluid1}, using analogs of the continuous map in \eqn{c11c}
for the random time changes in \eqn{e17}.  The new functions for applications of the CMT are
\beq
h_1 (x) (t) \equiv \mu \int_{0}^{t} (x(s)\wedge 1) \, ds  \qandq
h_2 (x) (t) \equiv \theta \int_{0}^{t} (x(s)- 1)^{+} \, ds
\eeqno
for $t \ge 0$.

Paralleling \eqn{c12}, here we have
\beq
(M_{A,n}, \Phi_{A,n}, M_{S,n},\Phi_{S,n}, M_{R,n},\Phi_{R,n}) \Rightarrow (B_1, \mu e, B_2, \mu e, B_3, \eta) \qinq D^6
\eeqno
as $n \ra \infty$.  Hence we can apply the CMT with composition just as before.  Paralleling \eqn{c13}, here we obtain first
\beq
(M_{A,n}\circ \Phi_{A,n}, M_{S,n} \circ \Phi_{S,n}, M_{R,n} \circ \Phi_{R,n}) \Rightarrow (B_1 \circ \mu e, B_2 \circ \mu e, B_3 \circ \eta)
\eeqno
in $D^3$ as $n \ra \infty$, and then
\beq
(M_{A,n}\circ \Phi_{A,n}- M_{S,n} \circ \Phi_{S,n}- M_{R,n} \circ \Phi_{R,n}) \Rightarrow (B_1 \circ \mu e- B_2 \circ \mu e- B_3 \circ \eta)
\eeqno
in $D$.  However,
\beql{e21}
B_1 \circ \mu e- B_2 \circ \mu e- B_3 \circ \eta \deq B_1 \circ \mu e- B_2 \circ \mu e- \eta  \deq \sqrt{2 \mu} B ~.
\eeq
Finally, the CMT with the integral representation \eqn{e11} and Theorem \ref{thCMT} completes the proof.
Note that the
limiting Brownian motion associated with $R$ does not appear, because $\Phi_{R,n}$ is asymptotically negligible.
That is why the infinitesimal variance is the same as before.~~~\bsq

\subsection{Finite Waiting Rooms}\label{secFinite}
\hsp
We can also obtain stochastic-process limits for the number of customers in the system
in associated $M/M/n/0$ (Erlang-$B$), $M/M/n/m_n$ and $M/M/n/m_n+M$ models, which have finite waiting rooms.
For the Erlang-$B$ model, there is no waiting room at all; for the other models there is a waiting room on size $m_n$
in model $n$, where $m_n$ is allowed to grow with $n$ so that
$m_n/\sqrt{n} \ra \kappa \ge 0$ as $n \ra \infty$
as in \eqn{QED2}.  The QED many-server heavy-traffic limit was stated as Theorem \ref{th3}.

The proof can be much the same as in \S \ref{secErlangA}.
The idea is to introduce the finite waiting room via a reflection map, as in \S\S 3.5, 5.2, 13.5, 14.2, 14.3 and 14.8 of Whitt \cite{W02},
corresponding to an upper barrier at $\kappa$, but the reflection map here is more complicated than for single-server queues
and networks of such queues, because it is not
applied to a free process.
We use an extension of Theorem \ref{thCMT} constructing a mapping from $D \times \RR$ into $D^2$, taking model data
into the content function and the upper-barrier regulator function.

\begin{theorem}{$($a continuous integral representation with reflection$)$}\label{thCMT2}
Consider the modified integral representation
\beql{lip11}
x(t) =   b + y (t) +  \int_{0}^{t} h(x (s)) \, ds - u(t) ~, \quad t \ge 0 ~,
\eeq
where $x(t) \le \kappa$,  $h: \RR \ra \RR$ satisfies $h(0) = 0$ and is a Lipschitz function as in {\em \eqn{lip2}}, and
$u$ is a nondecreasing nonnegative function in $D$ such that {\em \eqn{lip11}} holds and
\beq
\int_{0}^{\infty} 1_{\{x (t) < \kappa\}} \, d u(t) = 0 ~.
\eeqno
The modified integral representation in {\em \eqn{lip11}} has a unique solution $(x, u)$,
so that it constitutes a bonafide function $(f_1, f_2): D \times \RR \ra D \times D$
mapping $(y,b)$ into $x \equiv f_1 (y,b)$ and $u \equiv f_2 (y,b)$.
In addition, the function $(f_1, f_2)$ is continuous provided that the product topology is used for product spaces
and the function space $D$ $($in both the domain and range$)$ is endowed with
either: $(i)$ the topology of uniform
convergence over bounded intervals or $(ii)$ the Skorohod $J_1$ topology.
Moreover, if $y$ is continuous, then so are $x$ and $u$.
\end{theorem}

\paragraph{Proof.}
We only show the key step, for which we follow the argument in \S 3 of Mandelbaum and Pats \cite{MP98} and
\S 4 of Reed and Ward \cite{RW07}; see these sources for additional details and references.
The idea is to combine classical results for the conventional one-dimensional reflection map,
as in \S\S 5.2 and 13.5 of \cite{W02} with a modification of Theorem 4.1.
Let $(\phi_{\kappa}, \psi_{\kappa})$ be the one-sided reflection map with upper barrier at $\kappa$,
so that $\phi_{\kappa} (y) = y - \psi_{\kappa} (y)$, with $\phi_{\kappa} (y)$ being the content function
and $\psi_{\kappa} (y)$ being the nondecreasing regulator function; see \S\S 5.2 and 13.5 of \cite{W02}.  We observe that
the map in \eqn{lip11} can be expressed as $x = \phi_{\kappa} (w)$ and $u = \psi_{\kappa} (w)$, where
\beql{reed1}
w(t) \equiv \xi(b, y) (t) \equiv b + y (t) + \int_{0}^{t} h(\phi_{\kappa} (w(s))) \, ds, \quad t \ge 0 ~.
\eeq
This lets us represent the desired map as the composition of the maps $(\phi_{\kappa}, \psi_{\kappa})$ and $\xi$.
The argument to treat $\xi$ is essentially the same as in the proof of Theorem \ref{thCMT}, but
we need to make a slight adjustment; we could apply it directly if we had $h:D \ra D$ in Theorem \ref{thCMT}.
  Recall that $\phi_{\kappa}$ is Lipschitz continuous on $D([0,t]$
for each $t$ with the uniform norm, $\| \phi_{\kappa} (y_1) - \phi_{\kappa} (y_2)\|_{t} \le 2 \| y_1 - y_2\|_{t}$,
with modulus $2$ independent of $t$.
Hence, paralleling \eqn{c4}, we have
\beq
\| w_1 - w_2\|_{t} \le |b_1 - b_2| + \| y_1 - y_2\|_{t} + 2 c \int_{0}^{t} \| w_1 - w_2\|_{s} \, ds
\eeqno
for each $t > 0$.  Hence we can apply Gronwall's inequality in Lemma \ref{lemGronwall}
 to establish (Lipschitz) continuity of the map $\xi$ on $D([0,T])\times \RR$.
Combining this with the known (Lipschitz) continuity of the reflection map $(\phi_{\kappa}, \psi_{\kappa})$, we have the desired continuity
for the overall map in the uniform topology. We can extend to the $J_1$ topology as in the proof of Theorem \ref{thCMT}.
~~~\bsq

Now that we understand how we are treating the finite waiting rooms, the QED many-server heavy-traffic limit theorem
is as stated in Theorem \ref{th3}.
This modification alters the limiting diffusion process in Theorem \ref{th2} only by the addition of a reflecting upper barrier at
$\kappa$ for the sequence of models with waiting rooms of size $m_n$,
where $\kappa = 0$ for the Erlang-$B$ model.
When $\kappa = 0$, $X$ is a reflected OU (ROU) process.  Properties of the ROU process are contained in Ward and Glynn \cite{WG03b}.
Proofs for the two cases $\kappa > 0$ and $\kappa = 0$ by other methods are contained in
\S 4.5 of Whitt \cite{W05} and Theorem 4.1 of Srikant and Whitt \cite{SW96}.
General references on reflection maps are Lions and Sznitman \cite{LS84} and Dupuis and Ishii \cite{DI91}.

\paragraph{Proof of Theorem \ref{th3}.}  We briefly sketch the proof.
Instead of \eqn{ErlArep}, here we have representation

\begin{eqnarray}\label{RoomRep1}
Q_n (t) & \equiv & Q_n (0) + A(\lambda_n t) - D_n (t) - L_n (t) - U_n (t), \quad t \ge 0 ~, \nonumber \\
& = & Q_n (0) + A(\lambda_n t) - S\left(\mu\int_{0}^{t} (Q_n (s)\wedge n) \, d s \right) \nonumber \\
&& \quad - R \left(\theta \int_{0}^{t} (Q_n (s)- n)^{+} \, d s \right) - U_n (t) ~,
\end{eqnarray}
for $t \ge 0$, where $U_n (t)$ is
the number of arrivals in the time interval $[0,t]$
when the system is full in model $n$,
i.e., when $Q_n (t) = n + m_n$.  In particular,
\beql{RoomRep2}
U_n (t) \equiv \int_{0}^{t} 1_{\{Q_n (s) = n + m_n\}} \, d A(\lambda_n s), \quad t \ge 0 ~.
\eeq
To connect to Theorem \ref{thCMT2}, it is significant that $U_n$ can also be represented as
the unique nondecreasing nonnegative process such that $Q_n (t) \le n + m_n$, \eqn{RoomRep1} holds and
\beql{RoomRep3}
\int_{0}^{\infty} 1_{\{Q_n (t) < \kappa\}} \, d U_n (t) = 0 ~.
\eeq

We now construct a martingale representation, just as in \eqn{e8}--\eqn{mult7}.  The following is the
natural extension of Theorem \ref{thMartRep3}:

\begin{theorem}{\em $($first martingale representation for the scaled processes in the $M/M/n/m_n+M$ model$)$}\label{thMartRep4}
If $\kappa < \infty$, then the scaled processes have the martingale representation
\begin{eqnarray}\label{ee11}
X_n (t) & \equiv &  X_n (0) + M_{n,1} (t) - M_{n,2} (t) - M_{n,3} (t)  + \frac{(\lambda_n - \mu n) t}{\sqrt{n}} \nonumber \\
&&  -  \int_{0}^{t} \left[\mu (X_n (s)\wedge 0)+ \theta X_n (s)^{+}\right] \, d s - V_n (t), \quad t \ge 0  ~,
\end{eqnarray}
$M_{n,i}$ are the scaled martingales in {\em \eqn{mult7}} and
\beql{ee11a}
V_n (t) \equiv  \frac{U_n (t)}{\sqrt{n}}, \quad t \ge 0 ~,
\eeq
for $U_n$ in {\em \eqn{RoomRep1}}--{\em \eqn{RoomRep3}}.
The scaled processes $M_{n,i}$ are square-integrable martingales with
respect to the filtrations $\textbf{F}_n \equiv \{\sF_{n,t}: t \ge 0\}$ defined by
\begin{eqnarray}\label{ee12}
\sF_{n,t} & \equiv & \sigma \left( Q_n (0), A(\lambda_n s), S\left(\mu \int_{0}^{s} (Q_n (u) \wedge n) \, d u \right) \right., \nonumber \\
&&  \quad \left. R \left(\theta \int_{0}^{s} (Q_n (u)- n)^{+} \, d u  \right): 0 \le s \le t\right) ~, \quad t \ge 0 ~,
\end{eqnarray}
augmented by including all null sets.  Their associated
predictable
quadratic variations are
\begin{eqnarray}\label{ee13}
\langle M_{n,1} \rangle (t) & = & \frac{\lambda_n t}{n} ~,  \nonumber \\
\langle M_{n,2} \rangle (t) & = & \frac{\mu}{n} \int_{0}^{t} (Q_n (s)\wedge n) \, d s,  \nonumber \\
\langle M_{n,3} \rangle (t) & = & \frac{\theta}{n} \int_{0}^{t} (Q_n (s)- n)^{+} \, d s, \quad t \ge 0 ~,
\end{eqnarray}
where $E[\langle M_{n,i} \rangle (t)] < \infty$ for all $i$, $t \ge 0$ and $n \ge 1$.
\end{theorem}

By combining Theorems \ref{thCMT2} and \ref{thMartRep4}, we obtain the joint convergence
\beq
(X_n, V_n) \Rightarrow (X, U) \qinq D^2 \qasq n \ra \infty ~,
\eeqno
for $X_n$ and $V_n$ in \eqn{ee11} and \eqn{ee11a}, where the vector
$(X, U)$ is characterized by \eqn{ee6} and \eqn{ee6a}. That implies Theorem \ref{th3} stated in \S \ref{secIntro}.

\begin{remark}{\em $($correction$)$}
{\em
The argument in this section follows Whitt \cite{W05}, but provides more detail.
We note that the upper-barrier regulator processes are incorrectly expressed
in formulas (5.2) and (5.8) of \cite{W05}.
}
\end{remark}

\subsection{General Non-Markovian Arrival Processes}\label{secGen}
\hsp
In this section, following \S 5 of Whitt \cite{W05}, we show how to extend the many-server heavy-traffic limit from $M/M/n/m_n+M$ models
to $G/M/n/m_n+M$ models, where the arrival processes are allowed to be general stochastic point processes
satisfying a FCLT.  They could be renewal processes ($GI$) or even more general arrival processes.  The limit of
the arrival-process FCLT need not have continuous sample paths.  (As noted at the end of \S \ref{secSLLN}, this
separate argument is not needed if we do not use martingales.)

Let $\bar{A}_n \equiv \{\bar{A}_n (t): t \ge 0\}$ be the general arrival process in model $n$ and let
\beql{gen2}
A_n (t) \equiv \frac{\bar{A}_n (t) - \lambda_n t}{\sqrt{n}}, \quad t \ge 0 ~,
\eeq
be the associated scaled arrival process.  We assume that
\beql{gen1}
A_n \Rightarrow A \qinq D \qasq n \ra \infty ~.
\eeq
We also assume that, conditional on the entire arrival process, model $n$ evolves as the Markovian queueing process
with i.i.d. exponential service times and i.i.d. exponential times until abandonment.

Thus, instead of Theorem \ref{thMartRep4}, we have

\begin{theorem}{\em $($first martingale representation for the scaled processes in the $G/M/n/m_n+M$ model$)$}\label{thMartRep5}
Consider the family of $G/M/n/m_n+M$ models defined above, evolving as a Markovian queue conditional on the arrival process.
If $m_n < \infty$, then the scaled processes have the martingale representation
\begin{eqnarray}\label{eee11}
X_n (t) & \equiv & X_n (0) + A_{n} (t) - M_{n,2} (t) - M_{n,3} (t)  + \frac{(\lambda_n - \mu n) t}{\sqrt{n}} \nonumber \\
&& \quad  -  \int_{0}^{t} \left[\mu (X_n (s)\wedge 0)+ \theta X_n (s)^{+}\right] \, d s - V_n (t), \quad t \ge 0  ~,
\end{eqnarray}
where $A_n$ is the scaled arrival process in {\em \eqn{gen2}}, $M_{n,i}$ are the scaled martingales in {\em \eqn{mult7}} and
$V_n (t) \equiv  U_n (t)/\sqrt{n}$, $t \ge 0$,
for $U_n$ in {\em \eqn{RoomRep1}}-{\em \eqn{RoomRep3}}.
The scaled processes $M_{n,i}$ are square-integrable martingales with
respect to the filtrations $\textbf{F}_n \equiv \{\sF_{n,t}: t \ge 0\}$ defined by
\begin{eqnarray}\label{eee12}
\sF_{n,t} & \equiv & \sigma \left( Q_n (0), \{A_n (u): u \ge 0\}, S\left(\mu \int_{0}^{s} (Q_n (u) \wedge n) \, d u \right) \right., \nonumber \\
&& \quad \left. R \left(\theta \int_{0}^{s} (Q_n (u)- n)^{+} \, d u  \right): 0 \le s \le t\right) ~, \quad t \ge 0 ~,
\end{eqnarray}
augmented by including all null sets.  The associated
predictable
quadratic variations of the two martingales are
\begin{eqnarray}\label{eee13}
\langle M_{n,2} \rangle (t) & = & \frac{\mu}{n} \int_{0}^{t} (Q_n (s)\wedge n) \, d s,  \nonumber \\
\langle M_{n,3} \rangle (t) & = & \frac{\theta}{n} \int_{0}^{t} (Q_n (s)- n)^{+} \, d s, \quad t \ge 0 ~,
\end{eqnarray}
where $E[\langle M_{n,i} \rangle (t)] < \infty$ for all $i$, $t \ge 0$ and $n \ge 1$.
\end{theorem}

Here is the corresponding theorem for the $G/M/n/m_n + M$ model.

\begin{theorem} {\em $($heavy-traffic limit in $D$ for the $G/M/n/m_n+M$ model$)$}\label{th4}
Consider the sequence of $G/M/n/m_n+M$ models defined above, with the scaling in {\em \eqn{gen2}}, {\em \eqn{QED1}} and {\em \eqn{QED2}}.
Let $X_n$ be as defined in {\em \eqn{a2}}.
If
\beq
X_n (0) \Rightarrow X(0) \qinq \RR \qandq  A_n \Rightarrow A \qinq D \qasq n \ra \infty ~,
\eeqno
then
$X_n \Rightarrow X \qinq D \qasq n \ra \infty$,
where the limit process $X$
satisfies the stochastic integral equation
\begin{eqnarray}\label{eee6}
X (t) & = & X (0) + A(t) - \beta \mu t  +\sqrt{\mu} B (t) \nonumber \\
&& \quad -  \int_{0}^{t} \left[\mu(X (s)\wedge 0) + \theta (X(s)\vee 0)\right] \, d s - U (t)
\end{eqnarray}
for $t \ge 0$ with $B \equiv \{B(t): t \ge 0\}$ being a standard Brownian motion and $U$ being the unique nondecreasing nonnegative process in $D$
such that $X(t) \le \kappa$ for all $t$, {\em \eqn{eee6}} holds and
\beq
\int_{0}^{\infty} 1_{\{X (t) < \kappa\}} \, d U(t) = 0 ~.
\eeqno
\end{theorem}

\paragraph{Proof.} We start with the FCLT for the arrival process assumed in \eqn{gen1} and \eqn{gen2}.
We condition on possible realizations of the arrival process:
For each $n \ge 1$, let $\zeta_n$ be a possible realization of the scaled arrival process $A_n$
and let $\zeta$ be a possible realization of the limit process $A$.  Let $X^{\zeta_n}_{n}$ be the scaled process $X_n$
conditional on $A_n = \zeta_n$, and let $X^{\zeta}$ be the limit process $X$ conditional on $A = \zeta$.

Since $X_n$ and $A_n$ are random elements of $D$, these quantities $X^{\zeta_n}_{n}$ and $X^{\zeta}$ are well
defined via regular conditional probabilities; e.g., see \S 8 of Chapter V, pp 146-150, of Parthasarathy \cite{P67}.
In particular, we can regard $P(X_n \in \cdot | A_n = z)$ as a probability
measure on $D$ for each $z \in D$ and we can regard $P(X_n \in B | A_n = z)$ as a measurable function of $z$ in $D$
for each Borel set $B$ in $D$, where
\beq
P(X_n \in B) = \int_{B} P(X_n \in B | A_n = z) \, d P (A_n = z) ~.
\eeqno
And similarly for the pair $(X, A)$.

A minor modification of the previous proof of Theorem \ref{th3} establishes that
\beq
X^{\zeta_n}_{n} \Rightarrow  X^{\zeta} \qinq D \quad \mbox{whenever} \quad \zeta_n \ra \zeta \qinq D ~;
\eeqno
i.e., for each continuous bounded real-valued function $f$ on $D$,
\beql{gen12}
E[f(X^{\zeta_n}_n)] \ra E[f(X^{\zeta})] \qasq n \ra \infty
\eeq
whenever $\zeta_n \ra \zeta$ in $D$.
Now fix a continuous bounded real-valued function $f$ and let
\beql{gen12a}
h_n (\zeta_n) \equiv E[f(X^{\zeta_n}_n)] \qandq h (\zeta) \equiv E[f(X^{\zeta})] ~.
\eeq
Since we have regular conditional probabilities, we can regard the functions $h_n$ and $h$
as measurable functions from $D$ to $\RR$ (depending on $f$).

We are now ready to apply the generalized continuous mapping theorem, Theorem 3.4.4 of \cite{W02}.
Since $h_n$ and $h$ are measurable functions
such that $h_n (\zeta_n) \ra h(\zeta)$ whenever $\zeta_n \ra \zeta$
and since $A_n \Rightarrow A$ in $D$, we have
$h_n (A_n) \Rightarrow h(A)$ as $n \ra \infty$.
Since the function $f$ used in \eqn{gen12} and \eqn{gen12a} is bounded, these random variables are bounded.  Hence convergence
in distribution implies convergence of moments.  Hence,
 for that function $f$, we have
\beq
E[f(X_n)] \equiv E[h_n (A_n)] \ra E[h(A)] \equiv E[f(X)] \qasq n \ra \infty ~.
\eeqno
Since this convergence holds for all continuous bounded real-valued functions $f$ on $D$, we have shown
that $X_n \Rightarrow X$, as claimed.~~~\bsq

 \section{The Martingale FCLT}\label{secMartFCLT}
 \hsp
 We now turn to the martingale FCLT.  For our queueing stochastic-process limits, it is of interest
 because it provides
 one way to prove the FCLT for a Poisson process in Theorem \ref{thFCLTpoisson} and
because we can base our entire proof of Theorem \ref{th1} on the martingale FCLT.  However,
the gain in the proof of Theorem \ref{th1} is not so great.

We now state a version of the martingale FCLT for a sequence of local martingales $\{M_n: n \ge 1\}$
in $D^k$, based on Theorem 7.1 on p. 339 of Ethier and Kurtz \cite{EK86}, hereafter referred to as EK.
Another important reference is Jacod and Shiryayev \cite{JS87}, hereafter referred to as JS.
See Section VIII.3 of JS for related results; see other sections of JS for generalizations.

We will state a special case of Theorem 7.1 of EK in which the limit process
is multi-dimensional Brownian motion.  However, the framework always produces limits
with continuous sample paths and independent Gaussian increments.
Most applications involve convergence to Brownian motion.  Other situations are covered by JS, from which we see that
proving convergence to discontinuous processes is more complicated.

The key part of each condition below is the convergence of the
quadratic covariation processes.  Condition (i) involves the optional quadratic-covariation
(square-bracket) processes $\left[M_{n,i}, M_{n,j}\right]$, while condition (ii) involves
the predictable quadratic-covariation
(angle-bracket) processes $\langle M_{n,i}, M_{n,j} \rangle$.  Recall from \S \ref{secQuad} that
the square-bracket process is more general, being well defined for any local martingale (and thus any martingale),
whereas the associated angle-bracket process is well defined only for any locally square-integrable martingale
(and thus any square-integrable martingale).

Thus the key conditions below are the assumed convergence of the quadratic-variation processes in
conditions \eqn{ek2} and \eqn{ek6}.
The other conditions \eqn{ek1}, \eqn{ek4} and \eqn{ek6} are technical regularity conditions.
There is some variation in the literature concerning the extra technical regularity conditions; e.g.,
see Rebolledo \cite{R80} and JS.

Let $J$ be the \textbf{maximum-jump function}, defined for any $x \in D$ and $T > 0$ by
 \beql{mod3}
J(x,T) \equiv \sup{\{|x(t) - x(t-)|: 0 < t \le T \}} ~.
\eeq

\begin{theorem}{\em $($multidimensional martingale FCLT$)$}\label{thMart}
For $n \ge 1$, let $M_n \equiv (M_{n,1}, \ldots , M_{n,k})$ be a local martingale in $D^k$
with respect to a filtration $\textbf{F}_n \equiv \{\sF_{n,t}: t \ge 0\}$
satisfying $M_n (0) = (0, \ldots , 0)$.  Let $C \equiv (c_{i,j})$ be
a $k \times k$ covariance matrix, i.e., a nonnegative-definite symmetric matrix of real numbers.

\vspace{0.1in}
  \textbf{Assume that one of the following two conditions holds:}

\vspace{0.1in}
\textbf{$(i)$}  The expected value of the maximum jump in $M_n$ is asymptotically negligible; i.e., for each $T > 0$,
\beql{ek1}
\lim_{n \ra \infty}{\left\{ E\left[ J(M_n,T)\right]\right\}} = 0
\eeq
and, for each pair $(i,j)$ with $1 \le i \le k$ and $1 \le j \le k$, and each $t > 0$,
\beql{ek2}
\left[M_{n,i}, M_{n,j}\right] (t) \Rightarrow c_{i,j} t \qinq \RR \qasq n \ra \infty ~.
\eeq

\vspace{0.1in}
\textbf{$(ii)$}  The local martingale $M_n$ is locally square-integrable, so that the predictable quadratic-covariation processes
$\langle M_{n,i}, M_{n,j} \rangle$ can be defined.
The expected value of the maximum jump in $\langle M_{n,i}, M_{n,j} \rangle$
and the maximum squared jump of $M_n$ are asymptotically negligible;
i.e., for each $T > 0$
and $(i,j)$ with $1 \le i \le k$ and $1 \le j \le k$,
\beql{ek4}
\lim_{n \ra \infty}{\left\{ E\left[ J \left( \langle M_{n,i}, M_{n,j} \rangle, T \right) \right]\right\}}  = 0 ~,
\eeq
\beql{ek5}
\lim_{n \ra \infty}{\left\{ E\left[ J \left( M_n, T \right)^2 \right]\right\}} = 0 ~,
\eeq
and
\beql{ek6}
\langle M_{n,i}, M_{n,j} \rangle (t)  \Rightarrow c_{i,j} t \qinq \RR \qasq n \ra \infty
\eeq
for each $t > 0$ and for each $(i,j)$.

\vspace{0.1in}
  \textbf{Conclusion:}

\vspace{0.1in}
If indeed one of the the conditions $(i)$ or $(ii)$ above holds, then
\beq
M_n \Rightarrow  M \qinq D^k \qasq n \ra \infty ~,
\eeqno
where $M$ is a \textbf{$k$-dimensional $(0,C)$-Brownian motion}, having mean vector and covariance matrix
\beq
E[M(t)] = (0, \ldots , 0) \qandq E[M(t)M(t)^{tr}] = C t, \quad t \ge 0 ~,
\eeqno
where, for a matrix $A$, $A^{tr}$ is the transpose.
\end{theorem}

Of course, a common simple case arises when $C$ is a diagonal matrix; then the $k$ component marginal one-dimensional
Brownian motions are independent. When $C = I$, the identity matrix, $M$ is a standard
$k$-dimensional Brownian motion, with independent one-dimensional standard Brownian motions as marginals.

At a high level, Theorem \ref{thMart} says that, under regularity conditions, convergence of martingales in $D$ is implied
by convergence of the associated quadratic covariation processes.  At first glance, the result seems even stronger,
because we need convergence of only the one-dimensional quadratic covariation processes for a single time argument.
However, that is misleading, because the stronger weak convergence of these quadratic covariation processes
in $D^{k^2}$ is actually equivalent to the weaker required convergence in $\RR$ for each $t, i, j$ in conditions
\eqn{ek2} and \eqn{ek6}; see \cite{W07}.

\section{Applications of the Martingale FCLT}\label{secApp}
\hsp

In this section we make two applications of the preceding martingale FCLT.
First, we apply it to prove the FCLT for the scaled Poisson process,
Theorem \ref{thFCLTpoisson}.  Then we apply it to provide a third proof of Theorem \ref{th1}.
In the same way we could obtain alternate proofs of Theorems \ref{th3} and \ref{th2}.

\subsection{Proof of the Poisson FCLT}\label{secPoissonPf}

We now apply the martingale FCLT to prove the Poisson FCLT in Theorem \ref{thFCLTpoisson}.  To do so,
it suffices to consider the one-dimensional version in $D$, since the Poisson processes are mutually independent.
Let the martingales $M_n \equiv M_{A,n}$ be as defined in \eqn{c8}, i.e.,
\beq
M_n \equiv M_{A,n} (t) \equiv \frac{A(nt) - nt}{\sqrt{n}}, \quad t \ge 0 ~,
\eeqno
with their internal filtrations.
The limits in \eqn{ek1} and \eqn{ek5} hold, because
\beq
J(M_n, T)  \le \frac{1}{\sqrt{n}}
 \qandq
J(M_n^2, T)  \le \frac{1}{n}, \quad n \ge 1 ~.
\eeqno
For each $n \ge 1$, $M_n$ is square integrable, the optional quadratic variation process is
\beq
[M_n] (t)  \equiv [M_n, M_n] (t) = \frac{A(nt)}{n}, \quad t \ge 0  ~,
\eeqno
and the predictable quadratic variation process is
\beq
\langle M_n \rangle (t) \equiv \langle M_n, M_n  \rangle (t) = \frac{nt}{n} \equiv t, \quad t \ge 0 ~.
\eeqno
Hence both \eqn{ek4} and \eqn{ek6} hold trivially.
By the SLLN for a Poisson process,
$A(nt)/n \ra t$ w.p.1 as $n \ra \infty$
for each $t > 0$.  Hence both conditions (i) and (ii) in Theorem \ref{thMart} are satisfied,
with $C = c_{1,1} = 1$.~~~\bsq

\subsection[Third Proof of Theorem 1.1]{Completing the Proof of Theorem \ref{th1}}\label{secProofSecond}

The bulk of this paper has consisted of a proof of Theorem \ref{th1} based on
the first martingale representation in Theorem \ref{thMartRep}, which in turn is
based on the representation of the service-completion counting process as
a random time change of a rate-$1$ Poisson process, as in \eqn{b2}.
A second proof in \S \ref{secSLLN} established the fluid limit directly.

In this subsection we present a third proof of Theorem \ref{th1} based
on the second martingale representation in Theorem \ref{thMartRep2},
which in turn is based on a random thinning of rate-$1$ Poisson processes, as in \eqn{b4}.
This third proof also applies to the third martingale representation
in \S \ref{secThird}, which is based on constructing martingales for counting processes
associated with the birth-and-death process $\{Q(t): t \ge 0\}$ via
its infinitesimal generator.

Starting with the second martingale representation in Theorem \ref{thMartRep2} (or the third martingale representation
in Subsection \ref{secThird}), we cannot rely on the Poisson FCLT to obtain the required stochastic process limit
\beql{w1}
(M_{n,1}, M_{n,2}) \Rightarrow (\sqrt{\mu} B_1, \sqrt{\mu} B_2) \qinq D^2 \qasq n \ra \infty ~,
\eeq
in \eqn{c1}.  However, we can apply the martingale FCLT for this purpose, and we show how to do that now.

As in \S \ref{secPoissonPf},
we can apply either condition (i) or (ii) in Theorem \ref{thMart}, but it is easier to apply (ii), so we will.
The required argument looks more complicated because we have to establish the two-dimensional convergence
in \eqn{w1} in $D^2$ because the scaled martingales $M_{n,1}$ and $M_{n,2}$ are not independent.
Fortunately, however, they are orthogonal, by virtue of the following lemma.  That still means that we need to establish the two-dimensional limit
in \eqn{w1}, but it is not difficult to do so.

We say that two locally square-integrable martingales with respect to the filtration $\textbf{F}$, $M_1$ and $M_2$,
are \textbf{orthogonal} if the process $M_1 M_2$ is a local martingale with $M_1 (0) M_2 (0) = 0$.
Since $M_1 M_2 - \langle M_1, M_2 \rangle$ is a local martingale, orthogonality implies that
$\langle M_1, M_2 \rangle (t) = 0$ for all $t$.

\begin{lemma}{\em $($orthogonality of stochastic integrals with respect to orthogonal martingales$)$}\label{lemOrthog}
 Suppose that $M_1$ and $M_2$ are locally square-integrable martingales with respect to the filtration $\textbf{F}$,
 where $M_1 (0) = M_2 (0) = 0$,
 while $C_1$ and $C_2$ are locally-bounded $\textbf{F}$-predictable processes.  If $M_1$ and $M_2$ are orthogonal,
 $($which is implied by independence$)$, then the
 stochastic integrals $\int_{0}^{t} C_1 (s) \, dM_1 (s)$ and $\int_{0}^{t} C_2 (s) \, dM_2 (s)$ are orthogonal, which implies that
 \beq
 \left\langle \int C_1 (s) \, dM_1 (s) , \int C_2 (s) \, dM_2 (s)\right\rangle (t) = 0, \quad t \ge 0 ~.
 \eeqno
 and
 \beq
 \left[ \int C_1 (s) \, dM_1 (s) , \int C_2 (s) \, dM_2 (s)\right] (t) = 0, \quad t \ge 0 ~.
 \eeqno
\end{lemma}

Lemma \ref{lemOrthog} follows from the following representation
for the quadratic covariation of the stochastic integrals; see \S 5.9 of van der Vaart \cite{V06}.

\begin{lemma}{\em $($Quadratic covariation of stochastic integrals with respect to mar\-tingales$)$}\label{lemQVstochInt}
 Suppose that $M_1$ and $M_2$ are locally square-integrable martingales with respect to the filtration $\textbf{F}$,
 while $C_1$ and $C_2$ are locally-bounded $\textbf{F}$-predictable processes.
Then
\begin{eqnarray}
 && \left\{\int_{0}^{t} C_1 (s) \, dM_1 (s) \int_{0}^{t} C_2 (s) \, dM_2 (s) \right.\nonumber \\
&& \quad \quad - \left. \left[\int C_1 (s) \, dM_1 (s) , \int C_2 (s) \, dM_2 (s)\right] (t): t \ge 0\right\} \nonumber
\end{eqnarray}
and
\begin{eqnarray}
&& \left\{\int_{0}^{t} C_1 (s) \, dM_1 (s) \int_{0}^{t} C_2 (s) \, dM_2 (s) \right. \nonumber \\
&& \quad \quad  -  \left. \left\langle \int C_1 (s) \, dM_1 (s) , \int C_2 (s) \, dM_2 (s)\right\rangle (t): t \ge 0\right\} \nonumber
\end{eqnarray}
 are
local $\textbf{F}$-martingales, where
the quadratic covariation processes are
\beq
 \left[\int C_1 (s) \, dM_1 (s) , \int C_2 (s) \, dM_2 (s)\right] (t) = \int_{0}^{t} C_1 (s) C_2 (s) \, d [M_1, M_2] (s)
\eeqno
and
\beq
 \left\langle \int C_1 (s) \, dM_1 (s) , \int C_2 (s) \, dM_2 (s)\right\rangle (t)
= \int_{0}^{t} C_1 (s) C_2 (s) \, d \langle M_1, M_2 \rangle (s), \quad t \ge 0 ~.
\eeqno
\end{lemma}

As a consequence of the orthogonality provided by Lemma \ref{lemOrthog}, we have $[M_{n,1}, M_{n,2}] (t) = 0$ and
$\langle M_{n,1}, M_{n,2} \rangle (t) = 0$ for all $t$ and $n$
for the martingales in \eqn{w1}, which in turn come from Theorem \ref{thMartRep2}.
Thus the orthogonality trivially implies that
\beq
[M_{n,1}, M_{n,2}] \Rightarrow 0 \qandq \langle M_{n,1}, M_{n,2} \rangle (t) \Rightarrow 0 \qinq \RR \qasq n \ra \infty
\eeqno
for all $t \ge 0$.  We then have
\beq
 \langle M_{n,i}, M_{n,i} \rangle (t) \Rightarrow c_{i,i} t = \mu t \qinq \RR \qasq n \ra \infty
\eeqno
for each $t$ and $i = 1, 2$ by \eqn{xb110} in Theorem \ref{thMartRep2} and Lemma \ref{lemFluid1},
as in the previous argument used in the first proof of Theorem \ref{th1} in
\S\S \ref{secCont}-\ref{secSBcomp}.  As stated above, the bulk
of the proof is thus identical.
By additional argument, we can also show that
\beq
[M_{n,i}, M_{n,i}] (t) \Rightarrow c_{i,i} t = \mu t \qinq \RR \qasq n \ra \infty ~.
\eeqno
starting from \eqn{xb666}.

We have just shown that \eqn{ek6} holds.  It thus remains to show the other conditions in
Theorem \ref{thMart} (ii) are satisfied.
First, since we have a scaled unit-jump counting process,
condition \eqn{ek5} holds by virtue of the scaling in Theorem \ref{thMartRep2} and \eqn{b107}.
Next \eqn{ek4} holds trivially because the predictable quadratic variation processes
$\langle M_{n,1} \rangle$ and $\langle M_{n,2} \rangle$ are
continuous.  Hence this third proof is complete.~~~\bsq

In closing this section, we observe that this alternate method of proof also applies to
the Erlang-$A$ model in \S \ref{secErlangA} and the generalization with finite waiting rooms in Theorem \ref{th3}.

\section*{Acknowledgments}  The authors thank Jim Dai, Itay Gurvich
  and anonymous referees for helpful comments on this paper.
 This work was supported by
NSF Grant DMI-0457095.



\appendix{}

\begin{center}
\Large{\textbf{APPENDIX}}
\end{center}

\section{\texorpdfstring{Proof of Lemma \ref{lemInc}}{Proof of Lemma 3.1}}\label{secDetails}
\hsp
 Let $\Delta$ be the \textbf{jump function}, i.e., $\Delta X (t) \equiv X (t) - X(t-)$ and
let
\beq
\sum_{s \le t} (\Delta X (s))^2, \quad t \ge 0~,
\eeqno
represent the sum of all squared jumps.  Since $(\Delta X (s))^2 \ge 0$, the sum
is independent of the order.  Hence the sum is necessarily well defined, although possibly infinite.
(For any positive integer $n$, a function in $D$
has at most finitely many jumps with $|x(t) - x(t-)| > 1/n$ over any bounded interval; see Lemma 1 on p. 110 of Billingsley \cite{B68}.
Hence, the sample paths
of each stochastic process in $D$ have
only countably many jumps.)

The optional quadratic variation has a very simple structure for locally-bounded-variation processes, i.e.,
processes with sample paths of bounded variation over every bounded interval;
see (18.1) on p. 27 of Rogers and Williams \cite{RW87}.
\begin{lemma}\label{lemOQV}{\em $($optional quadratic covariation and variation for a locally-bounded-variation process$)$}
If stochastic processes $X$ and $Y$ almost surely have sample paths of bounded variation over each bounded interval, then
\beql{oqv1a}
[X, Y] (t) = \sum_{s \le t} (\Delta X (s)) (\Delta Y (s)), \quad t \ge 0~,
\eeq
and
\beql{oqv1b}
[X] (t) \equiv [X,X] (t) = \sum_{s \le t} (\Delta X (s))^2, \quad t \ge 0 ~.
\eeq
\end{lemma}

\begin{lemma}\label{lemQVcount}{\em $($optional quadratic variation for a counting process$)$}
If $N$ is a non-explosive unit-jump counting process with compensator $A$, both adapted to the filtration $\textbf{F}$,
 then $N - A$ is a locally square-integrable martingale
of locally bounded variation with square-bracket process
\begin{eqnarray}\label{qv1}
[N - A] (t) & = & \sum_{s \le t} (\Delta (N-A) (s))^2  \\
& = & N(t) - 2 \int_{0}^{t} \Delta A(s) \, d N(s) + \int_{0}^{t} \Delta A(s) \, d A(s), \quad t \ge 0 ~. \nonumber
\end{eqnarray}
If, in addition, the compensator $A$ is continuous, then
\beql{qv3}
[N - A] =  N ~,
\eeq
\end{lemma}

\paragraph{Proof.}  For \eqn{qv1}, we apply \eqn{oqv1b} in Lemma \ref{lemOQV}; see \S 5.8 of van der Vaart \cite{V06}.  For \eqn{qv3},
we observe that the last two terms in \eqn{qv1} become $0$ when $A$ is continuous.~~~\bsq

\paragraph{Proof of Lemma \ref{lemInc}.}  We will
directly construct the compensator of the submartingale $\{M(t)^2: t \ge 0\}$.
Since (i) $N$ is a non-explosive counting process, (ii)  $E[N(t)] < \infty$ and $E[A(t)] < \infty$ for all $t$
and (iii) $N$ and $A$ have nondecreasing sample paths,
the sample paths of the martingale $M \equiv N - A$ provided
by Theorem \ref{thDoobMeyer} are of bounded variation over bounded intervals.
Since $N(0) = 0$, $M(0) = 0$.
We can thus start by applying integration by parts, referred to as the
product formula on p. 336 of Br\'{e}maud \cite{B81}, to write
\begin{eqnarray}\label{ito1}
M (t)^2 & = &   \int_0^t M (s-) \, d M (s) + \int_0^t M (s) \, d M (s) \nonumber \\
        & = &  2 \int_0^t M (s-) \, d M (s) + \int_0^t (M(s) - M (s-)) \, d M (s)  \nonumber \\
        & = &  2 \int_0^t M (s-) \, d M (s) + [M ](t)  \nonumber \\
        & = &  2 \int_0^t M (s-) \, d M (s) + N (t), \quad t \ge 0 ~,
                \end{eqnarray}
The last step follows from \eqn{qv3} in Lemma \ref{lemQVcount}.

We now want to show that the stochastic integral $\int_0^t M (s-) \, d M (s)$ is a martingale.
To get the desired preservation of the martingale structure, we can apply the integration theorem, Theorem T6 on p. 10 of \cite{B81},
but we need additional regularity conditions.
At this point, we localize in order to obtain boundedness.

To obtain such bounded martingales associated with $N$ and $M \equiv N - A$, let
the stopping times be defined in the obvious way by
\beql{stop2}
\tau_m \equiv \inf\{t \ge 0: |M (t)| \ge m \quad \mbox{or}
\quad A (t) \ge m\} ~.
\eeq
Then $\{\tau_m \le t\} \in \mathcal{F}_t$, $t \ge 0$.
We now define the associated stopped processes:
Let
\beql{stop3}
M^m (t) \equiv M (t\wedge \tau_m), \quad N^m (t) \equiv  N(t\wedge\tau_m) \qandq
A^m (t)\equiv A (t\wedge\tau_m)
\eeq
for all $t \ge 0$ and $m \ge 1$.
 Then $N^m ( t) = M^m (t)+ A^m (t)$ for $t \ge 0$ and $\{M^m (t): t \ge 0 \}$ is a martingale with respect to $\{\sF_t\}$
 having compensator $\{A^m (t): t \ge 0 \}$
 for each $m \ge 1$, as claimed.  Moreover, all three stochastic processes $N^m$, $M^m$ and $A^m$ are bounded.
 The boundedness follows since $N$ has unit jumps and $A$ is continuous.

We then obtain the representation for these stopped processes corresponding to \eqn{ito1};
in particular,
\beql{ito4}
M^m (t)^2 =2 \int_0^t M^m (s-) \, d M^m (s) + N^m (t), \quad t \ge 0 ~,
\eeq
for each $m \ge 1$.
With the extra boundedness provided by the stopping times,
we can apply Theorem T6 on p. 10 of \cite{B81}
to deduce that the integral\break $\int_0^t M^m(s-) dM^m(s)$ is an $\textbf{F}$-martingale.
First, $M^m(t)$ is a martingale of integrable bounded variation with respect to $\textbf{F}$, as defined on p. 10 of \cite{B81}.
Moreover, the process $\{M^m(t-): t \ge 0\}$ is an $\textbf{F}$-predictable process such that
\beq
\int_{0}^{t} |M^m (s-)| \, d | M^m| (s) < \infty ~.
\eeqno
By Theorem T6 on p. 10 of \cite{B81},
the integral $\int_0^t M^m (s-)dM^m (s)$ is an $\textbf{F}$-martingale.
Thus $\{M^m (t)^2 - N^m (t): t \ge 0\}$ is an $\textbf{F}$-martingale.
But $\{N^m (t) - A^m (t): t \ge 0\}$ is also an $\textbf{F}$-martingale.
Adding, we see that
$\{M^m(t)^2- A^m (t): t \ge 0\}$ is an $\textbf{F}$-martingale for each $m$.
Thus, for each $m$, the
predictable quadratic variation of $M^m$ is
$\langle M^m \rangle (t) = A^m (t)$, $t \ge 0$.

Now we can let $m \uparrow \infty$ and
apply Fatou's Lemma to get
\begin{eqnarray*}
E[M (t)^2] &=&
E\Big[\lim\limits_{m\rightarrow\infty}M^m (t)^2\Big]\\
&\leq& \liminf_{m\rightarrow\infty}E \Big[M^m (t)^2\Big]
= \liminf_{m\rightarrow\infty}E \Big[A^m (t)\Big]
=  E\Big[A (t)\Big] < \infty.
\end{eqnarray*}
Therefore, $M$ itself is square integrable.
We can now apply the monotone convergence theorem in the conditioning framework,
as on p. 280 of \cite{B81}, to get
\beq
E[M^m (t+s)^2|\sF_t] \ra E[M (t+s)^2|\sF_t] \qandq E[A^m (t+s)|\sF_t] \ra E[A (t+s)|\sF_t]
\eeqno
as $m \ra \infty$ for each $t \ge 0$ and $s > 0$.  Then, since
\beq
E[M^m (t+s)^2 - A^m (t+s)|\sF_t] = M^m (t)^2 - A^m (t) \qforallq m \ge 1 ~,
\eeqno
$M^m (t) \ra M (t)$ and $A^m (t) \ra A (t)$, we have
\beq
E[M (t+s)^2 - A (t+s)|\sF_t] = M (t)^2 - A (t)
\eeqno
as well, so that $M^2 - A$ is indeed a martingale.
Of course that implies that $\langle M \rangle = A$, as claimed.
We get $[M] = N$ from Lemma \ref{lemQVcount}, as noted at the beginning of the proof.~~~\bsq

We remark that there is a parallel to Lemma \ref{lemQVcount} for the angle-bracket process,
applying to cases in which the compensator is not continuous.
In contrast to Lemma \ref{lemInc}, we now do not assume that $E[N(t)] < \infty$, so we need to localize.
\begin{lemma}\label{lemPQVcount}{\em $($predictable quadratic variation for a counting process$)$}
If $N$ is a non-explosive unit-jump counting process with compensator $A$, both adapted to the filtration $\textbf{F}$,
 then $N - A$ is a locally square-integrable martingale
of locally bounded variation with angle-bracket process
\begin{eqnarray}\label{qv21}
\langle N - A \rangle (t) &  = & \langle [N - A] \rangle (t)  \\
&  = & A(t) - \int_{0}^{t} \Delta A(s) \, d A(s) = \int_{0}^{t} (1 -\Delta A(s)) \, d A(s), \quad t \ge 0 \nonumber ~.
\end{eqnarray}
If, in addition, the compensator $A$ is continuous, then
\beql{qv22}
\langle N - A \rangle = A ~.
\eeq
\end{lemma}

\paragraph{Proof.}
For \eqn{qv21},
we exploit the fact that $\langle N - A \rangle    =  \langle [N - A] \rangle$; see p. 377
of Rogers and Williams \cite{RW87} and \S 5.8 of van der Vaart \cite{V06}.
The third term on the right in \eqn{qv1} is predictable and thus its own compensator.
The compensators of the first two terms in \eqn{qv1} are obtained by replacing $N$ by its compensator $A$.
See Problem 3 on p. 60 of Liptser and Shiryayev \cite{LS89}.~~~\bsq


\begin{thebibliography}{99}



\bibitem{A05}
Armony, M.  2005.
Dynamic routing in large-scale service systems with heterogeneous servers.
\emph{Queueing Systems} 51, 287--329.
\MR{2189596}

\bibitem{AGM06}
Armony M., I., Gurvich and C. Maglaras. 2006. Cross-selling in a call
center with a heterogeneous customer population. Working Paper,
New York University, New York, NY, and Columbia University, New
York, NY.


\bibitem{AGM06a} Armony M., I.  Gurvich and A. Mandelbaum. 2006.
Service level differentiation in call Centers with fully flexible
servers. \emph{Management Science}, forthcoming.

\bibitem{Atar05a} Atar R. 2005.
Scheduling control for queueing systems
with many servers: asymptotic optimality in heavy traffic.
\emph{Ann. Appl. Probab.} 15, 2606--2650.
\MR{2187306}



\bibitem{AMR04a} Atar R., A. Mandelbaum and M. I. Reiman. 2004a.
Scheduling a multi-class queue with many exponential servers.
\emph{Ann. Appl. Probab.} 14, 1084--1134.
\MR{2071417}

\bibitem{AMR04b} Atar R., A. Mandelbaum and M. I. Reiman. 2004b.
Brownian control problems for queueing systems in
the Halfin-Whitt regime.
\emph{Ann. Appl. Probab.} 14, 1084--1134.
\MR{2071417}


\bibitem{BW71}
Bickel P.J. and M.J. Wichura. 1971. Convergence criteria for multiparameter stochastic processes
and some applications. \emph{Ann. Math. Statist.} 42, 1656--1670.
\MR{0383482}

\bibitem{B68}
Billingsley, P. 1968.
\emph{Convergence of Probability Measures}, Wiley
(second edition, 1999).
\MR{1700749}



\bibitem{B67}
Borovkov, A. A. 1967.
On limit laws for service processes in multi-channel systems (in Russian).
\emph{Siberian Math J.} 8, 746--763.
\MR{0222973}

\bibitem{B84}
Borovkov, A. A. 1984.
\emph{Asymptotic Methods in Queueing Theory},
Wiley, New York.
\MR{0745620}

\bibitem{BMR04}
Borst, S., A. Mandelbaum and M. I. Reiman.  2004.
Dimensioning large call centers.
\emph{Oper. Res.} 52, 17--34.
\MR{2066238}

\bibitem{B81}
Br\'{e}maud, P. 1981.
\emph{Point Processes and Queues:  Martingale Dynamics}, Springer.
\MR{0636252}

\bibitem{CR81}
Cs\"{o}rg\'{o} M. and P. R\'{e}v\'{e}z. 1981. \emph{Strong Approximations in Probability
and Statistics}. Akademiai Kiado.


\bibitem{DT05}
Dai J. G. and T. Tezcan. 2005.
State space collapse in many server
diffusion limits of parallel server systems. Working Paper,
Georgia Institute of Technology, Atlanta, GA.

\bibitem{DT06}
Dai J. G. and T. Tezcan. 2006.
Dynamic control of $N$ systems
with many servers:  asymptotic optimality of a static priority policy
in heavy traffic.  Working Paper,
Georgia Institute of Technology, Atlanta, GA.

\bibitem{DT07}
Dai J. G. and T. Tezcan. 2007.
Optimal control of parallel server systems
with many servers
in heavy traffic. Working Paper,
Georgia Institute of Technology, Atlanta, GA.

\bibitem{DVJ03}
Daley, D. J. and D. Vere-Jones.  2003.
\emph{An Introduction to the Theory of Point Processes},
second ed., Springer.
\MR{0950166}

\bibitem{DI91}
Dupuis, P. and H. Ishii.  1991.
On when the solution to the Skorohod problem is Lipschitz continuous with applications.
\emph{Stochastics} 35, 31--62.
\MR{1110990}

\bibitem{EK86}
Ethier, S. N. and T. G. Kurtz.  1986.
\emph{Markov Processes;  Characterization and Convergence}, Wiley.
\MR{0838085}


\bibitem{Gans03}
Gans, N., G. Koole and A. Mandelbaum. 2003.
Telephone call
centers: tutorial, review and research prospects.  \emph{Manufacturing Service Oper. Management} \textbf{5}(2), 79--141.

\bibitem{GMR02}
Garnett, O., A. Mandelbaum and M. I. Reiman.  2002.
Designing a call center with impatient customers.
\emph{Manufacturing Service Oper. Management} 4, 208--227.

\bibitem{GW91}
Glynn, P. W. and W. Whitt.  1991.
A new view of the heavy-traffic limit theorem for the infinite-server queue.
\emph{Adv. Appl. Prob.} 23, 188--209.
\MR{1091098}

\bibitem{GW07a}
Gurvich, I. and W. Whitt.  2007a.
Queue-and-idleness-ratio controls in many-server servbice systems.
working paper, Columbia University.
Available at:  \url{http://www.columbia.edu/~ww2040}

\bibitem{GW07b}
Gurvich, I. and W. Whitt.  2007b.
Service-level differentiation in many-server service systems: a solution based on fixed-queue-ratio routing.
working paper, Columbia University.
Available at:  \url{http://www.columbia.edu/~ww2040}

\bibitem{GW07c}
Gurvich, I. and W. Whitt.  2007c.
Scheduling Flexible Servers with convex delay costs in many-server service systems.
\emph{Manufacturing and Service Operations Management}, forthcoming.
Available at:  \url{http://www.columbia.edu/~ww2040}

\bibitem{HW81}
Halfin, S. and W. Whitt. 1981.
Heavy-traffic limits for queues with many exponential servers.
\emph{Oper. Res.} 29, 567--588.
\MR{0629195}

\bibitem{HZ04}
Harrison, J. M. and A. Zeevi.  2005.
Dynamic scheduling of a multiclass queue in the Halfin and Whitt heavy traffic regime.
\emph{Oper. Res.} 52, 243--257.
\MR{2066399}

\bibitem{H74}
Helms, L. L. 1974.
Ergodic properties of several interacting Poisson particles.
\emph{Advances in Math.} 12, 32--57.
\MR{0345247}

\bibitem{I65}
Iglehart, D. L. 1965.
Limit diffusion approximations for the many server queue and
the repairman problem.  \emph{J. Appl. Prob.}  2, 429--441.
\MR{0184302}


\bibitem{JS87}
Jacod, J. and A. N. Shiryayev.  1987.
\emph{Limit Theorems for Stochastic Processes}, Springer.
\MR{0959133}

\bibitem{JMM04}
Jelenkovi\'{c}, P., A. Mandelbaum and P. Mom\v{c}ilovi\'{c}.  2004.
Heavy traffic limits for queues with many deterministic servers.
\emph{Queueing Systems} 47, 53--69.
\MR{2074672}

\bibitem{K02}
Kallenberg, O.  2002.
\emph{Foundations of Modern Probability}, second edition,
Springer.
\MR{1876169}

\bibitem{KS88}
Karatzas, I, and S. Shreve. 1988.
\emph{Brownian Motion and Stochastic Calculus}, Springer.
\MR{0917065}

\bibitem{KR07}
Kaspi, H. and K. Ramanan. 2007.
Law of large numbers limit for many-server queues.
Working paper.  The Technion and Carnegie Mellon University.


\bibitem{K02}
Khoshnevisan D. 2002. \emph{Multiparameter Processes: An Introduction to
Random Fields}, Springer.
\MR{1914748}

\bibitem{KLS86}
Kogan, Y., R. Sh. Liptser and A. V. Smorodinskii.  1986.
Gaussian diffusion approximations of closed Markov models
of computer networks.
\emph{Problems. Inform. Transmission} 22, 38--51.
\MR{0838688}

\bibitem{KP97}
Krichagina, E. V. and A. A. Puhalskii.  1997.
A heavy-traffic analysis of a closed queueing system with a $GI/\infty$ service center.
\emph{Queueing Systems} 25, 235--280.
\MR{1458591}

\bibitem{KW67}
Kunita, H. and S. Watanabe.  1967.
On square-integrble martingales.
\emph{Nagoya Math. J.} 30, 209--245.
\MR{0217856}


\bibitem{K78}
Kurtz, T.  1978.
Strong approximation theorems for density dependent Markov chains.
\emph{Stoch. Process Appl.}  6, 223--240.
\MR{0464414}

\bibitem{K80}
Kurtz, T.  1980.
Representations of Markov processes as multiparameter time changes.
\emph{Ann. Probability} 8, 682--715.
\MR{0577310}



\bibitem{K01}
Kurtz, T.  2001.
\emph{Lectures on Stochastic Analysis}, Department of Mathematics and Statistics,
University of Wisconsin, Madison, WI 53706-1388.

\bibitem{L48}
L\'{e}vy, P. 1948.
\emph{Processus Stochastiques et Mouvement Borownien},
Gauthiers-Villars, Paris.


\bibitem{LS84}
Lions, P. and A. Sznitman.  1984.
Stochastic differential equations with reflecting boundary conditions.
\emph{Commun. Pure Appl. Math.}  37, 511--537.
\MR{0745330}

\bibitem{LS89}
Liptser, R. Sh. and A. N. Shiryayev.  1989.
\emph{Theory of Martingales},
Kluwer (English translation of 1986 Russian edition).

\bibitem{L88}
Louchard, G. 1988. Large finite population queueing systems. Part
I: The infinite server model. \emph{Comm. Statist. Stochastic
Models} 4(3), 473--505.
\MR{0971602}


\bibitem{MMR98}
Mandelbaum, A., W. A. Massey and M. I. Reiman.  1998.
Strong approximations for Markovian service networks.
\emph{Queueing Systems} 30, 149--201.
\MR{1663767}

\bibitem{MP95}
Mandelbaum, A. and G. Pats.  1995.
State-dependent queues:  approximations and applications.
In \emph{Stochastic Networks},
F. P. Kelly, R. J. Williams (eds.),
Institute for Mathematics and its Applications, Vol. 71, Springer,
239--282.
\MR{1381015}

\bibitem{MP98}
Mandelbaum, A. and G. Pats.  1998.
State-dependent stochastic networks. Part I:  Approximations and Applications
with continuous diffusion limits.
\emph{Ann. Appl. Prob.} 8, 569--646.
\MR{1624965}

\bibitem{MZ05}
Mandelbaum, A. and S. Zeltyn.  2005.
The Erlang-$A$/Palm queue, with applications to call centers.
Working paper, The Technion, Haifa, Israel.
Available at:
\href{http://iew3.technion.ac.il/serveng/References/references.html}{\texttt{http://iew3.technion.ac.il/serveng/}}\allowbreak\href{http://iew3.technion.ac.il/serveng/References/references.html}{\texttt{References/references.html}}

\bibitem{MW93}
Massey, W. A. and W. Whitt.  1993.
Networks of infinite-server queues with nonstationary Poisson input.
\emph{Queueing Systems} 13, 183--250.
\MR{1218848}



\bibitem{P67}
Parthasarathy, K. R. 1967.
\emph{Probability Measures on Metric Spaces}, Academic Press.
\MR{0226684}


\bibitem{P94}
Puhalskii, A. A.  1994.
On the invariance principle for the first passage time.
\emph{Math. Oper. Res.} 19, 946--954.
\MR{1304631}

\bibitem{PR00}
Puhalskii, A. A. and M. I. Reiman.  2000.
The mutliclass $GI/PH/N$ queue in the Halfin-Whitt regime.
\emph{Adv. Appl. Prob.} 32, 564-595.
\MR{1778580}

\bibitem{R80}
Rebolledo, R.  1980.
Central limit theorems for local martingales.
\emph{Zeitschrift Wahrscheinlichkeitstheorie verw. Gebiete} 51, 269--286.
\MR{0566321}

\bibitem{R05}
Reed, J. 2005.
The $G/GI/N$ queue in the Halfin-Whitt regime.
working paper, Georgia Institute of Technology.

\bibitem{RW07}
Reed, J. and A. R. Ward. 2007.
Approximating the $GI/GI/1+GI$ queue with a nonlinear drift diffusion:
hazard-rate scaling in heavy traffic.
working paper, Georgia Institute of Technology.

\bibitem{R03}
Robert, P. 2003. \emph{Stochastic Networks and Queues}, Springer.
\MR{1996883}



\bibitem{RW87}
Rogers, L. C. G. and D. Williams.  1987.
\emph{Diffusions, Markov Processes and Martingales, Volume 2:  Ito Calculus},
Wiley.
\MR{0921238}

\bibitem{RW00}
Rogers, L. C. G. and D. Williams.  2000.
\emph{Diffusions, Markov Processes and Martingales, Volume 1:  Foundations},
Cambridge University Press.
\MR{1796539}

\bibitem{S56}
Skorohod, A. V.  1956.
Limit theorems for stochastic processes.
\emph{Prob. Theory Appl.} 1, 261--290.
\MR{0084897}


\bibitem{SW96}
Srikant, R. and W. Whitt.  1996.
Simulation run lengths to estimate blocking probabilities.
\emph{ACM Trans. Modeling Computer Simulations} 6, 7--52.

\bibitem{S63}
Stone, C. 1963.
Limit theorems for random walks, birth and death processes and diffusion processes.
\emph{Illinois J. Math.} 4, 638--660.
\MR{0158440}





\bibitem{T06}
Tezcan T. 2006.
Optimal control of distributed parallel server
systems under the Halfin and Whitt regime. Working Paper, University of Illinois
at Urbana-Champaign.
Available at:
\href{https://netfiles.uiuc.edu/ttezcan/www/TolgaTezcansubmittedMOR121906.pdf}{\texttt{https://netfiles.uiuc.edu/}}\allowbreak\href{https://netfiles.uiuc.edu/ttezcan/www/TolgaTezcansubmittedMOR121906.pdf}{\texttt{ttezcan/www/TolgaTezcansubmittedMOR121906.pdf}}

\bibitem{V06}
van der Vaart, A. W. 2006.
\emph{Martingales, Diffusions and Financial Mathematics}
Lecture Notes,
Available at:
\href{http://www.math.vu.nl/sto/onderwijs/mdfm/}{\texttt{http://www.math.vu.nl/sto/}}\allowbreak\href{http://www.math.vu.nl/sto/onderwijs/mdfm/}{\texttt{onderwijs/mdfm/}}

\bibitem{WG03a}
Ward, A. R. and P. W. Glynn.  2003a.
A diffusion approximation for a Markovian queue with reneging.
\emph{Queueing Systems} 43, 103--128.
\MR{1957808}

\bibitem{WG03b}
Ward, A. R. and P. W. Glynn.  2003b.
Properties of the reflected Ornstein-Uhlenbeck process.
\emph{Queueing Systems} 44, 109--123.
\MR{1993278}

\bibitem{WG05}
Ward, A. R. and P. W. Glynn.  2005.
A diffusion approximation for a $GI/GI/1$ queue with balking or reneging.
\emph{Queueing Systems} 50, 371--400.
\MR{2172907}



\bibitem{W82}
Whitt, W. 1982.
On the heavy-traffic limit theorem for $GI/G/\infty$ queues.
\emph{Adv. Appl. Prob.} 14, 171--190.
\MR{0644013}


\bibitem{W02}
Whitt, W. 2002.
\emph{Stochastic-Process Limits}, Springer.
\MR{1876437}

\bibitem{W02IS}
Whitt, W. 2002a.
\emph{Internet Supplement to Stochastic-Process Limits},
Available at:
\url{http://www.columbia.edu/~ww2040/supplement.html}


\bibitem{W04}
Whitt, W. 2004.
Efficiency-driven heavy-traffic approximations for many-server queues with abandonments.
\emph{Management Science} 50, 1449--1461.

\bibitem{W05}
Whitt, W. 2005.
Heavy-traffic limits for the $G/H^{*}_{2}/n/m$ queue.
\emph{Math. Oper. Res.} 30, 1--27.
\MR{2125135}

\bibitem{W07}
Whitt, W. 2007.
Proofs of the martingale functional central limit theorem.
\emph{Probability Surveys}, forthcoming.\textsf{}


\bibitem{ZM05}
Zeltyn S. and A. Mandelbaum.  2005.
  Call centers with impatient customers: many-server asymptotics of the M/M/n+G queue.
 \emph{Queueing Systems} 51, 361--402.
\MR{2189598}

\end{thebibliography}
\end{document}